\documentclass[a4paper,10pt]{article}
\usepackage[utf8x]{inputenc}
\usepackage{tracefnt,amsmath,tabu,array}
\usepackage{amssymb,graphicx,setspace,amsfonts,amsbsy}
\usepackage{pifont,latexsym,ifthen,amsthm,rotating,calc,textcase,booktabs,cancel,slashed}

\newtheorem{theorem}{Theorem}[section]
\newtheorem{lemma}[theorem]{Lemma}
\newtheorem{corollary}[theorem]{Corollary}
\newtheorem{proposition}[theorem]{Proposition}
\newtheorem{definition}[theorem]{Definition}

\newtheorem{remark}[theorem]{Remark}
\newcommand{\filledbox}{\leavevmode
  \hbox to.77778em{%
  \hfil\vbox to.675em{\hrule width.6em height.6em}\hfil}}

\newcommand{\Rm}{{\mathbb R}}

\begin{document}
\tabulinesep=1.0mm
\title{Modified wave operators and scattering for linear wave equations with a repulsive potential}

\author{Boya Fan and Ruipeng Shen\\
Centre for Applied Mathematics\\
Tianjin University\\
Tianjin, China
}

\maketitle

\begin{abstract}
 In this work we consider the wave equation with a repulsive potential, either on the half line $\Rm^+$ or the Euclidean space $\Rm^d$ with $d\geq 3$. We combine the operator theory and the inward/outward energy theory to deduce a modified wave operator for repulsive potentials decaying like $|x|^{-\beta}$ with $\beta>1/3$. In particular the regular wave operator without modification exists if $\beta>1$. This implies that the asymptotic behaviour of finite-energy solutions to the wave equation $u_{tt} - \Delta u + |x|^{-\beta} u =0$ is similar to that of the solutions to the classic wave equation if $\beta \in (1,2)$. 
\end{abstract}

\section{Introduction} 

\subsection{Topic and background}

The topic of this work is asymptotic behaviour of solutions to the wave equation with a radially symmetric repulsive potential 
\[ 
 \partial_t^2 u - \Delta u + V(x) u = 0, \qquad (x,t) \in \Rm^d \times \Rm. 
\]
Here ``repulsive'' means that the potential is positive, decreases as the radius grows and approaches zero as radius tends to infinity. This covers many important potentials in physics. For example, this includes the distractive Coulomb potential $V(x)=1/|x|$, which originates from a quantum mechanical description of the distractive Coulomb force between two particles with the same charge. There are typically two aspects to be discussed about the global behaviour of solutions to the wave equations given above. 
\begin{itemize}
 \item How to describe asymptotic behaviour of solutions as the time tends to infinity. This is usually done by the introduction of the conceptions of scattering and (modified) wave operators. 
 \item How to describe the decay of solutions as the time tends to infinity. This usually involves the dispersive and Strichartz estimates. 
\end{itemize}
In this work we focus on the first aspect, i.e. the scattering and modified wave operator; and we are particularly interested in the inverse power potential $V(x)=|x|^{-\beta}$ with $0<\beta<2$. Let us first make a brief review on previously known results related to this topic. 

\paragraph{Schr\"{o}dinger equations} The (modified) wave operators have been extensively studied in the case of Schr\"{o}dinger equations with a potential. Please refer to Beceanu \cite{BSchrWOAb}, Beceanu-Schlag \cite{BSSchrWO}, Christ-Kiselev \cite{CKSchrWO}, D'Ancona-Fanelli \cite{WMLpSchr}, H\"{o}rmander \cite{LHSchrWO}, Journ\'{e}-Soffer-Sogge \cite{JSSdecaySchr}, Reed-Simon \cite{Simon}, Weder \cite{WWkpSchr} and Yajima \cite{YWkpSchr, YLpSchr}, for example. 

\paragraph{Wave equations} There are very few literatures about the wave operators for wave equations with a potential. One previous result is for the wave equation with an inverse square potential 
\[
 \partial_t^2 u - \Delta u + \frac{a}{|x|^2} u = 0, \qquad (x,t) \in \Rm^n \times \Rm. 
\]
Mizutani \cite{WISnonradial} proves that the following wave operator is well-defined in the energy space $\dot{H}^{1} \times L^2$ for all parameters $a>-(n-2)^2/4$:
\[
  \mathop{{s-}\lim}\limits_{t\rightarrow +\infty} \vec{\mathbf{S}}_0 (-t) \vec{\mathbf{S}}_a (t).
\]
Here $\vec{\mathbf{S}}_a(t)$ is the wave propagation operator of the equation above. Namely, if $\vec{u}_0=(u_0,u_1)$ are the initial data and $u(t)$ is the corresponding solution, then we have $\vec{\mathbf{S}}_a(t) \vec{u}_0 = (u(\cdot, t), u_t(\cdot, t))$. In particular, $\vec{\mathbf{S}}_0(t)$ is the wave propagation operator of the classic wave equation $\partial_t^2 u - \Delta u = 0$. The author would like to mention that there are also many works about the dispersive/Strichartz estimates for free waves with a potential $V(x)$ whose decay rate is faster than or roughly equal to $|x|^{-2}$. Please see Beals-Strauss \cite{WPfastsmall}, Costin-Huang  \cite{WPfast1d}, Cuccagna \cite{WPfast3d}, Green \cite{WPfast2d}, Petkov \cite{WPcompact}, Vodev \cite{WPfasthigher} for compact-supported or fast decaying potentials; and D'ancona-Pierfelice \cite{WPKato1}, Pierfelice \cite{WPKato2} and Bui-Duong-Hong \cite{WPKato3} for  potentials $V(x)$ in the Kato class, i.e. 
 \[
  \|V\|_K = \sup_{x\in \Rm^n} \int_{\Rm^n} \frac{|V(y)|}{|x-y|^{n-2}} {\mathrm d} y < +\infty.  
 \]
The decay/dispersive estimates for inverse-square potentials $V(x)=a|x|^{-2}$ and other potentials with a similar decay at the infinity have been discussed by Donninger-Schlag \cite{WIS1d}, Donninger-Krieger \cite{WIS1d2}, Georgiev-Visciglia \cite{WISGV} and Planchon-Stalker-Tahvildar-Zadeh \cite{WISPST1, WISPST2}. The corresponding Stichartz estimates were given by Burq, Planchon, Stalker and Tahvildar-Zadeh in \cite{WISStrichartz1, WISStrichartz2}. Miao-Zhang-Zheng \cite{WISStrichartzRadial} showed that the range of admissible pairs for Strichartz estimates can be extended if the initial data possess additional angular regularity.
\subsection{Main idea}

In this work we let the space dimension $d\geq 3$, consider the linear wave equation with a repulsive potential
\begin{equation} \label{linear wave equation with potential}
 u_{tt} -\Delta u + V(x) u = 0,
\end{equation}
and discuss the asymptotic behaviour of the corresponding finite-energy solutions, especially its modified wave operator.  For convenience we may define 
\[
 \mathbf{H} = - \Delta + V(x). 
\]
We are particularly interested in the case with inverse power potential $V(x) = q(|x|)= |x|^{-\beta}$. Our main idea is to decompose this equation into a family of wave equations (with zero boundary condition) on the half line $\Rm^+$
\[
 w_{tt} - w_{rr} + \mu r^{-2} w + q(r) w = 0,
\]
and then discuss the modified wave operators for these wave equations, by making use of spectral theory of the corresponding ordinary differential self-adjoint operators.  The author would like to mention that Denisov \cite{WPDen} gives the existence of modified wave operator in $L^2(\Rm^+)$
\begin{equation} \label{strong limit Den}
 \mathop{{s-}\lim}\limits_{t\rightarrow +\infty} \exp (-{\mathrm i} t \sqrt{\mathbf{A}}) \mathbf{U}(t) 
\end{equation}
for the half-wave operator of the positive part $\mathbf{A}$ of the self-adjoint operators $-{\rm d}^2/{\rm d} x^2 + q(x)$ with $q \in L^2(\Rm^+)$ and $|q(x)|\lesssim C(1+x)^{-1/2}$. Here $\mathbf{U}(t)$ is the operator with a Fourier multiplier (in the classic Fourier analysis of half line, details can be found below)
\begin{equation} \label{multiplier beta bigger}
 \exp \left[\mathrm i \left(kt + \frac{1}{2k}\int_0^t q(s) {\rm d} s\right)\right].
\end{equation} 
The asymptotic completeness, i.e. the surjectivity of the modified wave operator given by \eqref{strong limit Den}, does not necessarily hold, due to the possible presence of point and singular continuous spectrum. In this work we focus on repulsive potentials. The corresponding self-adjoint operators are positive and always come with purely absolutely continuous spectrum, thus one may expect that the wave operator is actually a bijection. In order to study the asymptotic behaviour of energy solutions to the wave equation with a potential, it is more convenient to consider the strong limit of $\mathbf{U}^{-1} (t) \exp ({\mathrm i} t \sqrt{\mathbf{A}})$, because $\exp ({\mathrm i} t \sqrt{\mathbf{A}})$ is unitary in $L^2(\Rm^+)$ but is not unitary in $\dot{H}^1\times L^2$. Nevertheless, there are still a few challenges in the argument: 
\begin{itemize}
 \item Denisov's method highly depends on the calculation in Killip and Simon's remarkable work \cite{KS}, where the $L^2$ assumption on the potential $q(r)$ plays an essential role. However, we have to consider the self-adjoint operators $-{\rm d}^2/{\rm d}r^2 + \mu r^{-2} + q(r)$, which come with a strong singularity at zero when $\mu > 0$. Even if $\mu = 0$, a typical potential, for example the Coulomb potential $1/r$, is not contained in the $L^2$ space in a neighbourhood of zero, although it decays sufficiently fast near the infinity to fit in the $L^2$ space. 
 \item We would like to consider repulsive potentials  $q(x) \simeq x^{-\beta}$ with a lower decay rate $\beta < 1/2$ as well. Clearly this implies $q\notin L^2$ near the infinity. We would like to investigate whether a similar approximated wave propagation operator $\mathbf{U}(t)$ can be given in by a regular Fourier multiplier, and to determine the Fourier multipliers if they are different form those given above, at least for repulsive potential with a decay rate $\beta$ slightly smaller than $1/2$. 
 \item We also need to give a suitable approximated wave propagation operator for the wave propagation with a radial potential $V(x)=q(|x|)$ in $\Rm^d$. This means to combine infinitely many solutions to (possibly slightly different) wave equations on the half line, which is not as simple as at the first glance.  
\end{itemize}

\paragraph{The main ingredients} To overcome the difficulties mentioned above, our argument contains the following main ingredients 
\begin{itemize}
 \item Inward/outward energy theory describes the general rule of energy distribution, propagation and conversion of solutions to the wave equations on the half line with a repulsive potential. As an application we may show that the asymptotic behaviour of solutions is independent of the values of potential in a compact region as long as the potential is repulsive. This helps us remove the possible singularity of potential near zero. 
 \item We introduce an expansion formula of wave functions of the self-adjoint operators $\mathbf{A} = -{\rm d}^2/{\rm d} x^2 + q(x)$ for a suitable potential $q \in \mathcal{C}^1([0,+\infty))$ which decays like $x^{-\beta}$ near the infinity, then utilize it to introduce a modified wave operator $\mathbf{U}(t)$ with Fourier multiplier 
 \[
  \exp \left[{\mathrm i} \left(kt + \frac{1}{2k} \int_1^t q(s) {\rm d} s - \frac{1}{4k^3} q(t) \int_1^t q(s) {\rm d} s + \frac{1}{8k^3} \int_1^t q^2(s) {\rm d} s\right)\right],
 \]
for the half-wave operator $\exp ({\mathrm i} t \sqrt{\mathbf{A}})$ if $\beta>1/3$, and finally deduce the corresponding wave operator for the wave propagation operator of this equation. Please note that if $q(x)$ decays sufficiently fast near infinity, we may ignore the last two terms in the multiplier above, thus use the same multiplier as in \eqref{multiplier beta bigger}.
\item As for the higher dimensional case, we show that the Fourier transform of each term in the spherically harmonic decomposition, thus the Fourier transform $\hat{u}(\xi, t)$ of the free wave itself is an approximated solution to the ordinary differential equation (with parameter $|\xi|$)
\[
 v_{tt} + \left(|\xi|^2 + q(t) - \frac{q'(t) \int_1^t q(s) {\rm d}s}{2|\xi|^2 + 2q^2(t)}\right) v = 0;
\]
and then apply the inverse Fourier transform and a result on the asymptotic behaviour ($t\rightarrow +\infty$) of approximated solutions to the ordinary differential equation  above to deduce the modified wave operator in the higher dimensional case. Again if $q$ decays fast with a decay rate $\beta > 1/2$, we may substitute the equation above by $v_{tt} + (|\xi|^2+q(t)) v =0$, thus the solution itself is an approximated solution to the wave equation with a time-dependent potential
\[
 u_{tt} - \Delta u + q(t) u = 0.
\]
\end{itemize}

\subsection{Assumptions and notations} \label{sec: notation}

Before we give the main results, we first make a few definitions for potentials and introduce a few notations. We start from potentials defined on the half-line $\Rm^+$.
\begin{definition}
 We call a potential $q(x)$ repulsive if and only if
 \begin{itemize}
  \item $q \in \mathcal{C}^1(0,+\infty)$;
  \item $q(x)>0$ and $q'(x)<0$;
  \item $q(x)$ converges to zero as $x$ tends to infinity.
 \end{itemize}
 \end{definition} 
 \begin{definition} 
  We say $q(x)$ satisfies the growth condition at zero if there exists a constant $0<\kappa<2$, such that $|q(x)|\lesssim x^{-\kappa}$ and $|q'(x)|\lesssim x^{-\kappa-1}$ for all $x$ near zero. 
 \end{definition} 
 \begin{definition}
  We say $q(x)$ satisfies the decay condition if there exists a constant $\beta>0$, such that $|q(x)| \lesssim x^{-\beta}$ for large $x\gg 1$. If $\beta\leq 1/2$, we also assume that $|q'(x)| \lesssim x^{-\beta-1}$ for large $x\gg 1$. The constant $\beta$ is called a decay rate of $q(x)$. 
 \end{definition}
 We say the a radial potential $V(x) = q(|x|)$ in $\Rm^d$ is repulsive or satisfies the growth/decay condition if $q(x)$ satisfies the corresponding assumption above.

\begin{definition} 
 In this work we mainly consider three types of repulsive potentials: Besides the repulsive condition we also assume: 
 \begin{itemize}
  \item Type I: $q(x)\in \mathcal{C}^1 ([0,+\infty))$ satisfies the decay condition with a decay rate $\beta>0$; 
  \item Type II: $q(x)$ satisfies growth condition at zero and the decay condition with a decay rate $\beta>1/3$; 
  \item Type III: $q(x) = \mu x^{-2} + q_0(x)$, where $q_0(x)$ is a type II repulsive potential and $\mu \geq 3/4$. 
 \end{itemize}
\end{definition}

\begin{remark}
 The physical meaning of the repulsive assumption is clear. Let us consider the simple case of two particles, for instance. If $q'(x)<0$, then the force between these particles push them away from each other. Please note that a type I repulsive potential with a decay rate $\beta>1/3$ must be a type II repulsive potential as well. 
\end{remark}

\paragraph{Self-adjoint operators} Let $q_0(x)$ be a repulsive potential satisfying the growth condition at zero and $\lambda \in \{0\}\cup [3/4,+\infty)$ be a real number. Then the operator $\mathbf{A} = -{\rm d}^2/{\rm d}x^2 + q_0(x) + \lambda x^{-2}$ with zero boundary condition at zero is a self-adjoint operator. More details about its actual domain can be found at the beginning of Section 2. Let $V(x)$ be a type II radial repulsive potential. The self-adjointness of $-\Delta + V(x)$ is briefly discussed in Section 5. 

\paragraph{Fourier transforms} In this work we use the following notation for the Fourier transform. In the case of $\Rm^d$, $\mathcal{F}_0 f$ and $\hat{f}$ both represent the classic Fourier transform of $f$. As for $\Rm^+$, there are two cases: Case one, the classic Fourier transform $\mathcal{F}_0$ and its inverse on half-line defined by 
\begin{align*}
 &(\mathcal{F}_0 f) (k) = \int_0^\infty (\sin kr) f(r) {\rm d} r, \qquad k>0;& &(\mathcal{F}_0^{-1} g)(r) = \frac{2}{\pi} \int_0^\infty (\sin kr) g(k) {\rm d} k.&
\end{align*}
This can be viewed as a version of the regular Fourier transform for odd function on $\Rm^+$. Case two, the Fourier transform $\mathcal{F}$ associate a self-adjoint operator $\mathbf{A} = -{\rm d}^2 /{\rm d} x^2 + q(x)$ with zero boundary condition, where $q$ is a type I repulsive potential. We first define wave functions $u(x,k)$ to be the solution to the equation 
\[
 -u''(x) + q(x) u(x) = k^2 u(x)
\]
with boundary condition $u(0) = 0$ and $u'(0)=1$.  Then the Fourier transform and its inverse can be defined by 
\begin{align*}
 &(\mathcal{F} f)(k)= \int_0^\infty u(x,k) f(x) {\rm d}x; & & (\mathcal{F}^{-1} g)(x) = \int_0^\infty u(x,k) g(k) {\rm d} \rho (E).
\end{align*}
Here $\rho$ is the spectral measure associated to $\mathbf{A}$. More details will be given at the beginning of Section 3. Please note that if $q(x) \equiv 0$, then we have $u(x,k) = (1/k) \sin kx$, thus $(\mathcal{F}_0 f)(k) = k (\mathcal{F} f)(k)$. 

\paragraph{Linear propagation operators} Given a linear wave equation $\partial_t^2 u + \mathbf{A} u = 0$, where $\mathbf{A} = -{\rm d}^2/{\rm d}x^2 +q(x)$ is a positive self-adjoint operator, we may give its solution with initial data $(u_0,u_1)$ by functional calculus
\[
 \begin{pmatrix} u(\cdot,t)\\ u_t(\cdot,t) \end{pmatrix} = \vec{\mathbf{S}}_q (t) \begin{pmatrix} u_0\\ u_1 \end{pmatrix} = \begin{pmatrix} \cos (t\mathbf{A}^{1/2}) & \mathbf{A}^{-1/2} \sin (t\mathbf{A}^{1/2}) \\ -\mathbf{A}^{1/2} \sin (t \mathbf{A}^{1/2}) & \cos (t \mathbf{A}^{1/2}) \end{pmatrix} \begin{pmatrix} u_0\\ u_1 \end{pmatrix}.
\]

\paragraph{Sobolev spaces} Let $\mathbf{A} = -{\rm d}^2/{\rm d}x^2 +q(x)$ be a positive self-adjoint operator. We define the corresponding Sobolev spaces by 
\begin{align*}
 &\dot{\mathcal{H}}_{\mathbf{A}}^s = \{f: \mathbf{A}^{s/2} f \in L^2(\Rm^+)\}& &\|f\|_{\dot{\mathcal{H}}_{\mathbf{A}}^s} = \|\mathbf{A}^{s/2} f\|_{L^2(\Rm^+)}. 
\end{align*}
Please note that we refrain from using the standard notation $\dot{H}^s$ because $\dot{H}^s$ is also frequently used in this work for the standard Sobolev space defined by $-{\rm d}^2/{\rm d} x^2$ or $-\Delta$ instead. For $s>0$, we also have $\mathcal{H}_{\mathbf{A}}^s = \dot{\mathcal{H}}_{\mathbf{A}}^s \cap L^2(\Rm^+)$. In particular we have for type I, II or III repulsive potentials $q(x)$ that
\[
 \dot{\mathcal{H}}_{\mathbf{A}}^1 = \left\{u \in \dot{H}^1(\Rm^+): \int_0^\infty \left(|u'(x)|^2 + q(x)|u(x)|^2 \right){\rm d} x < +\infty, \; u(0)=0 \right\}.
\]
Similarly if $V$ is a radially symmetric type I, II or III potential in $\Rm^d$, then we have 
\[
 \dot{\mathcal{H}}_V^1(\Rm^d) = \left\{u \in \dot{H}^1(\Rm^d): \int_{\Rm^d} \left(|\nabla u(x)|^2 + V(x)|u(x)|^2 \right){\rm d} x < +\infty \right\}.
\]
\paragraph{The notation $\lesssim$} In this work the notation $A\lesssim B$ means that there exists a constant $c$ so that $A \leq c B$. We may also add subscript(s) to the notation $\lesssim$ to emphasize that the constant $c$ depends on the subscript(s) but nothing else. In particular, the notation $\lesssim_1$ means that the implicit constant $c$ is an absolute constant. The meaning of notations $\simeq$ and $\gtrsim$ is similar.

\subsection{Main results}

Now we give the main results of this work. We start by the half-line case. 

\begin{theorem}[Wave operator on the half-line] \label{thm main one-dimensional}
Assume that $q(x)$ is a type I, II or III repulsive potential with a decay rate $\beta>1/3$. Let $\mathbf{A} = -{\rm d}^2/{\rm d} x^2 + q(x)$ be the self-adjoint operator with zero boundary condition and $\vec{\mathbf{S}}_q$ be its wave propagation operator. We define
\begin{align*} 
 &\vec{\mathbf{U}}(t) = \mathcal{F}_0^{-1} \begin{pmatrix} \cos \eta (k,t) & k^{-1} \sin \eta(k,t) \\ -k \sin \eta(k,t) & \cos \eta(k,t) \end{pmatrix} \mathcal{F}_0; &
 & \eta(k,t) = kt + P(k,t). 
\end{align*}
with phase shift function 
\[
 P(k,t) = \frac{1}{2k} \int_1^t q(s) {\rm d} s - \frac{1}{4k^3} q(t) \int_1^t q(s) {\rm d} s + \frac{1}{8k^3} \int_1^t q^2(s) {\rm d} s.
\]
Then the modified wave operator defined by the strong limit in $\dot{H}^1(\Rm^+)\times L^2(\Rm^+)$
\begin{align*}
 \vec{\mathbf{W}} \doteq \mathop{{s-}\lim}\limits_{t\rightarrow +\infty} \vec{\mathbf{U}}(t)^{-1} \vec{\mathbf{S}}_q (t) 
\end{align*}
exists and is a unitary bijection from the energy space $\dot{\mathcal{H}}_{\mathbf{A}}^1 \times L^2(\Rm^+)$ to $\dot{H}^1(\Rm^+) \times L^2(\Rm^+)$. In addition, we may use simpler phase shift function in $\vec{\mathbf{U}}(t)$ if the decay rate of $q$ is higher. More precisely, the same conclusion holds in the following situations: 
\begin{itemize} 
 \item If the decay rate $\beta > 1/2$, we may substitute $P(k,t)$ by $\displaystyle P_1(k,t) =  \frac{1}{2k} \int_1^t q(s) {\rm d} s$; 
 \item If $q$ also satisfies $q \in L^1(1,+\infty)$, then we may substitute $\vec{\mathbf{U}}(t)$ by the classic wave propagation operator 
 \[
  \vec{\mathbf{S}}_0(t) = \mathcal{F}_0^{-1} \begin{pmatrix} \cos kt & k^{-1} \sin kt \\ -k \sin kt & \cos kt \end{pmatrix} \mathcal{F}_0.
 \]
\end{itemize} 
\end{theorem}


\begin{theorem}[Wave operator in the high-dimensional case] \label{thm main high-dimensional}
Assume that $d\geq 3$ and $V(x) = q(|x|)$ is a Type II repulsive potential with a decay rate $\beta>1/3$. If $\beta \leq 1/2$, we also assume that $q$ is $\mathcal{C}^2$ when $x$ is large and satisfies either $q''(x)>0$ or $|q''(x)| \lesssim x^{-2}$ for large $x$. Let $\vec{\mathbf{S}}_V$ be the corresponding wave propagation operator of the wave equation $u_{tt} - \Delta u + V(x) u = 0$ and define
\begin{align*} 
 &\vec{\mathbf{U}}(t) = \mathcal{F}_0^{-1} \begin{pmatrix} \cos \eta (|\xi|,t) & |\xi|^{-1} \sin \eta(|\xi|,t) \\ -|\xi| \sin \eta(|\xi|,t) & \cos \eta(|\xi|,t) \end{pmatrix} \mathcal{F}_0; &
 & \eta(|\xi|,t) = |\xi|t + P(|\xi|,t). 
\end{align*}
Here the phase shift function $P(|\xi|,t)$ is defined by 
\[
 P(|\xi|,t) = \left\{\begin{array}{ll} \frac{1}{2|\xi|} \int_1^t q(s) {\rm d} s - \frac{1}{4|\xi|^3} q(t) \int_1^t q(s) {\rm d} s + \frac{1}{8|\xi|^3} \int_1^t q^2(s) {\rm d} s, & \beta \leq 1/2; \\
 \frac{1}{2|\xi|} \int_1^t q(s) {\rm d} s, & \beta>1/2. \end{array}\right.
\]
Then the (modified) wave operator defined by the strong limit in the space $\dot{H}^1(\Rm^+) \times L^2(\Rm^+)$
\begin{align*}
 \vec{\mathbf{W}} \doteq \mathop{{s-}\lim}\limits_{t\rightarrow +\infty} \vec{\mathbf{U}}(t)^{-1} \vec{\mathbf{S}}_V (t) 
\end{align*}
exists and is a unitary bijection from the energy space $\dot{\mathcal{H}}_{V}^1 \times L^2(\Rm^d)$ to $\dot{H}^1(\Rm^d) \times L^2(\Rm^d)$. Furthermore, if $q\in L^1(1,+\infty)$ is a type II repulsive potential, then the same result holds if we substitute $\vec{\mathbf{U}}(t)$ by the classic wave propagation operator 
\[ 
 \vec{\mathbf{S}}_0 (t) = \mathcal{F}_0^{-1} \begin{pmatrix} \cos t|\xi| & |\xi|^{-1} \sin t|\xi| \\ -|\xi| \sin t|\xi| & \cos t|\xi| \end{pmatrix} \mathcal{F}_0 = \begin{pmatrix} \cos t\sqrt{-\Delta} & (-\Delta)^{-1/2} \sin t \sqrt{-\Delta} \\ -\sqrt{-\Delta} \sin t\sqrt{-\Delta} & \cos t\sqrt{-\Delta} \end{pmatrix}.
\]
\end{theorem}

\begin{remark}
 For Schr\"{o}dinger equation with a suitable potential $V(x)$, H\"{o}rmander \cite{LHSchrWO} shows that the modified wave operator can be given by the Fourier multiplier $\exp[{\mathrm i} W(\xi,t)]$ where $W$ is the solution to Hamilton-Jacobi equation 
 \[
  W_t (\xi,t) = |\xi|^2 + V(W_\xi). 
 \]
 As for the case of wave equation, the function $\eta(\xi,t) = |\xi| t + P(|\xi|,t)$ given above is actually an approximated solution to the equation 
 \[
  \eta_t^2 = |\xi|^2 + V(\eta_\xi).
 \]
\end{remark}

Next we give a result on the dispersion rate of free waves, as an application of the theorems given above:

\begin{corollary}[dispersion rate in higher dimensions] \label{main dispersion high}
Assume that $u$ is a finite-energy solution to the wave equation $u_{tt} - \Delta u +V(x) u = 0$ in $\Rm^d$ with $d\geq 3$, where $V(x) = q(|x|)$ is a type II repulsive potential satisfying the decay condition with a decay rate $\beta > 1/3$ and 
\begin{align*}
& \lim_{t\rightarrow +\infty} Q_1(t) = +\infty, & & Q_1(t) = \int_1^t q(s) {\rm d} s.
\end{align*}
Then given any $\varepsilon > 0$ there exists two constants $0<c_1<c_2<+\infty$ such that the following inequalities 
 \begin{align*}
  &\int_{|x|<t-c_2 Q_1(t)} e(x,t) {\rm d} x < \varepsilon;& &\int_{|x|>t-c_1 Q_1(t)} e(x,t) {\rm d} x < \varepsilon
 \end{align*}
 hold for sufficiently large time $t\gg 1$. Here $e(x,t)$ is the energy density function 
 \[
  e(x,t) = |\nabla u(x,t)|^2 + |u_t(x,t)|^2 + V(x) |u(x,t)|^2. 
 \]
 In addition, if $\ell(t)$ is a function satisfying $\ell(t)/Q_1 (t) \rightarrow 0$ as $t\rightarrow +\infty$, then 
 \[
  \lim_{t\rightarrow +\infty} \left(\sup_{r>0} \int_{r<|x|<r+\ell(t)} e(x,t) {\rm d} x\right) = 0. 
 \]
\end{corollary}

\begin{remark}
Let us consider a typical examples. For $1/3<s<2$ we consider the wave equation in $\Rm^d$ with $d\geq 3$
\begin{equation} \label{classic examples s} 
 u_{tt} - \Delta u + \frac{u}{|x|^s} = 0
\end{equation} 
Here the potential $V(x) = |x|^{-s}$ is a type II repulsive potential. The behaviour of $Q_1(t)$ when $t\rightarrow +\infty$ can be described as below. 
\begin{itemize} 
 \item $Q_1(t)$ is uniformly bounded if $s\in (1,2)$;
 \item $Q_1(t)\simeq \ln t$ when $s=1$ and $Q_1(t) \simeq t^{1-s}$ if $s\in (1/3,1)$. 
\end{itemize} 
The asymptotic behaviour of free waves are similar to the classic waves if $s>1$, according to Theorem \ref{thm main high-dimensional}. The majority of energy concentrates in a sphere shell of a constant thickness in this case. In other words, the scattering does not happen in the radial direction at all. Another classic example is the Klein-Gordon equation, which can be also be given by \eqref{classic examples s} with $s=0$. In this case the scattering also fully happens in the radial direction. Corollary \ref{main dispersion high} implies that the wave equation with potential $|x|^{-s}$ when $1/3<s\leq 1$ does scatter in the radial direction, but the dispersion is at a lower rate than the Klein-Gordon equation. 
\end{remark}

\subsection{The Structure of this work}

This paper is organized as the follows. In Section 2 we introduce inward/outward energy theory, which shows that the asymptotic behaviour of free waves only depends on the behaviour of the potential for large $x$ as long as the potential is repulsive. This will help us overcome the difficulty brought by the singularity of the potential near the origin. In Section 3 we assume that the potential is sufficiently good and investigate the behaviour of the corresponding wave functions, which is an important aspect of the spectral theory of one-dimensional differential self-adjoint operators $-{\rm d}^2/{\rm d} x^2 + q(x)$. In Section 4 we utilize the results of Section 3 and prove the existence of (modified) wave operators for type I repulsive potentials in the half-line case and then then incorporate the result of Section 2 to deal with type II and III repulsive potentials. Next in Section 5 we combine the half-line results and the spherical harmonic decomposition to study the asymptotic behaviour of solutions to wave equations with repulsive potentials in higher dimensional case $d\geq 3$. Finally in Section 6 we give the dispersion rate of wave equations with suitable repulsive potentials, as an application of our main result on modified wave operators. 

\section{Inward/outward energy theory} \label{sec: inward} 

In this section we consider the energy distribution theory of solutions to the following abstract wave equation on the half line $\Rm^+$
\begin{equation} \label{wave equation A}
 \left\{\begin{array}{l} w_{tt} + \mathbf{A} w = 0; \\ w(0) = w_0; \\ w_t(0) = w_1.\end{array}\right.
\end{equation}
Here $\displaystyle \mathbf{A} = - {\rm d}^2/{\rm d} x^2 + q(x)$ is a self-adjoint operator with potential $q(x) = \lambda x^{-2} + q_0(x)$ satisfying: 
\begin{itemize} 
 \item The real number $\lambda$ satisfies either $\lambda=0$ or $\lambda\geq 3/4$; 
 \item $q_0$ is a repulsive potential, i.e. $q_0 \in \mathcal{C}^1(\Rm^+)$ satisfies $q_0(x)>0$ and $q'(x)<0$ for all $x\in \Rm^+$ with 
 \[
  \lim_{x\rightarrow +\infty} q_0(x) = 0. 
 \]
 \item $q_0$ satisfies the growth condition near zero ($\kappa \in (0,2)$)
 \[
  q_0(x) \lesssim x^{-\kappa}; \quad |q'_0(x)| \lesssim x^{-\kappa-1}, \qquad \forall x \ll 1. 
 \]
\end{itemize}
We show that for the potentials above, all energy will eventually moves to the infinity at roughly the light speed. Intuitively this implies that the behaviour of solutions when $t$ tends to infinity only depends on the behaviour of the potential $q(x)$ near infinity. Our eventual goal of this section is to show that if two potentials as mentioned above coincide when $x$ is large, then the intervened wave operator exists between finite-energy solutions to the wave equation with these two potentials. 

\begin{remark} 
 To make the operator $\mathbf{A} = - {\rm d}^2/{\rm d} x^2 + q(x)$ self-adjoint in $L^2(\Rm^+)$, boundary condition at $x=0$ is necessary. Roughly speaking we give zero boundary condition. If $\lambda\geq 3/4$, then by Reed-Simon \cite{Simon} the operator $\mathbf{A}$ is self-adjoint with a core $C_0^\infty (\Rm^+)$. The case $\lambda=0$ is slightly different. For the convenience of readers we record the actual domain of these self-adjoint operators, which are also helpful in the further discussion. If $\lambda =0$, then 
 \[
  D(\mathbf{A}) = \left\{w\in AC^2(\Rm^+): \begin{array}{l} w''(x) \in L^1(0,1)\cap L^2(1,+\infty);  w\in L^2(\Rm^+); \\ w(0) = 0; -w''(x) +q w \in L^2(\Rm^+)\end{array} \right\}.
 \]
 If $\lambda \geq 3/4$, then
  \[
  D(\mathbf{A}) = \left\{w\in AC^2(\Rm^+): \begin{array}{l} w''(x) \in L^p(0,1)\cap L^2(1,+\infty), \forall 1\leq p<2;  w\in L^2(\Rm^+); \\ w(0) = 0; w'(0) = 0; -w''(x) + q w \in L^2(\Rm^+)\end{array} \right\}.
 \] 
More details can be found in Appendix, where we verify that $\mathbf{A}$ is indeed a self-adjoint operator with the domain given above. 
In later section discussing higher dimensional case, the additional potential $\lambda x^{-2}$ will naturally emerge if we transform the high dimensional case to half line case by spherically harmonic decomposition. We will also explain why the assumption $\lambda>3/4$ is natural there. 
\end{remark} 

Please note that these operators $\mathbf{A}$ are all positive self-adjoint operators and zero is not an eigenvalue of $\mathbf{A}$. Thus as usual we may define the corresponding Sobolev spaces
\begin{align*}
& \dot{\mathcal{H}}_{\mathbf{A}}^s = \left\{w: \mathbf{A}^{s/2} w \in L^2 \right\};& &\mathcal{H}_{\mathbf{A}}^s = \left\{w: (\mathbf{A}+1)^{s/2} w \in L^2\right\}. 
\end{align*}
In particular we have (please note that $u(0)=0$ for $u\in \dot{H}^1(\Rm^+)$)
\[
 \dot{\mathcal{H}}_{\mathbf{A}}^1 = \left\{w \in \dot{H}^1(\Rm^+): \int_0^\infty \left(|w'(x)|^2 + q(x) |w(x)|^2\right) {\rm d} x < +\infty \right\}.
\]
The solution to the wave equation \eqref{wave equation A} can be given by the functional calculus 
\begin{align*}
 w(t) & = \cos t\sqrt{\mathbf{A}} w_0 + \frac{\sin t \sqrt{\mathbf{A}}}{\sqrt{\mathbf{A}}} w_1; \\
 w_t(t) & = -\sqrt{\mathbf{A}}\sin t\sqrt{\mathbf{A}} w_0 + \cos t\sqrt{\mathbf{A}} w_1. 
\end{align*}
Thus we have the conservation law
\[
 \|w(\cdot,t)\|_{\dot{\mathcal{H}}_{\mathbf{A}}^s}^2 + \|w_t(\cdot, t)\|_{\dot{\mathcal{H}}_{\mathbf{A}}^{s-1}}^2 = \|w_0\|_{\dot{\mathcal{H}}_{\mathbf{A}}^s}^2 + \|w_1\|_{\dot{\mathcal{H}}_{\mathbf{A}}^{s-1}}^2.
\]
In particular, when $s=1$ we obtain the energy conservation law
\[
 E = \int_0^\infty \left(|w_x (x,t)|^2 + |w_t(x,t)|^2 + q(x) |w(x,t)|^2 \right) {\rm d} x. 
\]
For convenience we define the energy density function 
\[
 e(x,t) = |w_x (x,t)|^2 + |w_t(x,t)|^2 + q(x) |w(x,t)|^2. 
\]

\begin{remark} \label{C1 continuity}
 All the norms in the definition of $D(\mathbf{A})$ are all bounded by $\|w\|_{\mathcal{H}_{\mathbf{A}}^2}$. In both cases we also have 
\begin{align*}
 \|w_x\|_{L^2(\Rm^+)} + \sup_{x\in \Rm^+} |w(x)| + \sup_{x\in \Rm^+} |w_x(x)| + |w_x(0)| \lesssim \|w\|_{\mathcal{H}_{\mathbf{A}}^2}.
\end{align*}
In addition, if $\lambda\geq 3/4$, then the following inequalities hold:
\[
 |w(x)| \lesssim \left\{\begin{array}{ll} |x|^{3/2}\|w\|_{\mathcal{H}_{\mathbf{A}}^2}, & \lambda>3/4;\\
 |x|^{3/2}(|\ln x|^{1/2}+1)\|w\|_{\mathcal{H}_{\mathcal{A}}^2}, & \lambda=3/4. \end{array}\right. \quad x\in (0,1).
\]
More details can be found in the appendix. 
\end{remark}

\begin{remark}
 Inward/outward energy theory can also be applied on wave equations with a defocusing nonlinear term. Please refer to Miao-Shen \cite{wavebook} and Shen \cite{shen3dnonradial}. 
\end{remark}

\subsection{Preliminary results}

\begin{lemma}[Morawetz estimates]
 Let $w$ be a solution to \eqref{wave equation A} with a finite energy $E$. Then 
 \[
  \int_{-\infty}^{\infty} |w_x(0,t)|^2 {\rm d} t + \int_{-\infty}^{\infty} \int_0^\infty  \left(- q'(x)\right)|w(x,t)|^2 {\rm d} x {\rm d} t \leq 2E. 
 \]
\end{lemma}
\begin{proof}
 Let us temporarily assume $(w_0,w_1)\in \mathcal{H}_{\mathbf{A}}^2 \times \mathcal{H}_{\mathbf{A}}^1$. Thus we have 
 \begin{align*}
  &(w(\cdot, t),w_t(\cdot, t)) \in \mathcal{C}(\Rm; \mathcal{H}_{\mathbf{A}}^2 \times \mathcal{H}_{\mathbf{A}}^1);& &w_{tt}(\cdot, t) = - \mathbf{A} w \in \mathcal{C}(\Rm; L^2(\Rm^+)).&
 \end{align*}
 Now we consider the inner product 
 \[
  J(t) ={\rm Re} \int_0^\infty 2w_t (x,t) \overline{w_x (x,t)} {\rm d} x
 \]
 By the embedding $\mathcal{H}_{\mathbf{A}}^1 \hookrightarrow \dot{H}^1(\Rm^+)$, the strong derivative of $w_x$ at time $t$ is exactly $w_{tx}$. Thus 
 \begin{align*}
  J'(t) & = {\rm Re} \int_0^\infty \left(2w_{tt} \overline{w_x} + 2w_t \overline{w_{tx}} \right){\rm d} x \\
  & = \lim_{r\rightarrow 0^+, R\rightarrow +\infty} {\rm Re} \int_r^R  \left(2w_{tt} \overline{w_x} + 2w_t \overline{w_{tx}} \right){\rm d} x \\
  & = \lim_{r\rightarrow 0^+, R\rightarrow +\infty} {\rm Re} \int_r^R \left[2\left(w_{xx} - \frac{\lambda}{x^2} w - q_0(x) w\right) \overline{w_x} + \frac{\partial}{\partial x} |w_t|^2 \right]{\rm d} x\\
  & = \lim_{r\rightarrow 0^+, R\rightarrow +\infty} \int_r^R \left[\frac{\partial}{\partial x}\left(|w_x|^2 + |w_t|^2\right) - \left(\frac{\lambda}{x^2} + q_0(x)\right)\frac{\partial}{\partial x} |w|^2 \right] {\rm d} x\\
  & =-  |w_x(0,t)|^2 + \lim_{r\rightarrow 0^+, R\rightarrow +\infty} \int_r^R \left(- \frac{2\lambda}{x^3} + q'_0(x) \right)|w(x,t)|^2 {\rm d} x\\
  & =- |w_x(0,t)|^2  - \int_0^\infty \left(\frac{2\lambda}{x^3} - q'_0(x)\right)|w(x,t)|^2 {\rm d} x. 
 \end{align*}
 Here we use the following asymptotic behaviour 
 \begin{align*}
  &\lim_{x\rightarrow \infty} v_x(x)  = 0, \quad v\in \mathcal{H}_{\mathbf{A}}^2=D(\mathbf{A}); &  &\lim_{x\rightarrow \infty} v(x)  = 0, \quad v \in \mathcal{H}_{\mathbf{A}}^1; \\
  & |v(x)| \lesssim \left\{\begin{array}{ll} x, & \lambda=0, \\ x^{3/2}|\ln x|^{1/2}, & \lambda =3/4; \\ x^{3/2}, & \lambda > 3/4, \end{array}\right. \quad x\ll 1, v\in \mathcal{H}_{\mathbf{A}}^2;& & \lim_{x \rightarrow 0^+} v_x(x) \;\hbox{exist}, \quad v\in \mathcal{H}_{\mathbf{A}}^2;\\
  & \lim_{x\rightarrow 0^+} v(x) = 0, \quad v\in \mathcal{H}_{\mathbf{A}}^1. 
 \end{align*}
An integration of $-J'(t)$ for $t\in [t_1,t_2]$ yields
\begin{equation} \label{Morawetz id 1}
 \int_{t_1}^{t_2} |w_x(0,t)|^2 {\rm d} t + \int_{t_1}^{t_2} \int_0^\infty  \left(- q'(x)\right)|w(x,t)|^2 {\rm d} x {\rm d} t = J(t_1)- J(t_2). 
\end{equation}
Next we observe that 
\[
 |J(t)| \leq 2\|w_t\|_{L^2(\Rm^+)} \|w_x\|_{L^2(\Rm^+)} \leq 2\|w_t\|_{L^2(\Rm^+)} \|w\|_{\dot{\mathcal{H}}_{\mathbf{A}}^1(\Rm^+)} \leq E. 
\]
Inserting this upper bound into \eqref{Morawetz id 1} and letting $t_1\rightarrow -\infty$, $t_2\rightarrow +\infty$, we obtain the desired inequality when the initial data satisfy the stronger assumption $(w_0,w_1)\in \mathcal{H}_{\mathbf{A}}^2 \times \mathcal{H}_{\mathbf{A}}^1$. Finally we observe that $\mathcal{H}_{\mathbf{A}}^2 \times \mathcal{H}_{\mathbf{A}}^1$ is dense in the energy space $\dot{\mathcal{H}}_{\mathbf{A}}^1\times L^2(\Rm^+)$ and conclude the proof. Please note that when we only assume that $w$ comes with a finite energy, $w_x(0,t)$ may be ill-defined at a particular time $t$, but it is still well-defined as an $L^2(\Rm)$ function by continuity.  
\end{proof}

Next we introduce a few technical lemma concerning the regularity of solutions to \eqref{wave equation A}. 
\begin{lemma}
 Assume that $J=[a,b]\subset (0,+\infty)$ be a bounded closed interval. Let $\mathcal{C}_b^k(J)$ be the space of continuous functions with $k$ continuous derivatives with norm 
 \[
  \|w\|_{\mathcal{C}_b^k(J)} = \max\left\{|w^{(j)}(x)|: j=0,1,\cdots, k; \; x\in J\right\}. 
 \]
 Then the following embeddings hold
 \begin{align*}
  \mathcal{H}_{\mathbf{A}}^{k+1} \hookrightarrow \mathcal{C}_b^k (J), \qquad k=0,1,2. 
 \end{align*}
\end{lemma}
\begin{proof}
 When $k=0$, then we have $\mathcal{H}_{\mathbf{A}}^1 \hookrightarrow H^1(\Rm^+) \hookrightarrow \mathcal{C}_b^0 (J)$. The case $k=1$ follows the fact $\mathcal{H}_{\mathbf{A}}^2 = D(\mathbf{A}) \subset AC^2(J)$ and Remark \ref{C1 continuity}. Finally by the definition of $D(\mathbf{A})$, if $w\in \mathcal{H}_{\mathbf{A}}^3$, then the second derivative is defined at least almost everywhere and 
 \[
  -w_{xx} + q(x) w = \mathbf{A} w \in \mathcal{H}_{\mathbf{A}}^1 \hookrightarrow \mathcal{C}_b^0 (J). 
 \]
 Thus $w_{xx}$ can be defined at all points $x\in J$ with 
 \[
  \max_{x\in J} |w_{xx}| \lesssim_{q,J}  \max_{x\in J} |w(x)| +  \max_{x\in J} |(\mathbf{A} w)(x)| \lesssim  \|w\|_{\mathcal{H}_{\mathbf{A}}^1} + \|\mathbf{A} w\|_{\mathcal{H}_{\mathbf{A}}^1} \lesssim \|w\|_{\mathcal{H}_{\mathbf{A}}^3}. 
 \]
 This finishes the proof. 
\end{proof}

\begin{corollary} \label{smoothness away from t}
 Let $(w_0,w_1) \in \mathcal{H}_{\mathbf{A}}^3 \times \mathcal{H}_{\mathbf{A}}^2$.  Then the corresponding solutions $w$ satisfies $w(x,t) \in \mathcal{C}^2(\Rm^+ \times \Rm)$. 
\end{corollary}
\begin{proof}
This immediately follows the previous Lemma. Given a bounded closed interval $J = [a,b]$, we have 
\begin{align*}
 u(t) &\in \mathcal{C}(\Rm; \mathcal{H}_{\mathbf{A}}^3) \hookrightarrow \mathcal{C} (\Rm; \mathcal{C}_b^2(J)); \\
 u_t(t) & \in \mathcal{C}(\Rm; \mathcal{H}_{\mathbf{A}}^2) \hookrightarrow \mathcal{C} (\Rm; \mathcal{C}_b^1(J)); \\
 u_{tt}(t) & \in \mathcal{C}(\Rm; \mathcal{H}_{\mathbf{A}}^1) \hookrightarrow \mathcal{C} (\Rm; \mathcal{C}_b^0(J));
\end{align*}
which immediately gives $\mathcal{C}^2$ continuity away from the origin $x = 0$. 
\end{proof}

\subsection{Energy flux formula}

Now we introduce the definitions of inward/outward energy. Let 
\[
 e_\pm (x,t) = \frac{1}{2} |w_x(x,t)\mp w_t(x,t)|^2  + \frac{q(x)}{2} |w(x,t)|^2;
\]
be the density function and define the outward energy in an interval $J = [a,b]$ at time $t$ to be 
\[
 E_+(J; t) = \int_J e_+ (x,t) {\rm d} x. 
\]
The inward energy can be defined in the same manner. In particular we use the notation $E_-(t)$ and $E_+(t)$ for the inward/outward energy in the total half line. The energy conservation law implies that the sum of inward/outward energies in $\Rm^+$ is exactly the total energy. 
\[
 E_+(t) + E_-(t) = \int_0^\infty \left(e_\pm (x,t) + e_\pm (x,t)\right) {\rm d} x = E. 
\]
Unlike the total energy, the inward/outward energies are not conserved quantities. For convenience we also introduce the notations for the non-directional energy function and the Morawetz density function. 
\begin{align*}
& e'(x,t)  = q(x) |w(x,t)|^2; &
 & M(x,t)  =  - \frac{q'(x)}{2} |w(x,t)|^2 \geq 0.
\end{align*}
The major tool of our inward/outward energy theory is the following energy flux formula. 

\begin{proposition}[Inward/outward energy flux] \label{energy flux formula A}
Let $\Omega\subset [0,+\infty)\times \Rm$ be a region whose boundary $\partial \Omega$ is a simple curve consisting of finite line segments paralleled to the coordinates axes or the lines $x = \pm t$ with clockwise orientation. Then for any finite-energy solution $w$ of \eqref{wave equation A}, we have 
\begin{align*}
  \frac{1}{2} \int_{\partial \Omega} \left(|w_x+w_t|^2 + e'(x,t)\right) {\rm d}x  + \left(+|w_x+w_t|^2 - e'(x,t)\right) {\rm d}t & =  - \iint_{\Omega} M(x,t) {\rm d}x {\rm d}t; & &(\hbox{Inward})\\
  \frac{1}{2} \int_{\partial \Omega} \left(|w_x-w_t|^2 +  e'(x,t)\right) {\rm d}x + \left(-|w_x-w_t|^2 + e'(x,t)\right) {\rm d}t & =  +\iint_{\Omega} M(x,t) {\rm d}x {\rm d}t. & & (\hbox{Outward}) 
 \end{align*}
\end{proposition}
\begin{remark}
 If part of the boundary lies on the $t$-axis, then the corresponding path integral can be understood in a natural way, i.e. this part of integral is equal to
\begin{align*}
 &\frac{1}{2}\int_{t_1}^{t_2} |w_x(0,t)|^2 {\rm d} t;\quad (\hbox{Inward})& & -\frac{1}{2}\int_{t_1}^{t_2} |w_x(0,t)|^2 {\rm d} t.\quad (\hbox{Outward})
\end{align*}
In particular, if $\lambda \geq 3/4$, then the integral must be zero thus can be simply ignored. 
\end{remark} 
\begin{proof}
 We first consider the case with good data $(w_0,w_1)\in \mathcal{H}_{\mathbf{A}}^3 \times \mathcal{H}_{\mathbf{A}}^2$. In this case the solution $w$ is at least $\mathcal{C}^2$ away from the $t$-axis by Corollary \ref{smoothness away from t}. If $\Omega$ is away from the $t$-axis, we may apply Green's formula for the inward case and obtain 
 \begin{align*}
  \hbox{LHS} & = \frac{1}{2} \iint_\Omega \left[\frac{\partial}{\partial t} \left(|w_x+w_t|^2 + e'(x,t)\right) - \frac{\partial}{\partial x} \left(|w_x+w_t|^2 - e'(x,t)\right)\right] {\rm d} x {\rm d} t\\
  & = \frac{1}{2} \iint_\Omega \left[\left(\frac{\partial}{\partial t} -\frac{\partial}{\partial x}\right) |w_x+w_t|^2 +\left(\frac{\partial}{\partial t} + \frac{\partial}{\partial x}\right)  e'(x,t)\right] {\rm d} x {\rm d} t\\
  & = \frac{1}{2} \iint_\Omega \left[\left(\frac{\partial}{\partial t} -\frac{\partial}{\partial x}\right) |w_x+w_t|^2 + q(x) \left(\frac{\partial}{\partial t} + \frac{\partial}{\partial x}\right) |w(x,t)|^2\right] {\rm d} x {\rm d} t \\
  & \qquad + \frac{1}{2} \iint_{\Omega} \frac{\partial}{\partial x} \left( q(x) \right) |w(x,t)|^2 {\rm d} x {\rm d} t\\
  & = {\rm Re} \iint_\Omega \left[\left(w_{tt} - w_{xx} \right)\left(\overline{w_t} + \overline{w_x}\right) +  q(x) w(x,t) \left(\overline{w_t} + \overline{w_x}\right) \right] {\rm d} x {\rm d} t \\
  & \qquad +  \iint_{\Omega}  \frac{q'(x)}{2}  |w(x,t)|^2 {\rm d} x {\rm d} t\\
  & = {\rm Re} \iint_\Omega \left(w_{tt} - w_{xx}  + q(x) w \right)\left(\overline{w_t} + \overline{w_x}\right)  {\rm d} x {\rm d} t +  \iint_{\Omega} \frac{q'(x)}{2} |w(x,t)|^2 {\rm d} x {\rm d} t\\
  & = - \iint_{\Omega} M(x,t) {\rm d} x {\rm d} t. 
 \end{align*}
 The proof of outward case is similar. Now we consider the case when part of the boundary of $\Omega$ lies on the $t$-axis. We first utilize the region $\Omega_r = \Omega \cap ([r,+\infty) \times \Rm)$ and then let $r$ tends to zero. The Morawetz estimates guarantees that 
 \[
  \lim_{r\rightarrow 0^+} \iint_{\Omega_r} M(x,t) {\rm d} x {\rm d} t = \iint_{\Omega} M(x,t) {\rm d} x {\rm d} t. 
 \]
 In addition, Remark \ref{C1 continuity} implies that the following limits/estimates hold uniformly for $t$ in any bounded closed interval as $x \rightarrow 0^+$. 
  \begin{align*}
  & |w(x)| \lesssim \left\{\begin{array}{ll} x, & \lambda=0, \\ x^{3/2}|\ln x|^{1/2}, & \lambda =3/4; \\ x^{3/2}, & \lambda > 3/4; \end{array}\right.& & |w_x(x,t)|  \lesssim 1;\\
  & \lim_{x\rightarrow 0^+} w_t (x,t) = 0; & &\lim_{x \rightarrow 0^+} w_x(x,t)  = w_x(0,t).
 \end{align*}
 Therefore the path integral for $\Omega_r$ converges to the corresponding path integral for $\Omega$ as $r\rightarrow 0^+$. In particular, we have 
 \[
  \lim_{r\rightarrow 0^+} \int_{t_1+\kappa_1 r}^{t_2+\kappa_2 r}  \left(\pm |w_x\pm w_t|^2 \mp e'(r,t)\right) {\rm d}t = \pm \int_{t_1}^{t_2} |w_x(0,t)|^2 {\rm d} t, \qquad \kappa_1, \kappa_2\in \{+1,0,-1\}.
 \]
 This completes the proof if the initial data are sufficiently good. For general initial data $(w_0,w_1)\in \dot{\mathcal{H}}_{\mathbf{A}}^1 \times L^2(\Rm^+)$, we only need to apply the approximation techniques by the fact that $\mathcal{H}_{\mathbf{A}}^3 \times \mathcal{H}_{\mathbf{A}}^2$ is dense in the energy space $\dot{\mathcal{H}}_{\mathbf{A}}^1 \times L^2(\Rm^+)$. 
\end{proof} 

\begin{remark}
 By the Morawetz estimate, if $(w_{0,k}, w_{1,k}) \in \mathcal{H}_{\mathbf{A}}^3 \times \mathcal{H}_{\mathbf{A}}^2$ satisfy 
 \[
  \lim_{k \rightarrow +\infty} \left\|(w_{0,k},w_{1,k})-(w_0,w_1)\right\|_{\dot{\mathcal{H}}_{\mathbf{A}}^1 \times L^2(\Rm^+)} = 0,
 \]
 then the corresponding solutions $w^k$ and $w$ satisfies the following limit uniformly in $t$. 
 \[
  \lim_{k\rightarrow \infty} \int_0^\infty \left(|w_x^k - w_x|^2 + |w_t^k - w_t|^2 + \frac{|w-w^k|^2}{x^2} + q(x) |w-w^k|^2 \right) {\rm d} x = 0.
 \]
 Here we utilize the energy conservation law and the Hardy's inequality. In addition, the Morawetz estimate implies 
 \[
  \lim_{k \rightarrow \infty} \iint_{\Rm^+\times \Rm} \left( - q'(x)\right) |w^k (x,t)- w(x,t)|^2 {\rm d} x {\rm d} t. 
 \]
 If $\lambda = 0$, then we also have 
 \[
  \lim_{k \rightarrow \infty} \int_{-\infty}^\infty |w^k(0,t)-w(0,t)|^2 {\rm d} t = 0. 
 \]
 Therefore the double integral and the path integrals along $t$-axis or the lines $t = t_0$ must converges to the corresponding integral of $u$ as $k\rightarrow +\infty$. The other cases are a little subtle. Let us consider the path integral for inward energy along $x = t + \tau$: 
 \[
  \int_{t_1}^{t_2} |w_x + w_t|^2 {\rm d} t
 \]
 This integral is possibly meaningless for a particular value of $\tau$ if only the finiteness of energy is assumed, since $w_x$ and $x_t$ are merely defined almost everywhere in $\Rm^+$ for each given $t$. Nevertheless, we have 
 \begin{align*}
 \lim_{k \rightarrow +\infty} & \int_{\tau_1}^{\tau_2} \int_{t_1}^{t_2} \left(|(w_x-w_x^k)(t+\tau,t)|^2 + |(w_t-w_t^k)(t+\tau,t)|^2\right) {\rm d} t {\rm d} \tau\\
 & \leq \lim_{k \rightarrow +\infty} \int_{t_1}^{t_2} \int_0^\infty \left(|(w_x-w_x^k)(x,t)|^2 + |(w_t-w_t^k)(x,t)|^2\right) {\rm d} x {\rm d} t = 0. 
 \end{align*}
 It follows that at least for a subsequence of $k$ we have 
 \begin{equation} \label{wk convergence}
  \lim_{k \rightarrow +\infty} \int_{t_1}^{t_2} \left|(w_x^k + w_t^k)(t+\tau,t)\right|^2 {\rm d} t  =  \int_{t_1}^{t_2} \left|(w_x + w_t)(t+\tau,t)\right|^2 {\rm d} t, \quad {\rm a.e.} \; \tau \in [\tau_1,\tau_2]. 
 \end{equation} 
 Thus the energy flux formula holds in the almost everywhere sense. Indeed the limit above holds almost everywhere without extracting a subsequence. An application of the energy flux formula on the solution $w^k- w^j$ and region 
 \[
  \Omega = \left\{(x,t): 0\leq x \leq \tau + t, t_1\leq t\leq t_2\right\}
 \]
 shows that 
 \[
  \lim_{j,k\rightarrow +\infty} \int_{t_1}^{t_2} \left|(\partial_x + \partial_t) (w^j -w^k) (t+\tau,t)\right|^2 {\rm d} t = 0, \qquad \forall \tau\geq -t_1,
 \]
 since all other integrals converge to zero in the energy flux formula. By this convergence we may actually redefine $w_r,w_t$ in a set of measure of zero so that \eqref{wk convergence} holds for all $\tau, t_1, t_2$, without affecting the values of other type integrals in the energy flux formula. The situations for integrals along $x=x_0$ and $x = s -t$ are similar. 
\end{remark}

\begin{remark}[The triangle law] 
 We explain why the identities in Proposition \ref{energy flux formula A} is actually an energy flux formula and give the physical explanation of each integral. We do this by considering a typical example. Let us consider the case $\lambda =0$ and apply Proposition \ref{energy flux formula A} in the inward case on the triangle region 
 \[
  \Omega =\{(x,t)\in [0,+\infty)\times \Rm: x+t \leq s, t\geq t_0\}, \qquad s>t_0. 
 \]
 We obtain 
 \[
   - \int_{0}^{s-t_0} e_-(x,t_0) {\rm d}x + \frac{1}{2} \int_{t_0}^s |w_x(0,t)|^2 {\rm d} t + \int_{t_0}^s e'(s-t,t) {\rm d} t = - \iint_\Omega M(x,t) {\rm d} x {\rm d}t. 
 \]
 Moving the negative term to the other side of identity yields 
 \begin{equation} \label{triangle law}
  E_-([0,s-t_0],t_0) = \frac{1}{2} \int_{t_0}^s |w_x(0,t)|^2 {\rm d} t + \int_{t_0}^s e'(s-t,t) {\rm d} t + \iint_\Omega M(x,t) {\rm d} x {\rm d}t. 
 \end{equation}
 The left hand side is the amount of inward energy contained in the interval $[0,s-t_0]$ at time $t_0$. Clearly there is no energy contained in the single point $0$ at the top of the triangle at time $s$. This energy loss can be divided into there parts, each of them is represented by one term in the right hand of \eqref{triangle law}. The first term is the amount of energy carried by the waves reflected by the boundary point $x=0$. Intuitively when the waves reach the boundary and are reflected, they transform from inward waves to outward waves. The second term is the amount of energy leak at the boundary $x+t = \tau$. Here we use the word ``leak'' because if the energy strictly moved inward, no energy would leave through this boundary. Finally the third term represents the amount of inward energy loss in the space-time region $\Omega$ as inward energy transforms to outward one everywhere and at every time due to the potential effect. 
\end{remark}

\subsection{General theory of inward/outward energy}

\begin{lemma}  \label{Morawetz representation}
 Let $w$ be a finite-energy solution to \eqref{wave equation A}. Given any time $t_0 \in \Rm$, we have the following Morawetz integral representation of inward/outward energy
\begin{align*}
 E_-(t_0) & = \frac{1}{2} \int_{t_0}^\infty |w_x(0,t)|^2 {\rm d} t + \int_{t_0}^\infty \int_0^\infty M(x,t) {\rm d} x {\rm d} t;\\
 E_+(t_0) & = \frac{1}{2} \int_{-\infty}^{t_0} |w_x(0,t)|^2 {\rm d} t + \int_{-\infty}^{t_0} \int_0^\infty M(x,t) {\rm d} x {\rm d} t. 
\end{align*}
\end{lemma}
\begin{proof}
 Let us consider the inward energy case, the outward case can be dealt with in the same way. By approximation techniques and the Morawetz inequality it suffices to consider initial data  $(w_0,w_1) \in \mathcal{H}_{\mathbf{A}}^1 \times (\dot{\mathcal{H}}_{\mathbf A}^{-1} \cap L^2)$, because these data are dense in the energy space. We claim that letting $s \rightarrow +\infty$ in \eqref{triangle law} yields the desired result. It is clear that we only need to show the limit of the second term in the right hand side is zero.  Since the limits of all other three terms exist and are finite, thanks to the Morawetz inequality, the limit 
 \[
  \lim_{s\rightarrow+\infty} \int_{t_0}^s e'(s-t,t) {\rm d} t 
 \]
 also exists. In order to show this limit is zero, it suffices to show this in the average sense, i.e. 
 \begin{equation} \label{to prove average}
  \lim_{R \rightarrow +\infty} \frac{1}{R}\int_{t_0+R}^{t_0+2R} \left(\int_{t_0}^s e'(s-t,t) {\rm d} t\right)  {\rm d} s = 0. 
 \end{equation} 
 We may rewrite this integral in the following form by a change of variable 
 \begin{align*}
  &\frac{1}{R} \int_{\Omega(t_0,R)} e'(x,t) {\rm d} x {\rm d} t;& &\Omega(t_0,R) = \{(x,t)\in \Rm^+\times \Rm: R\leq x+t-t_0 \leq 2R, t\geq t_0\}.
 \end{align*}
 We further split this region into three parts $\Omega(t_0,R) = \Omega_1(t_0,R) \cup \Omega_2(t_0,R) \cup \Omega_3(t_0,R)$ with 
 \begin{align*}
  \Omega_1(t_0,R) & = \left\{(x,t) \in \Omega(t_0,R): |x|<r_1\right\}\\
  \Omega_2(t_0,R) & = \left\{(x,t) \in \Omega(t_0,R): r_1\leq |x|\leq R_1\right\}; \\
  \Omega_3(t_0,R) & = \left\{(x,t) \in \Omega(t_0,R): |x|>R_1\right\}. 
 \end{align*}
 Here $0<r_1<R_1<+\infty$ are constants. We give upper a bound of integral in each part separately 
 \begin{align*}
  \frac{1}{R} \int_{\Omega_1(t_0,R)} e'(x,t) {\rm d} x {\rm d} t & = \frac{1}{R} \int_{\Omega_1(t_0,R)} \left(\lambda \frac{|w(x,t)|^2}{x^2} + q_0(x) |w(x,t)|^2\right) {\rm d} x {\rm d} t\\
  & \lesssim \frac{1}{R} \int_{\Omega_1(t_0,R)} \left(\lambda r_1 \frac{|w(x,t)|^2}{x^3} + x^{-\kappa} (x E)\right) {\rm d} x {\rm d} t\\
  & \lesssim \frac{r_1}{R} \int_{\Omega_1(t_0,R)} M(x,t) {\rm d} x {\rm d} t + r_1^{2-\kappa} E\\
  & \lesssim \left(\frac{r_1}{R} + r_1^{2-\kappa}\right)E. 
 \end{align*}
 The implicit constant above depends only on the potential $q_0(x)$. Here we use the growth condition $|q_0(x)| \lesssim x^{-\kappa}$ and the pointwise estimate
 \[
  |w(x,t)| \leq x^{1/2} \|w(\cdot,t)\|_{\dot{H}^1(\Rm^+)} \leq x^{1/2} \|w(\cdot,t)\|_{\dot{\mathcal{H}}_{\mathbf{A}}^1}. 
 \]
 Next we consider the region $\Omega_{2}(t_0,R)$. Since $q'(x)<0$ and $q(x)>0$ are both continuous in $[r_1,R_1]$, we may find a constant $c=c(r_1,R_1,q)>0$ such that 
 \[
  q(x) \leq c(-q'(x)), \; x\in [r_1,R_1] \quad \Rightarrow \quad e'(x,t) \leq 2c M(x,t), \; x\in [r_1,R_1]. 
 \] 
 It follows that 
 \begin{align*}
  \frac{1}{R} \int_{\Omega_2(t_0,R)} e'(x,t) {\rm d} x {\rm d} t \leq \frac{2c}{R} \int_{\Omega_2(t_0,R)} M(x,t) {\rm d} x {\rm d} t \leq \frac{2c}{R} E. 
 \end{align*}
 Finally we consider $\Omega_3 (t_0,R)$. We have 
 \begin{align*}
  \frac{1}{R} \int_{\Omega_3(t_0,R)} e'(x,t) {\rm d} x {\rm d} t & \leq \frac{1}{R} \int_{t_0}^{t_0+2R} \int_{R_1}^\infty q(x) |w(x,t)|^2 {\rm d} x {\rm d} t\\
  & \leq \frac{q(R_1)}{R} \int_{t_0}^{t_0+2R} \int_{R_1}^\infty |w(x,t)|^2 {\rm d} x {\rm d} t\leq 2q(R_1) \|(w_0,w_1)\|_{L^2\times \dot{\mathcal{H}}_{\mathbf{A}}^{-1}}. 
 \end{align*}
 In summary, we may temporarily fix $r_1<R_1$ and let $R\rightarrow +\infty$ to deduce
 \[
  \limsup_{R\rightarrow +\infty}\frac{1}{R} \int_{\Omega(t_0,R)} e'(x,t) {\rm d} x {\rm d} t \leq C(q_0) r_1^{2-\kappa} E + 2q(R_1) \|(w_0,w_1)\|_{L^2\times \dot{\mathcal{H}}_{\mathbf{A}}^{-1}}. 
 \]
 Letting $r_1\rightarrow 0^+$ and $R_1\rightarrow +\infty$ verifies \eqref{to prove average} and finishes the proof. 
\end{proof}

It immediately follows form Lemma \ref{Morawetz representation} and the identity $E_-(t) +E_+(t) = E$ that 
\begin{corollary} \label{cor limit asymptotic}
 Let $w$ be a finite-energy solution to \eqref{wave equation A}. Then the inward energy $E_-(t)$ is a decreasing function of $t$; while the outward energy $E_+(t)$ is an increasing function of $t$. We also have the limits of inward/outward energy. 
  \begin{align*}
   &\lim_{t\rightarrow \pm \infty} E_\pm (t) = E;& &\lim_{t\rightarrow \mp \infty} E_\pm (t) = 0. 
  \end{align*}
 In particular we have 
 \[
  \lim_{t\rightarrow \pm \infty} \int_0^\infty q(x) |w(x,t)|^2 {\rm d} x = 0. 
 \]
 In addition, the following Morawetz identity holds
\[
   \frac{1}{2} \int_{-\infty}^{\infty} |w_x(0,t)|^2 {\rm d} t + \int_{-\infty}^{\infty} \int_0^\infty M(x,t) {\rm d} x {\rm d} t = E.
 \]
\end{corollary}

\begin{proposition} \label{asymptotic b}
 Let $w$ be a finite-energy solution to \eqref{wave equation A}. Then almost all energy moves away at roughly the light speed as time tends to infinity. More precisely, given any constant $c\in (0,1)$, we have 
  \begin{equation} \label{first limit asymptotic b}
   \lim_{t\rightarrow \pm \infty} \int_{0}^{c|t|} e(x,t) {\rm d} x = 0. 
  \end{equation}
  In addition, we have 
  \begin{align*}
   &\lim_{t\rightarrow \pm \infty} \int_0^\infty \frac{|w(x,t)|^2}{x^2} {\rm d} x = 0;& &\lim_{t\rightarrow \pm \infty} \left(\sup_{x>0} \frac{|w(x,t)|^2}{x}\right) = 0.
  \end{align*}
\end{proposition}
\begin{proof}
 We prove the negative time direction as an example. Again it suffices to consider initial data $(w_0,w_1)\in \mathcal{H}_{\mathbf{A}}^1 \times (L^2 \cap \dot{\mathcal{H}}_{\mathbf{A}}^{-1})$. Applying the triangle law on the region
 \[
  \Omega(t,r) = \left\{(x,t')\in \Rm^+ \times \Rm: x+t'-t \leq r, t'\geq t\right\},\qquad c|t|<r < \frac{1+c}{2} |t|;
 \]
 we obtain 
 \begin{align*}
  E_-([0,r];t) & = \frac{1}{2} \int_t^{t+r} |w_x(0,t')|^2 {\rm d} t' + \int_{t}^{t+r} e'(t+r-t',t') {\rm d} t' + \int_{\Omega(t,r)} M(x,t') {\rm d} x {\rm d} t'\\
  & \leq \frac{1}{2} \int_t^{\frac{1-c}{2}t} |w_x(0,t')|^2 {\rm d} t'+ \int_{t}^{t+r} e'(t+r-t',t') {\rm d} t' + \int_{\Omega(t,\frac{1+c}{2}|t|)} M(x,t') {\rm d} x {\rm d} t'. 
 \end{align*}
 It immediately follows that 
 \begin{align*}
  E_-([0,c|t|];t) &\leq \frac{2}{(1-c)|t|} \int_{c|t|}^{\frac{1+c}{2}|t|} E_-([0,r];t) {\rm d} r \\
  & \leq \frac{1}{2} \int_t^{\frac{1-c}{2}t} |w_x(0,t')|^2 {\rm d} t' + \int_{\Omega(t,\frac{1+c}{2}|t|)} M(x,t') {\rm d} x {\rm d} t'\\
  & \qquad + \frac{2}{(1-c)|t|} \int_{c|t|}^{\frac{1+c}{2}|t|} \left(\int_{t}^{t+r} e'(t+r-t',t') {\rm d} t'\right) {\rm d} r. 
 \end{align*}
 The first two terms in the right hand clearly converges to zero as $t\rightarrow -\infty$ by the integrability of $M(x,t)$ and $|w_x(0,t)|^2$, while the third term also converges to zero by a similar argument to the proof of Lemma \ref{Morawetz representation}. Thus we have 
 \[
  \lim_{t\rightarrow -\infty} E_-([0,c|t|];t) = 0. 
 \]
 Combining this with the already known fact $E_+(t)\rightarrow 0$, we finish the proof of the first conclusion. Now let us prove the second one. We utilize Hardy's inequality and obtain  
 \begin{align*}
  \int_0^\infty \frac{|w(x,t)|^2}{x^2} {\rm d} x & \leq \int_0^{|t|/2} \frac{|w(x,t)|^2}{x^2} {\rm d} x + \int_{|t|/2}^\infty \frac{|w(x,t)|^2}{x^2} {\rm d} x \\
  & \lesssim \int_0^{|t|/2} |w_x(x,t)|^2 {\rm d} x + |t|^{-2} \int_{|t|/2}^\infty |w(x,t)|^2 {\rm d} x\\
  & \lesssim \int_0^{|t|/2} e(x,t) {\rm d} x + |t|^{-2} \|(w_0,w_1)\|_{L^2 \times \dot{\mathcal{H}}_{\mathbf{A}}^{-1}}^2. 
 \end{align*}
 This clearly converges to zero as $t\rightarrow \pm \infty$, where we utilize the first limit \eqref{first limit asymptotic b}. Finally we consider the third inequality. On one hand, we may apply the Cauchy-Schwarz and deduce
 \begin{align*}
  \sup_{x\in (0,|t|/2)} \frac{|w(x,t)|^2}{x} \leq \sup_{x\in (0,|t|/2)} \int_0^x |w_x(x,t)|^2 {\rm d} x \leq \int_0^{|t|/2} e(x,t) {\rm d} x \rightarrow 0, \quad \hbox{as}\;\; t\rightarrow \pm \infty. 
 \end{align*}
 On the other hand, since $w_x(\cdot,t), w(\cdot,t)$ are both uniformly bounded in $L^2(\Rm^+)$ for all $t\in \Rm$, $|w(x,t)|$ is uniformly bounded for all $(x,t) \in \Rm^+\times \Rm$, thus 
 \[
  \sup_{x>|t|/2} \frac{|w(x,t)|^2}{x} \lesssim \sup_{x>|t|/2} \frac{1}{x} \rightarrow 0, \quad \hbox{as}\;\; t\rightarrow \pm \infty. 
 \]
 We combine these two cases and finish the proof. 
\end{proof}

\subsection{Equivalence of asymptotic behaviours}
 
In this subsection we consider two self-adjoint operators $\mathbf{A}_1 = - {\rm d}^2/{\rm d}x^2 + q_1(x)$ and $\mathbf{A}_2 = - {\rm d}^2/{\rm d}x^2 + q_2(x)$ as described at the beginning of this chapter. We prove that if these two potential are different only in a compact set, then their free waves share a 
similar asymptotic behaviour as time tends to infinity. More precisely 

\begin{proposition}\label{local difference A} 
 Assume that $\mathbf{A}_1 = - {\rm d}^2/{\rm d}x^2 + q_1(x)$ and $\mathbf{A}_2 = - {\rm d}^2/{\rm d}x^2 + q_2(x)$ are two self-adjoint operators with potentials as described at the beginning of this chapter. If there exists a number $R>0$ such that 
 \[
  q_1(x) = q_2(x), \qquad x>R,
 \]
 then the wave operator defined by the following strong limit in $\dot{\mathcal{H}}_2^1 \times L^2$
 \[
  \mathbf{T} \doteq \mathop{\rm{s-}\lim}\limits_{t\rightarrow +\infty} \vec{\mathbf{S}}_{2} (-t) \vec{\mathbf{S}}_{1} (t) 
 \]
 is a well-defined linear operator. In addition, $\mathbf{T}: \dot{\mathcal{H}}_1^1 \times L^2 \rightarrow \dot{\mathcal{H}}_2^1 \times L^2$ is an isometric homeomorphism. Here $\dot{\mathcal{H}}_j^1$ and $\vec{\mathbf{S}}_j$ are the corresponding homogeneous Sobolev spaces and wave propagation operators associated to $\mathbf{A}_j$. 
 \end{proposition} 

\begin{remark}
 Please note the assumption $q_1(x) = q_2(x)$ for $x>R$ does not imply the corresponding values of $\lambda$'s are the same, because the value of $\lambda$ determines the behaviour of potential near zero but has very little effect at the infinity. 
\end{remark}

\begin{proof}
Without loss of generality we assume $R=1$. We first show that the wave operator is well-defined. Let $(u_0,u_1)\in \dot{\mathcal{H}}_1^1 \times L^2$ be initial data and $\vec{u} = \vec{\mathbf{S}}_{1} (u_0,u_1)$ be the corresponding free wave. Since the energy space $\dot{\mathcal{H}}_{2}^1\times L^2$ is a complete Hilbert space, it suffices to show $\vec{\mathbf{S}}_{2}(-t) \vec{u}(t)$ is a Cauchy sequence, i.e. 
 \begin{equation} \label{Cauchy local difference} 
  \lim_{t_1, t_2\rightarrow +\infty} \left\|\vec{\mathbf{S}}_{2}(-t_1) \vec{u}(t_1) - \vec{\mathbf{S}}_{2}(-t_2) \vec{u}(t_2)\right\|_{\dot{\mathcal{H}}_{2}^1\times L^2(\Rm^+)} = 0. 
 \end{equation}
 By the inward/outward energy theory Corollary \ref{cor limit asymptotic} and Proposition \ref{asymptotic b}, we have 
 \begin{align*}
  \lim_{t\rightarrow +\infty} \left(E_{1}^-(t) + E_{1}^{+} ([0,2]; t) + \int_{0}^\infty \frac{|w(x,t)|^2}{|x|^2} {\rm d} x\right) = 0. 
 \end{align*}
 Thus given $\varepsilon > 0$, there exists a large time $t_0$, such that 
 \begin{equation} \label{initial data of u} 
  E_{1}^-(t) + E_{1}^{+} ([0,2]; t) + \int_{0}^\infty \frac{|u(x,t)|^2}{|x|^2} {\rm d} x < \varepsilon, \qquad \forall t\geq t_0. 
 \end{equation} 
 Here $E_{1}^\pm$ are the inward/outward energies of $u$. For any time $t_1>t_0$, we let $v$ be the free wave $v = \mathbf{S}_{2}(t-t_1) \vec{u}(t_1)$. It immediately follows that the inward/outward energies of $v$ (with potential $q_2$) satisfy
 \begin{equation} \label{initial data of v} 
  E_{2}^-(t_1) + E_{2}^{+} ([0,2]; t_1) \lesssim \varepsilon. 
 \end{equation} 
 Here we meed to utilize the Hardy's inequality to deduce 
 \[
  \int_0^1 q_2(x) |u(x,t_1)|^2 {\rm d} x \lesssim \int_0^1 \frac{|u(x,t_1)|^2}{x^2} {\rm d} x \lesssim \varepsilon. 
 \]
 It follows from Lemma \ref{Morawetz representation} and \eqref{initial data of u}, \eqref{initial data of v} that 
 \begin{equation} \label{evolution of u and v}
  \int_{t_1}^\infty \left(|u_x(0,t)|^2+|v_x(0,t)|^2\right) {\rm d} t + \int_{t_1}^\infty \int_{0}^\infty \left(M_1 (x,t) + M_2 (x,t)\right) {\rm d} x {\rm d} t \lesssim \varepsilon. 
 \end{equation} 
 Here $M_{1}$ and $M_{2}$ are Morawetz density function associated to $u$ and $v$, respectively. The notations $e'_{j}$ below are defined in the same manner.  Our assumption on the potentials guarantees that there exists a constant $\mu>0$ such that 
 \[
   q_j (x) \leq \mu (-q'_j (x)),\; \forall x \in [1,2], \qquad j=1,2;
 \]
 which implies that 
 \[
  e'_j (x,t) \leq 2\mu M_j (x,t), \qquad x\in [1,2], \; t\in \Rm.
 \]
 From this ineqauality and \eqref{evolution of u and v} we deduce that 
 \begin{equation} \label{cylinder e prime decay}
  \int_{t_1}^\infty \int_{1\leq |x|\leq 2} e'_{j}(x,t) {\rm d} x {\rm d} t \lesssim \varepsilon. 
 \end{equation} 
 An application of the inward/outward energy formula in the cylinder region $[0,r] \times [t_1,+\infty)$ and an integration for $r\in (1,2)$ yields
 \begin{equation} \label{cylinder Lpm decay} 
 \int_{t_1}^\infty \int_1^2 \left(|u_x \pm u_t|^2 + |v_x \pm v_t|^2\right) {\rm d} x {\rm d} t \lesssim \varepsilon. 
 \end{equation}
 Here we use the estimates \eqref{initial data of u}-\eqref{cylinder e prime decay}, which imply that all the other terms in the energy flux formula is dominated by $\varepsilon$, at least in an average sense. Please note that we may ignore the inward/outward energy $E_\pm([0,r];t)$ as $t\rightarrow +\infty$ by the inward/outward energy theory given in Corollary \ref{cor limit asymptotic} and Proposition \ref{asymptotic b}. We then apply the inward/outward energy formula in the finite cylinder region $[0,r] \times [t_1,t_2]$ and integrate for  $r\in (1,2)$ to obtain 
 \begin{equation} \label{inner estimate} 
  \int_1^2 \left(E_{1}^\pm ([0,r]; t_2) + E_{2}^\pm ([0,r];t_2) \right) {\rm d} r \lesssim \varepsilon, \quad \forall t_2 > t_1. 
 \end{equation}
 We then utilize \eqref{initial data of u} to deduce 
 \begin{align*}
  E_{2}^\pm (u; [0,r]; t_2) & \leq E_{1}^\pm (u; [0,r];t_2) + \int_0^1 q_2(x) |u(x,t_2)|^2 {\rm d} x\\
   & \leq E_{1}^\pm (u; [0,r];t_2) + C\int_0^1 \frac{|u(x,t_2)|^2}{x^2} {\rm d} x\\
   & \leq E_{1}^\pm (u; [0,r];t_2) + C\varepsilon.
 \end{align*}
 Here $E_{j}^\pm (u; [0,r];t_2)$ are the inward/outward energy of $u$ with the potential $q_j$. Thus we may integrate and obtain 
 \begin{equation} \label{inner estimate 2} 
  \int_1^2 E_{2}^\pm (u; [0,r];t_2) {\rm d} r \lesssim \varepsilon. 
 \end{equation} 
 Next we observe that $u-v$ solves the wave equation 
 \[
  w_{tt} -  w_{xx} + q_2(x) w = 0 
 \]
 in the exterior region $\Omega = (1,+\infty) \times \Rm$. This enable us to apply the energy flux formula on $u-v$ in the region $[r,+\infty) \times [t_1,t_2]$ for $r>1$ and to deduce 
 \begin{align*}
  E_{2}^\pm (u-v; [r,+\infty); t_2) = & E_{2}^\pm (u-v; [r, +\infty); t_1) \pm \int_{t_1}^{t_2} \int_{r}^\infty M(u-v; x,t) {\rm d} x {\rm d} t\\
   & \pm \frac{1}{2} \int_{t_1}^{t_2} \left(|(\partial_x \mp \partial_t)(u-v)(r,t)|^2 - e'(u-v;r,t) \right) {\rm d} t.
 \end{align*}
 Here $M(u-v; x,t)$ and $e'(u-v; x,t)$ are the corresponding Morawetz and energy density function of $u-v$, with the potential $q_2$. We recall that $\vec{u}\equiv \vec{v}$ at time $t_1$, add the inward/outward part up and integrate to deduce 
 \begin{align*}
 \sum_{\pm} \int_1^2 E_{2}^\pm (u-v; [r,+\infty); t_2) {\rm d} r  \leq \frac{1}{2} \int_{t_1}^{t_2} \int_1^2 \left|(\partial_x -\partial_t)(u-v)(x,t)\right|^2  {\rm d} x {\rm d} t. 
 \end{align*}
 Since the terms $M(u-v;x,t)$, $|(\partial_x \pm \partial_t)(u-v)|^2$ and $e'(u-v;x,t)$ are all quadratic forms of $\vec{u}-\vec{v}$, the estimate \eqref{cylinder Lpm decay} guarantees that 
 \[
  \sum_{\pm} \int_1^2 E_{2}^\pm (u-v; [r,+\infty); t_2) {\rm d} r \lesssim \varepsilon. 
 \]
 Because the density function of $E_{2}^\pm (u-v)$ is also a quadratic forms of $\vec{u}-\vec{v}$, we may also combine \eqref{inner estimate} and \eqref{inner estimate 2} to deduce 
 \[
  \sum_{\pm} \int_1^2 E_{2}^\pm (u-v; [0,r]; t_2) {\rm d} r \lesssim \varepsilon. 
 \]
 A combination of the inner and outer estimates given above immediately yields 
 \[
  \|\vec{u}(t_2)-\vec{v}(t_2)\|_{\dot{\mathcal{H}}_{2}^1 \times L^2}^2 = E_{2} (\vec{u}(t_2) -\vec{v}(t_2)) = \sum_{\pm} E_{2}^\pm (\vec{u}(t_2) -\vec{v}(t_2)) \lesssim \varepsilon. 
 \]
 Here the implicit constant does not depends on the time $t_2>t_1>t_0$. In view of the identity 
 \[
  \vec{u}(t_2)-\vec{v}(t_2) = \vec{u}(t_2) - \vec{\mathbf{S}}_{2}(t_2-t_1) \vec{u}(t_1), 
 \]
 and the fact $\vec{\mathbf{S}}_{2}(t)$ preserves the $\dot{\mathcal{H}}_{2}^1 \times L^2$ norm, we obtain that 
 \begin{align*}
  \left\|\vec{\mathbf{S}}_{2} (-t_2) \vec{u} (t_2) - \vec{\mathbf{S}}_{2} (-t_1) \vec{u}(t_1)\right\|_{\dot{\mathcal{H}}_{2}^1 \times L^2}  = \left\|\vec{\mathbf{S}}_{2}(-t_2)(\vec{u}-\vec{v})(t_2)\right\|_{\dot{\mathcal{H}}_{2}^1\times L^2}  \lesssim \varepsilon^{1/2}
 \end{align*}
 holds for all sufficiently large time $t_2>t_1$. Since $\varepsilon$ is an arbitrary positive constant, we immediately obtain \eqref{Cauchy local difference}. Thus, the wave operator 
 \[
  \mathbf{T} \doteq \mathop{\rm{s-}\lim}\limits_{t\rightarrow +\infty} \vec{\mathbf{S}}_{2} (-t) \vec{\mathbf{S}}_{1} (t) 
 \] 
 is well-defined. It is not difficult to see that its inverse can be given be by the strong limit
  \[
  \mathbf{T}^{-1} \doteq \mathop{\rm{s-}\lim}\limits_{t\rightarrow +\infty} \vec{\mathbf{S}}_{1} (-t) \vec{\mathbf{S}}_{2} (t). 
 \] 
 Here we use the fact 
 \[
  \left\|(u,u_t)\right\|_{\dot{\mathcal{H}}_1^1 \times L^2} \simeq \left\|(u,u_t)\right\|_{\dot{\mathcal{H}}_2^1 \times L^2},
 \]
 thus $\vec{\mathbf{S}}_{k} (-t) \vec{\mathbf{S}}_{j} (t)$ are uniformly bounded operators from $\dot{\mathcal{H}}_j^1 \times L^2(\Rm^+)$ to $\dot{\mathcal{H}}_k^1 \times L^2(\Rm^+)$. Finally we show that $\mathbf{T}$ preserves the norm. Clearly we have
 \begin{align*}
  \|\mathbf{T} (u_0,u_1)\|_{\dot{\mathcal{H}}_{2}^1 \times L^2}  = \lim_{t\rightarrow +\infty} \left\|\vec{\mathbf{S}}_{2} (-t) \vec{u}(t)\right\|_{\dot{\mathcal{H}}_{2}^1\times L^2} = \lim_{t\rightarrow +\infty} \left\|\vec{u}(t)\right\|_{\dot{\mathcal{H}}_{2}^1\times L^2}.
 \end{align*}
 We claim that 
 \[
  \lim_{t\rightarrow +\infty} \left\|\vec{u}(t)\right\|_{\dot{\mathcal{H}}_{2}^1\times L^2} = \lim_{t\rightarrow +\infty} \left\|\vec{u}(t)\right\|_{\dot{\mathcal{H}}_{1}^1\times L^2} = \|(u_0,u_1)\|_{\dot{\mathcal{H}}_{1}^1 \times L^2}. 
 \]
 Indeed, we have 
 \[
  \left\|\vec{u}(t)\right\|_{\dot{\mathcal{H}}_{2}^1\times L^2}^2 = \left\|\vec{u}(t)\right\|_{\dot{\mathcal{H}}_{1}^1\times L^2}^2 + \int_0^1 (q_2(x)-q_1(x)) |u(x,t)|^2 {\rm d} x.
 \]
 Proposition \ref{asymptotic b} then implies 
 \begin{align*}
  \int_{0}^1 |q_2(x)-q_1(x)| |u(x,t)|^2 {\rm d} x \lesssim \int_{0}^1 \frac{|u(x,t)|^2}{|x|^2} {\rm d} x \rightarrow 0.
 \end{align*}
 This immediately finishes the proof.
\end{proof}

\section{Approximation of wave functions}

In this section we investigate the asymptotic behaviour of wave functions when $q$ is a type I repulsive potential which decays like $x^{-\beta}$ for large $x$. This estimate plays an essential role in the subsequent section when we discuss the modified wave operator. 

\subsection{Basic spectrum theory}

Let $q(x)$ be a type I repulsive potential, we now make a very brief review on the spectrum theory of the self-adjoint ordinary differential operators 
\[
 \mathbf{A} = - \frac{{\rm d}^2}{{\rm d}x^2} + q(x)
\]
with boundary condition $w(0) = 0$. For more details of this theory, one may refer to  book \cite{selfadjointbook}. Our assumption on $q$ given above implies that $\mathbf{A}$ is a positive operator thus the spectrum $\sigma(\mathbf{A}) \subseteq [0,+\infty)$ and zero is not an eigenvalue of $\mathbf{A}$. 

\paragraph{Wave functions and spectral measure} Let $u(x,k)$ be the solution to second-order differential equation
\begin{align*}
 -u''(x) + q(x) u(x) = E u(x), \quad u(0) = 0, \quad u'(0) = 1,
\end{align*}
where $E = k^2$. In the argument below it suffices to consider $k>0$ since $\mathbf{A}$ comes with only positive spectrum. We usually call these solutions $u(x,k)$ wave functions of the self-adjoint operator $\mathbf{A}$. Let $R>0$ be a real number. We may consider the truncated version $\mathbf{A}_R$ of the self-adjoint operator $\mathbf{A}$ defined by 
\begin{align*}
 \mathbf{A}_R & = -{\rm d}^2/{\rm d} x^2 + q(x); \\
 D(\mathbf{A}_R) & = \left\{u\in AC^2([0,R]): u''(x) \in L^2([0,R]), u(0) = u(R) = 0\right\}.
\end{align*}
The spectrum of the self-adjoint operator $\mathbf{A}_R$ is of pure point type with single eigenvalues $0<k_{0,R}^2<k_{1,R}^2<k_{2,R}^2<k_{3,R}^2<\cdots$ and eigenfunctions $u(x,k_{j,R})$. Here $u(x,k)$ are the wave functions defined above. The spectrum measure ${\rm d} \rho_R$ of $\mathbf{A}_R$ can be given by 
\[
 \rho_R(E) = \sum_{k_{j,R}^2 \leq E} \frac{1}{\int_0^R |u(x,k_{j,R})|^2 {\rm d} x}. 
\]
The spectral measure ${\rm d}\rho (E)$ of $\mathbf{A}$ can be introduced by a limit of the spectral measure ${\rm d} \rho_R$ as $R\rightarrow +\infty$. This limit is not unique for a general potential $q(x)$ but our assumption on $q(x)$ above guarantees that the corresponding spectral measure ${\rm d} \rho(E)$ is unique. The details of the spectral measure ${\rm d}\rho$ will be discussed later in this section. 

\paragraph{Fourier transforms} Now we are able to define the Fourier transform $\mathcal{F}$ and its inverse $\mathcal{F}^{-1}$ associated to $\mathbf{A}$ by the wave functions $u(x,k)$ and spectrum measure ${\rm d}\rho(E)$ given above. More precisely we have
\begin{align*}
 (\mathcal{F} f)(E) & = \int_0^\infty f(x) u(x,k) {\rm d} x; \\
 (\mathcal{F}^{-1} g)(x) & = \int_0^\infty g(E) u(x,k) {\rm d} \rho(E). 
\end{align*}
The Fourier transform is an isometric bijection between $L^2(\Rm^+)$ and $L^2(\Rm^+; {\rm d}\rho(E))$. In particular we have the generalized Plancherel identity  
\[
 \int_0^\infty |(\mathcal{F}f)(E)|^2 {\rm d} \rho (E) = \int_{0}^\infty |f(x)|^2 {\rm d} x. 
\]
In the argument below, for convenience we will slightly abuse the notation and let $(\mathcal{F} f)(k) = (\mathcal{F} f) (E)$ with $E=k^2$ thus 
\begin{align*}
 (\mathcal{F} f)(k) & = \int_0^\infty f(x) u(x,k) {\rm d} x; \\
 (\mathcal{F}^{-1} g)(x) & = \int_0^\infty g(k) u(x,k) {\rm d} \rho(E). 
\end{align*}

\begin{remark}
 Spectral properties of Schr\"{o}dinger operators $-\Delta + V(x)$ with $V(x)=O(|x|^{-\beta})$ for $\beta>2$ has been studied in previous literature. Please refer to Jensen-Kato \cite{JKSchrOper} and Jenson \cite{JSchrOper, JSchrOperd4}, for instance. 
\end{remark}

\subsection{General eigenfunctions}

 Assume that $q(x)$ is a Type I repulsive potential with $q(x) \lesssim x^{-\beta}$ for large $x$. Given $k>0$, we give an expansion formula of any solutions to the differential equation 
 \begin{equation} \label{general eigenvalue problem}
  -w''(x) + q(x) w(x) = k^2 w(x).
 \end{equation} 
 We start by defining
 \[
  \varphi(k,x) = \exp \left({\mathrm i} k x - \sum_{j=1}^N \frac{{\mathrm i}c_j}{k^{2j-1}} \int_0^x q^j(t) {\rm d} t \right). 
 \]
 Here $N=N(q)$ is a positive integer satisfying $N \beta \geq 1$ and $c_j$'s are real parameters to be determined.  A straightforward calculation shows 
 \begin{align}
  \varphi_x (k,x) & = \left({\mathrm i}k - \sum_{j=1}^N \frac{{\mathrm i}c_j}{k^{2j-1}} q^j(x)\right) \varphi(k,x); \nonumber\\
  \varphi_{xx} (k,x) & = \left({\mathrm i}k - \sum_{j=1}^N \frac{{\mathrm i}c_j}{k^{2j-1}} q^j(x)\right)^2 \varphi(k,x) + \left(-{\mathrm i} \sum_{j=1}^N \frac{j c_j}{k^{2j-1}} q^{j-1}(x) q'(x)\right) \varphi(k,x) \nonumber\\
  & = - \left(k - \sum_{j=1}^N \frac{c_j}{k^{2j-1}} q^j(x)\right)^2 \varphi(k,x) + b(k,x) \varphi_x(k,x). \label{expresion of varphixx}
 \end{align}
 Here 
 \[
  b(k,x) = \frac{\displaystyle -\sum_{j=1}^N \frac{j c_j}{k^{2j-1}} q^{j-1}(x) q'(x)}{\displaystyle k - \sum_{j=1}^N \frac{c_j}{k^{2j-1}} q^j(x)}
 \]
 is a real-valued function and satisfies 
 \[
  \int_R^\infty |b(k,x)| {\rm d} x \lesssim R^{-\beta}, \qquad R>R_0.
 \]
 Given any compact interval $[r_1,r_2]\subset \Rm^+$, the implicit constant and $R_0$ in the inequality can be chosen independent of $k \in [r_1,r_2]$. In fact this is the case for all similar estimates in this subsection. We choose $c_j$'s such that 
 \begin{equation} \label{expansion square}
  \left(k - \sum_{j=1}^N \frac{c_j}{k^{2j-1}} q^j(x)\right)^2 = k^2 - q(x) + \sum_{j=N+1}^{2N} \frac{c'_j}{k^{2j-2}} q^j(x)\doteq k^2 - q(x) + c(k,x).
 \end{equation} 
 Here we list the choice of the first three parameters $c_1 = 1/2$, $c_2 = 1/8$ and $c_3 = 1/16$. It is not difficult to see that all the parameters $c_j$ are positive absolute constants.  By our assumption $N\beta \geq 1$, we also have 
 \[ 
  \int_R^\infty |c(k,x)| {\rm d} x \lesssim R^{-\beta}, \qquad R>R_0.
 \]
 Inserting \eqref{expansion square} into \eqref{expresion of varphixx}, we obtain that both $\varphi(k,x)$ and $\overline{\varphi(k,x)}$ solve the second order differential equation 
 \[
  \mathbf{H} w =  -w''(x) + b(k,x) w'(x) + \left(q(x) -k^2 - c(k,x)\right) w = 0,
 \]
 with Wronskian 
\[
 W(k,x) = \begin{vmatrix} \varphi(k,x) & \overline{\varphi(k,x)} \\ \varphi_x(k,x) & \overline{\varphi_x(k,x)} \end{vmatrix} = -2 {\mathrm i}k + \sum_{j=1}^N  \frac{2{\mathrm i}c_j}{k^{2j-1}} q^j(x).
\]
Let $u$ be a solution of \eqref{general eigenvalue problem}. Then $u$ satisfies 
\begin{align*}
  \mathbf{H} u = b(k,x) u_x (x) - c(k,x) u(x).
\end{align*}
Thus we expect the solution $(u,u_x)$ to satisfy the integral equation
\begin{align*}
 u(x) = & \varphi(k,x)\left[A - \int_x^\infty \overline{\varphi(k,s)} \frac{b(k,s) u_x(s) - c(k,s) u(s)}{W(k,s)}{\rm d} s\right] \\
 &\quad + \overline{\varphi(k,x)} \left[B + \int_x^\infty \varphi(k,s) \frac{b(k,s) u_x(s) - c(k,s) u(s)}{W(k,s)}{\rm d} s\right]; \\
 u_x(x) = & \varphi_x(k,x)\left[A- \int_x^\infty \overline{\varphi(k,s)} \frac{b(k,s) u_x(s) - c(k,s) u(s)}{W(k,s)}{\rm d} s\right]\\
 & \quad + \overline{\varphi_x (k,x)}\left[B+ \int_x^\infty \varphi(k,s) \frac{b(k,s) u_x(s) - c(k,s) u(s)}{W(k,s)}{\rm d} s\right].
\end{align*}
To see that there is such a solution for any complex numbers $A, B$, we consider the space $X = \mathcal{C}([R,+\infty))^2$ of a pair of bounded continuous functions with norm 
\[
 \|(f,g)\|_X = \sup_{x\in [R,+\infty)} \left( |f(x)|+|g(x)|\right)
\]
and a map $\mathbf{T}: X\rightarrow X$ defined by 
\[
 \mathbf{T}(f,g) = \begin{pmatrix}  \displaystyle \varphi(k,x)\left[A - \int_x^\infty \overline{\varphi(k,s)} \frac{b(k,s) g(s) - c(k,s) f(s)}{W(k,s)}{\rm d} s\right] \\
 \displaystyle \qquad + \overline{\varphi(k,x)} \left[B + \int_x^\infty \varphi(k,s) \frac{b(k,s) g(s) - c(k,s) f(s)}{W(k,s)}{\rm d} s\right] \\
 \displaystyle \varphi_x(k,x)\left[A- \int_x^\infty \overline{\varphi(k,s)} \frac{b(k,s) g(s) - c(k,s) f(s)}{W(k,s)}{\rm d} s\right]\\
 \displaystyle \qquad + \overline{\varphi_x (k,x)}\left[B+ \int_x^\infty \varphi(k,s) \frac{b(k,s) g(s) - c(k,s) f(s)}{W(k,s)}{\rm d} s\right]
 \end{pmatrix}.
\]
In view of the the integral estimate of $b(k,x)$ and $c(k,s)$, we obtain the following inequalities for all sufficiently large $R>R_1$:
\begin{align*}
 \|\mathbf{T}(f,g)\|_X & \leq C_1(|A|+|B|) + C_2 R^{-\beta} \|(f,g)\|_X; \\
 \|\mathbf{T}(f_1,g_1) - \mathbf{T}(f_2,g_2)\|_X & \leq C_2 R^{-\beta} \|(f_1,g_1)-(f_2,g_2)\|_X.
\end{align*}
Here the constants $C_1$, $C_2$ and $R_1$ do not depend on $k \in [r_1,r_2]$. We choose a sufficiently large number $R>R_1$ such that $C_2 R^{-\beta}<1/2$. As a result, the map $\mathbf{T}$ becomes a contraction map on $X$, whose unique fixed point immediately gives a solution $(u,u_x)$ to the integral equation, thus a solution to \eqref{general eigenvalue problem}. From the inequality above we see
\[
 \sup_{x\geq R} |u(x)| + |u_x(x)| \leq 2 C_1 (|A|+|B|).  
\]
Inserting this into the integral equation we obtain ($x>R$)
\[
 (u,u_x) = \left(A\varphi(k,x) + B\overline{\varphi(k,x)}, A\varphi_x(k,x) + B\overline{\varphi_x(k,x)}\right) + O\left((|A|+|B|)x^{-\beta}\right). 
\]
Neither implicit constant in the remainder nor $R$ depends on $k\in [r_1,r_2]$. Finally we observe that these solutions span a two-dimensional linear space, thus they are exactly all solutions to \eqref{general eigenvalue problem}. 
\subsection{Wave functions}

In this subsection we investigate the asymptotic behaviour of wave functions $u(k,x)$. We prove 
\begin{lemma} \label{asymptotic ukx}
 Each wave functions $u(k,x)$ can be given by
 \begin{align*}
  (u(k,x),& u_x(k,x))\\
   & = \frac{1}{2 {\mathrm i}} \left(\overline{A(k)}\varphi(k,x) - A(k) \overline{\varphi(k,x)}, \overline{A(k)}\varphi_x(k,x) -A(k)  \overline{\varphi_x(k,x)}\right) + O\left(|A(k)| x^{-\beta}\right)\\
  & = {\rm Im} \left(\overline{A(k)}\varphi(k,x), \overline{A(k)}\varphi_x(k,x) \right) + O\left(|A(k)| x^{-\beta}\right).
 \end{align*}
 Here $A(k)$ is a continuous function of $k>0$ such that $A(k)\neq 0$ for all $k>0$. The remainder term satisfies 
 \[
  O\left(A(k) x^{-\beta}\right) \leq C_3 |A(k)| x^{-\beta}, \qquad x > R.
 \]
 Here $C_3$ and $R$ can be chosen uniformly for all $k$ in any given compact interval $[r_1,r_2] \in (0,+\infty)$. 
\end{lemma}
\begin{proof}
 First of all, the result of last subsection gives the approximation 
 \begin{align*}
 (u(k,x),u_x(k,x)) = & \left(A_\ast (k)\varphi(k,x) + B_\ast (k)\overline{\varphi(k,x)}, A_\ast(k)\varphi_x(k,x) + B_\ast(k)\overline{\varphi_x(k,x)}\right) \\
 & \qquad + O\left((|A_\ast(k)|+|B_\ast(k)|)x^{-\beta}\right).
\end{align*}
Here $A_\ast(k)$ and $B_\ast(k)$ are complex-valued functions of $k>0$. For convenience of further calculation, we let $A_\ast(k) = \overline{A}(k)/2\mathrm i$ and $B_\ast (k) = B(k)/2\mathrm i$ and obtain 
 \begin{align*}
 (u(k,x),u_x(k,x)) = & \frac{1}{2\mathrm i}\left(\overline{A(k)}\varphi(k,x) + B(k)\overline{\varphi(k,x)}, \overline{A(k)}\varphi_x(k,x) + B(k)\overline{\varphi_x(k,x)}\right) \\
 & \qquad + O\left((|A(k)|+|B(k)|)x^{-\beta}\right).
\end{align*}
 Since all wave functions are real-valued functions, the function
\[
 2u(k,x) - 2{\rm Im} \left(\overline{A(k)} \varphi(k,x)\right) = -{\mathrm i} (A(k)+B(k)) \overline{\varphi(k,x)} + O\left((|A(k)|+|B(k)|)x^{-\beta}\right)
\]
is also real-valued. Since we have 
\[
  \lim_{x\rightarrow +\infty} \left(k x - \sum_{j=1}^N \frac{c_j}{k^{2j-1}} \int_0^x q^j(t) {\rm d} t\right) = +\infty,
\]
we may find a sequence $x_n \rightarrow +\infty$, such that 
\[
 -{\mathrm i} (A(k)+B(k)) \overline{\varphi(k,x_n)} = {\mathrm i}\left|A(k) + B(k)\right|, \qquad \forall n \in \mathbb{N}.
\]
If $B(k) + A(k)$ is a nonzero complex number, this gives a contradiction as $n\rightarrow +\infty$. Thus we may plug $B(k) = -A(k)$ and rewrite 
\[
(u(k,x),u_x(k,x)) = {\rm Im} \left(\overline{A(k)}\varphi(k,x), \overline{A(k)}\varphi_x(k,x) \right) + O\left(|A(k)| x^{-\beta}\right).
\]
The estimate on the remainder term immediately follows from the corresponding upper bound estimate given in the previous subsection. Next we show that $A(k)$ is a continuous function of $k$. If this were false, then we might find a sequence $k_j \rightarrow k_0$ with $k_j\in [k_0/2,2k_0]$ but $A_j = A(k_j)$ satisfies 
\[
 |A_0 - A_j| > \delta > 0, \quad \forall j\geq 1. 
\]
This also implies that there exists a constant $c$ independent of $j\geq 1$ such that 
\[
 |A_0| + |A_j| \leq c |A_0 - A_j|, \qquad \forall j\geq 1. 
\]
Thus for sufficiently large $x>R$, we have 
\[ 
 u(k_j,x) - u(k_0,x) = {\rm Im} \left(\overline{A_j} \varphi(k_j,x) - \overline{A_0} \varphi(k_0,x) \right) + O(|A_j-A_0| x^{-\beta}). 
\]
Here the implicit constant in the remainder term does not depend on $j \geq 1$. Thus we have 
\begin{align*}
 \frac{u(k_j,x)-u(k_0,x)}{|A_j-A_0|} =  {\rm Im} \left(\frac{\overline{A_j}}{|A_j-A_0|} \left(\varphi(k_j,x)- \varphi(k_0,x)\right) + \frac{\overline{A_j-A_0}}{|A_j-A_0|} \varphi(k_0,x) \right) + O(x^{-\beta}). 
\end{align*}
Without loss of generality, we may assume that the limit
\[
 \lim_{j\rightarrow +\infty} \frac{\overline{A_j - A_0}}{|A_j-A_0|} 
\]
exists, by extracting a subsequence if neceaary. Again we may choose a large $x>R$, such that $|O(x^{-\beta})|\ll 1$ and
\[
 \lim_{j\rightarrow +\infty} \frac{\overline{A_j-A_0}}{|A_j-A_0|} \varphi(k_0,x) = \mathrm i 
\]
Now we recall that $\varphi(k,x)$ and $u(k,x)$ are both continuous functions of $k>0$ for a fixed number $x>0$, let $j\rightarrow +\infty$ in the identity above and obtain a contradiction. Finally we show that $A(k) =0$ can never happen. If $A(k)=0$, then we would have 
\[
 u(k,x) = 0, \qquad \forall x \gg 1. 
\]
This means that $u\equiv 0$ thus can never happen.  
\end{proof}

\begin{remark} 
 According to Lemma \ref{asymptotic ukx} We may rewrite the wave function and its derivative in the following form 
 \begin{align*}
  u(x,k) & = |A(k)|\left[\sin\left(kx - \sum_{j=1}^N \frac{c_j}{k^{2j-1}} \int_0^x q^j(t) {\rm d} t - \arg A(k)\right) + O(x^{-\beta})\right];\\
  u_x(x,k) & = k |A(k)| \left[\cos\left(kx - \sum_{j=1}^N \frac{c_j}{k^{2j-1}} \int_0^x q^j(t) {\rm d} t - \arg A(k)\right) + O(x^{-\beta})\right].
 \end{align*}
 The asymptotic behaviour of the wave function given above gives the locally uniform limit
 \[
  \lim_{R\rightarrow +\infty} \frac{1}{R}\int_0^R |u(k,x)|^2 {\rm d} x = \frac{|A(k)|^2}{2}.  
 \]
 In addition, for each $k$ in a compact subset $[r_1,r_2]$ of $\Rm^+$, when $x > R(r_1,r_2)$ is large, there is exactly one zero of $u(x,k)$ around each $x$ satisfying 
 \[
  kx - \sum_{j=1}^N \frac{c_j}{k^{2j-1}} \int_0^x q^j(t) {\rm d} t - \arg A(k) \in \pi \mathbb{N}. 
 \]
 Combining this with the fact that the $j$-th  eigenfunction $u(x,k_{j,R})$ for truncated self-adjoint operator $\mathbf{A}_R$ comes with exactly $j$ zeros in the interval $(0,R)$, as given in Section 1.3 of \cite{selfadjointbook}, we may deduce that two consecutive eigenvalues $k_{j,R}^2, k_{j+1,R}^2 \in [r_1^2,r_2^2]$ of $\mathbf{A}_R$ satisfy 
 \[
  k_{j+1,R} - k_{j,R} \approx \pi/R, \qquad R > R(r_1,r_2). 
 \]
 By the way we introduce the spectrum measure ${\rm d}\rho(E)$ from its truncated version ${\rm d}\rho_R (E)$, we may deduce that the spectrum measure $\rho(E)$ must be purely absolutely continuous and given by
 \begin{equation} \label{measure identity 2}
  {\rm d} \rho (E) = \frac{E^{-1/2}}{\pi |A(k)|^2} {\rm d} E = \frac{2}{\pi |A(k)|^2} {\rm d} k. 
 \end{equation} 
Indeed, the well known Lavine's theorem (see Theorem XIII.29 in \cite{Simon}) implies that a self-adjoint operator $-\Delta + V$ with a suitable repulsive potential $V$ always comes with a purely absolutely continuous spectrum. 
\end{remark}

\begin{remark}
 The function $A(k)$ is actually an alternative version of modified Jost function $j_m (k)$. If $N=1$, then we have $A(k) = k^{-1} j_m(k)$. For more details about the modified Jost function, please refer to \cite{WPDen, KS}. 
\end{remark}

\section{Modified wave operators in half-line}

In this section we consider the modified wave operator for repulsive potentials with decay rate higher than $1/3$. 

\subsection{Half wave operators} 

Let us start by introducing a few notations 
\begin{align*}
 Q_j (x) = \int_0^x q^j(x) {\rm d} x.
\end{align*}
Let $P(k,t)$ be the phase shift function 
\[
 P(k,t) = \frac{1}{2k} Q_1(t) - \frac{1}{4k^3} q(t) Q_1(t) + \frac{1}{8k^3} Q_2(t) + \frac{1}{16k^5} \int_0^\infty q^3(x) {\rm d} x. 
\]
We define the approximated half-wave operator 
\begin{align*}
 \mathbf{U}(t) f & = \mathcal{F}_0^{-1} \left({\mathrm e}^{{\mathrm i}k t + {\mathrm i} P(k,t)} \mathcal{F}_0 f\right)\\
 & = \frac{2}{\pi} \int_0^\infty \sin kx \left({\mathrm e}^{{\mathrm i}k t + {\mathrm i} P(k,t)}(\mathcal{F}_0 f) (k)\right) {\rm d} k;
\end{align*}
and the operator
\begin{equation} \label{def of W}
  \mathbf{W} = \mathcal{F}_0^{-1} \frac{1}{\;\overline{A(k)}\;} \mathcal{F}
 \end{equation}
By the spectrum measure given by \eqref{measure identity 2} one may verify that the operator $\mathbf{W}$ is actually a unitary operator from $L^2(\Rm^+)$ to itself, and from $\dot{\mathcal{H}}_{\mathbf{A}}^1$ to $\dot{H}^1(\Rm^+)$. Our main result of this subsection is 

\begin{proposition}
 Assume that $q(x)$ is a type I repulsive potential with decay rate $\beta>1/3$. Let $\mathbf{A} = -{\rm d}^2/{\rm d} x^2 + q(x)$ be the self-adjoint operator on $L^2(\Rm^+)$ with zero boundary condition. Then $\mathbf{W}$ defined above is exactly the (modified) wave operator defined by the strong limit 
 \[
  \mathop{{\rm s-}\lim}\limits_{t\rightarrow +\infty} \mathbf{U}(t)^{-1} {\mathrm e}^{{\mathrm i} t \sqrt{\mathbf{A}}} 
 \]
in $L^2(\Rm^+)$. A similar result holds from $\dot{\mathcal{H}}_{\mathbf{A}}^1$ to $\dot{H}^1(\Rm^+)$. 
\end{proposition}

\begin{proof}
Let us first consider the strong limit in the $L^2(\Rm^+)$. Since all the operators involved, i.e. the half wave operator ${\mathrm e}^{{\mathrm i} t \sqrt{\mathbf{A}}} $, $\mathbf{U}(t)$ and $\mathbf{W}$ are unitary operators in $L^2(\Rm^+)$, it suffices to prove the limit 
\begin{equation} \label{half wave to prove m1}
 \lim_{t\rightarrow +\infty} \left\|\mathbf{U}(t) \mathbf{W} u_0 - {\mathrm e}^{{\mathrm i} t \sqrt{\mathbf{A}}}  u_0 \right\|_{L^2(\Rm^+)} = 0
\end{equation}
for $u_0$ in a dense subset of $L^2(\Rm^+)$. Let us consider $u_0$ such that its Fourier transform $(\mathcal{F} u_0)(k)$ is supported in a finite interval $[a,b]\subset (0,+\infty)$ and $(\mathcal{F} u_0)(k)/A(k)$ is a smooth function of $k$. To see these initial data are dense, one need to use the fact that $A(k)$ is a continuous function with no zeros. We claim that for any $\tau > \max\{1-\beta,0\}$, there holds
\begin{equation} \label{outside decay Ut}
 \lim_{t\rightarrow +\infty} \left\|\mathbf{U}(t) v_0\right\|_{L^2((0,t-t^\tau)\cup(t+t^\tau,+\infty))} = 0, \qquad \forall v_0\in L^2(\Rm^+). 
\end{equation} 
Indeed, if $\mathcal{F}_0 v_0$ is smooth and supported in a bounded interval $[a,b]\subset (0,+\infty)$, then an integration by parts gives the following decay for large time $t$ and $x\in (0,t-t^\tau)\cup(t+t^\tau,+\infty)$
\begin{align*}
 (\mathbf{U}(t) v_0)(x) 
 \lesssim \frac{1}{|x-t|}, 
\end{align*}
which gives the desired decay in $L^2$. One may immediately verify the general case by observing these data are dense in $L^2(\Rm^+)$. For convenience we let $J_t = (t-t^\tau,t+t^\tau)$. Here $\tau$ is a constant slightly greater than $\max\{1-\beta,0\}$. Next we claim that the limit \eqref{half wave to prove m1} holds as long as 
\begin{equation} \label{half wave to prove m2}
 \lim_{t\rightarrow +\infty} \left\|\mathbf{U}(t) \mathbf{W} u_0 - {\mathrm e}^{{\mathrm i} t \sqrt{\mathbf{A}}}  u_0 \right\|_{L^2(J_t)} = 0.
\end{equation}
Indeed, in view of \eqref{outside decay Ut}, the limit \eqref{half wave to prove m2} implies that 
\[
 \lim_{t\rightarrow +\infty}\left\| {\mathrm e}^{{\mathrm i} t \sqrt{\mathbf{A}}} u_0 \right\|_{L^2(J_t)} = \lim_{t\rightarrow +\infty} \left\|\mathbf{U}(t) \mathbf{W} u_0 \right\|_{L^2(J_t)} = \|u_0\|_{L^2(\Rm^+)}.
\]
Thus 
\[
 \lim_{t\rightarrow +\infty}\left\| {\mathrm e}^{{\mathrm i} t \sqrt{\mathbf{A}}} u_0 \right\|_{L^2(\Rm^+ \setminus J_t)} = 0.
\]
A combination of this with \eqref{outside decay Ut} and \eqref{half wave to prove m2} then verifies \eqref{half wave to prove m1}. Next we prove \eqref{half wave to prove m2}. We start by writting $u(x,t) \doteq {\mathrm e}^{{\mathrm i} t \sqrt{\mathbf{A}}}  u_0$ in details for $x\in J_t$ when $t$ is large by using the wave functions expression given in Lemma \ref{asymptotic ukx} and the spectrum measure \eqref{measure identity 2}:
\begin{align*}
  u(x,t) & = \mathcal{F}^{-1} \left({\mathrm e}^{{\mathrm i} t k} \mathcal{F} u_0 \right)\\
  & = \frac{1}{2 \mathrm i} \int_{a^2}^{b^2}  \left(\overline{A(k)} \varphi(k,x) - A(k) \overline{\varphi(k,x)} + O(x^{-\beta})\right) {\mathrm e}^{{\mathrm i} t k} \mathcal{F} u_0 {\rm d} \rho(E)\\
  & = \frac{1}{\pi \mathrm i} \int_a^b \left(\overline{A(k)} \varphi(k,x) - A(k) \overline{\varphi(k,x)} + O(x^{-\beta})\right) {\mathrm e}^{{\mathrm i} t k} \mathcal{F} u_0 \frac{{\rm d} k}{|A(k)|^2}\\
  & = \frac{1}{\pi \mathrm i} \int_a^b \left(\frac{\mathcal{F} u_0}{A(k)} \varphi(k,x)  {\mathrm e}^{{\mathrm i} t k} - \frac{\mathcal{F} u_0}{\;\overline{A(k)}\;} \overline{\varphi(k,x)}  {\mathrm e}^{{\mathrm i} t k}\right) {\rm d} k + O(x^{-\beta}).
 \end{align*}  
We claim that the remainder term can be ignored, because we have 
\[
 \|x^{-\beta}\|_{L^2(J_t)} \simeq t^{\frac{\tau}{2}-\beta}.
\]
The exponent $\tau/2 -\beta$ is negative if $\tau$ is slightly greater than $\max\{1-\beta,0\}$ and $\beta>1/3$. Similarly the first term in the integral can also be ignored because an integration by parts shows that the integral of the first term decays like $1/(t+x)$. Thus we may rewrite 
\begin{equation} \label{approximation wave 1}
 u(x,t) \approx \frac{\mathrm i}{\pi} \int_a^b \frac{\mathcal{F} u_0}{\;\overline{A(k)}\;} \exp\left[\mathrm i \left(kt - kx + \frac{1}{2k} Q_1(x) + \frac{1}{8k^3} Q_2(x) + \frac{1}{16k^5} Q_3(x) \right)\right] {\rm d} k.
\end{equation}
Here the notation ``$u(x,t) \approx\cdots$'' means that the $L^2(J_t)$ norm of the difference vanishes as $t\rightarrow +\infty$. Next we let 
\begin{align*}
 Q_1(x) & = Q_1(t) + (x-t) q(t) + Q_1^\ast (x,t);\\
 Q_2(x) & = Q_2(t) + Q_2^\ast (x,t); \\
 Q_3(x) & = \int_0^\infty q^3(s) {\rm d} s + Q_3^\ast (x). 
\end{align*}
When $x\in J_t$ we utilize the decay of $q$ and $q'$ to deduce 
\begin{align*}
 Q_1^\ast (x,t) & = \int_t^x q'(s) (x-s) {\rm d} s = O(t^{2\tau-1-\beta}); \\
 Q_2^\ast (x,t) & = \int_t^x q^2(s) {\rm d} s = O(t^{\tau-2\beta}); \\
 Q_3^\ast (x) & = \int_x^\infty q^3(s) {\rm d} s = O(t^{1-3\beta}). 
\end{align*}
Since $\beta>1/3$ and $\tau$ is slightly larger than $\max\{0,1-\beta\}$, all the exponents of $t$ above are negative. Please note that if $\beta >1/2$, we do not assume the decay of $q'(x)$ but we still have 
\[
 |Q_1^\ast (x,t)| = \left|\int_t^x q(s) {\rm d} s - (x-t) q(t)\right| \leq t^{\tau-\beta}.
\]
Thus $Q_1^\ast(x,t)$ can also be bounded by a negative power of $t$. We then write the exponent in \eqref{approximation wave 1} as 
\begin{align*}
 kt & - kx + \frac{1}{2k} Q_1(x) + \frac{1}{8k^3} Q_2(x) + \frac{1}{16k^5} Q_3(x) \\
& = (t-x) \left(k - \frac{1}{2k} q(t)\right) + \frac{1}{2k} Q_1(t) + \frac{1}{8k^3} Q_2(t) + \frac{1}{16k^5} \int_0^\infty q^3(s) {\rm d} s\\
& \qquad + \frac{1}{2k} Q_1^\ast (x,t) + \frac{1}{8k^3} Q_2^\ast (x,t) + \frac{1}{16k^5} Q_3^\ast (x) \\
& \doteq R(k,x,t) + \frac{1}{2k} Q_1^\ast (x,t) + \frac{1}{8k^3} Q_2^\ast (x,t) + \frac{1}{16k^5} Q_3^\ast (x).
\end{align*}
We then write the exponential part in \eqref{approximation wave 1} accordingly 
\begin{align*}
&\exp  \left[\mathrm i \left(kt - kx + \frac{1}{2k} Q_1(x) + \frac{1}{8k^3} Q_2(x) + \frac{1}{16k^5} Q_3(x) \right)\right] \\
 & = \exp \left({\mathrm i} R(k,x,t)\right) \left(\sum_{j_1=0}^{N} \frac{{\mathrm i}^{j_1} Q_1^\ast(x,t)^{j_1}}{j_1! (2k)^{j_1}}\right) \left(\sum_{j_2=0}^{N} \frac{{\mathrm i}^{j_2} Q_2^\ast(x,t)^{j_2}}{j_2! (8k^3)^{j_2}}\right) \left(\sum_{j_3=0}^{N} \frac{{\mathrm i}^{j_3} Q_3^\ast(x)^{j_3}}{j_3! (16k^5)^{j_3}}\right) + O(t^{-\beta}). 
\end{align*}
Here we choose $N$ to be very large integers such that the remainder term of the power expansions are all dominated by $t^{-\beta}$, where we use the power-like decay of $Q_1^\ast, Q_2^\ast, Q_3^\ast$. Again the remainder $t^{-\beta}$ can be ignored. Thus we have 
\begin{align} \label{expansion u}
 u(x,t) \approx \frac{\mathrm i}{\pi} \sum_{j_1,j_2,j_3=0}^N \frac{Q_1^\ast(x,t)^{j_1} Q_2^\ast(x,t)^{j_2} Q_3^\ast(x,t)^{j_3}}{(-\mathrm i)^{j_1+j_2+j_3} j_1! j_2! j_3!} \int_a^b \frac{(\mathcal{F} u_0)\exp ({\mathrm i}R(k,x,t))}{(2k)^{j_1} (8k^3)^{j_2} (16k^5)^{j_3}\overline{A(k)}} {\rm d}k. 
\end{align}
We claim that the $L^2(J_t)$ norm of each integral in the sum above is uniformly bounded. In fact we may let 
\begin{align*}
 R_1(k,t) & = \frac{1}{2k} Q_1(t) + \frac{1}{8k^3} Q_2(t) + \frac{1}{16k^5} \int_0^\infty q^3(s) {\rm d} s;\\
 F_{j_1 j_2 j_3} (k,t) & = \frac{\exp \left({\mathrm i} R_1(k,t)\right) (\mathcal{F} u_0)(k)}{(2k)^{j_1} (8k^3)^{j_2} (16k^5)^{j_3}\overline{A(k)}}
\end{align*}
and rewrite the integral as 
\begin{align*}
 \int_a^b \frac{(\mathcal{F} u_0)\exp ({\mathrm i}R(k,x,t))}{(2k)^{j_1} (8k^3)^{j_2} (16k^5)^{j_3}\overline{A(k)}} {\rm d}k = \int_a^b F_{j_1 j_2 j_3} (k,t) \exp \left[{\mathrm i}(t-x) \left(k - \frac{1}{2k} q(t)\right)\right]{\rm d} k.
\end{align*}
For sufficiently large $t$, we then apply a change of variables 
\[
 k_1 = k - \frac{1}{2k} q(t) \quad \Longleftrightarrow \quad k = \frac{k_1 + \sqrt{k_1^2 + 2q(t)}}{2} 
\]
and rewrite the integral in the form of 
\begin{equation} \label{after change} 
\int_{a-q(t)/2a}^{b-q(t)/2b} \exp \left({\mathrm i}k_1(t-x)\right) \left(\frac{1}{2} + \frac{k_1}{2\sqrt{k_1^2+2q(t)}}\right) F_{j_1 j_2 j_3} \left(\frac{k_1 + \sqrt{k_1^2 + 2q(t)}}{2}, t\right) {\rm d} k_1. 
\end{equation}
Its $L^2(\Rm^+)$ norm can be dominated by
\[
 \left\|\exp \left({\mathrm i}k_1 t\right) \left(\frac{1}{2} + \frac{k_1}{2\sqrt{k_1^2+2q(t)}}\right) F_{j_1 j_2 j_3} \left(\frac{k_1 + \sqrt{k_1^2 + 2q(t)}}{2}, t\right)\right\|_{L_{k_1}^2([a-q(t)/2a, b-q(t)/2b])},
\]
which is uniformly bounded for all large $t$, since all the functions above are uniformly bounded, so is the length of the interval $[a-q(t)/2a, b-q(t)/2b]$. Combining this uniform bound of $L^2(J_t)$ norm and the fact that $Q_1^\ast, Q_2^\ast, Q_3^\ast$ decays like a negative power of $t$, we conclude that all the terms in the expansion \eqref{expansion u} except for $j_1=j_2=j_3=0$ can be ignored as $t\rightarrow \infty$. Thus we may recall \eqref{after change} and write 
\[
 u(x,t) \approx \frac{\mathrm i}{\pi} \int_{a-\delta}^{b} \exp \left({\mathrm i}k_1(t-x)\right) \left(\frac{1}{2} + \frac{k_1}{2\sqrt{k_1^2+2q(t)}}\right) F_{000} \left(\frac{k_1 + \sqrt{k_1^2 + 2q(t)}}{2}, t\right) {\rm d} k_1.
\]
Here we slightly enlarge the integral interval. This will not affect the value of the integral because $(\mathcal{F} u_0)(k)$ is zero if $k_1$ is not in the interval $[a-q(t)/2a, b-q(t)/2b]$. Next we observe that 
\[
 \frac{1}{k} = \frac{1}{k_1} - \frac{1}{2k_1^3} q(t) + O(t^{-2\beta}).
\]
Thus 
\begin{align*}
 R_1(k, t) & = \frac{1}{2k} Q_1(t) + \frac{1}{8k^3} Q_2(t) + \frac{1}{16k^5} \int_0^\infty q^3(s) {\rm d} s\\
& = \frac{1}{2k_1} Q_1(t) - \frac{1}{4k_1^3} q(t) Q_1(t)  + \frac{1}{8k_1^3} Q_2(t) + \frac{1}{16k_1^5} \int_0^\infty q^3(s) {\rm d} s +O(t^{-\nu});\\
& = P(k_1,t) + O(t^{-\nu}).
\end{align*}
Here $\nu = \nu(\beta)>0$ is a constant. As a result, we have
\begin{align*}
 F_{000} (k,t) & = \frac{\exp \left({\mathrm i} R_1(k,t)\right) (\mathcal{F} u_0)(k)}{\overline{A(k)}} \\
 & = \frac{\big(\exp \left({\mathrm i} P(k_1,t)\right)+ O(t^{-\nu})\big) (\mathcal{F} u_0)\left(\frac{k_1 + \sqrt{k_1^2 + 2q(t)}}{2}\right)}{\bar{A}\left(\frac{k_1 + \sqrt{k_1^2 + 2q(t)}}{2}\right)}.
\end{align*}
It follows that 
\[
  \lim_{t\rightarrow +\infty} \sup_{k_1\in [a-\delta,b]} \left|F_{000} \left(\frac{k_1 + \sqrt{k_1^2 + 2q(t)}}{2},t\right) - \frac{\exp \left({\mathrm i} P(k_1,t)\right) (\mathcal{F} u_0)\left(k_1\right)}{\overline{A\left(k_1\right)}}\right| = 0.
\]
Here we use the assumption/fact that $\mathcal{F} u_0$ and $1/A(k)$ are all continuous functions in $\Rm^+$. In addition we have
\[
 \frac{1}{2} + \frac{k_1}{2\sqrt{k_1^2+2q(t)}} = 1 + O(q(t)/k_1^2) = 1 + O(t^{-\beta});
\]
In summary the difference 
\begin{align*}
 \exp \left({\mathrm i}k_1t\right)  \left(\frac{1}{2} + \frac{k_1}{2\sqrt{k_1^2+2q(t)}}\right) F_{000} & \left(\frac{k_1 + \sqrt{k_1^2 + 2q(t)}}{2}, t\right) \\
& - \exp \left({\mathrm i}k_1t\right) \frac{\exp \left({\mathrm i} P(k_1,t)\right) (\mathcal{F} u_0)\left(k_1\right)}{\overline{A\left(k_1\right)}}
\end{align*}
converges uniformly for $k_1\in [a-\delta,b]$ to zero as $t\rightarrow +\infty$, thus also converges to zero in $L^2([a-\delta,b])$. By the Parseval equality we have 
\begin{align}
 u(x,t) & \approx \frac{\mathrm i}{\pi} \int_{a-\delta}^b \exp \left({\mathrm i}k_1(t-x)\right) \frac{\exp \left({\mathrm i} P(k_1,t)\right) (\mathcal{F} u_0)\left(k_1\right)}{\overline{A\left(k_1\right)}} {\rm d} k_1\nonumber\\
 & \approx \frac{\mathrm i}{\pi} \int_a^b \exp \left[{\mathrm i} \left(-k x + k t + P(k,t)\right)\right] \frac{(\mathcal{F} u_0)(k)}{\overline{A(k)}} {\rm d} k. \label{approximation u 2}
\end{align}
Next we claim that the following limit holds as $t\rightarrow +\infty$ for all $f(k) \in L^2([a,b])$:
\begin{equation} \label{u add 2}
 \left\|\int_a^b \exp \left[{\mathrm i} \left(k x + k t + P(k,t)\right)\right] f(k) {\rm d} k\right\|_{L^2(J_t)} \rightarrow 0. 
\end{equation}
Indeed, this is clear for $f\in C_0^\infty (a,b)$ since an integration by parts shows that the integral decays no slower than $t^{-1}$. The general case follows a smooth approximation technique and the Parseval equality 
\[
 \left\|\int_a^b \exp \left[{\mathrm i} \left(k x + k t + P(k,t)\right)\right] f(k) {\rm d} k\right\|_{L^2(J_t)} \lesssim \left\| \exp \left[{\mathrm i} \left(k t + P(k,t)\right)\right] f(k)\right\|_{L^2([a,b])} =\|f\|_{L^2([a,b])}.
\]
Combining \eqref{approximation u 2} and \eqref{u add 2}, we have 
\[
 u(x,t) \approx \frac{2}{\pi} \int_a^b \sin kx \left({\mathrm e}^{{\mathrm i}k t + {\mathrm i} P(k,t)}(\mathcal{F}_0 \mathbf{W} u_0) (k)\right) {\rm d} k = \left(\mathbf{U}(t) \mathbf{W} u_0\right)(x).
\]
This verifies \eqref{half wave to prove m2} thus completes the proof of $L^2(\Rm^+)$ case. Next we show that a similar result also holds from $\dot{\mathcal{H}}_{\mathbf{A}}^1$ to $\dot{H}^1(\Rm^+)$. More precisely the following limit holds in the space $\dot{H}^1$:
 \begin{equation} \label{to prove H1}
  \mathbf{W} u_0 = \lim_{t\rightarrow +\infty} \mathbf{U}(t)^{-1} {\mathrm e}^{{\mathrm i} t \sqrt{\mathbf{A}}} u_0, \qquad \forall u_0 \in \dot{\mathcal{H}}_{\mathbf{A}}^1. 
 \end{equation}
 Again we observe that $\mathbf{W}$ is a unitary operator from $\dot{\mathcal{H}}_{\mathbf{A}}^1$ to $\dot{H}^1(\Rm^+)$ and that the operator norm of $\mathbf{U}(t)^{-1} {\mathrm e}^{{\mathrm i} t \sqrt{\mathbf{A}}}$ from the space $\dot{\mathcal{H}}_{\mathbf A}^1$ to $\dot{H}^1(\Rm^+)$ is bounded by $1$ for all $t$. Thus it suffices to prove \eqref{to prove H1} for initial data $u_0$ in a dense subset of $\dot{\mathcal{H}}_{\mathbf A}^1$. Let us consider the initial data whose Fourier transform $\mathcal{F} u_0$ is supported in a compact interval $[a,b]\subset (0,+\infty)$. Let $\mathbf{P}_{J}$ be the classic frequency cut-off operator  
\[
 \mathbf{P}_J = \mathcal{F}_0^{-1} \chi_J \mathcal{F}_0,
\]
where $\chi_J$ is the characteristic function of $J$. Then our $L^2$ theory above immediately gives 
\[
 \left\|\mathbf{P}_{[0,b]} \left(\mathbf{U}(t) \mathbf{W} u_0 - {\mathrm e}^{{\mathrm i} t \sqrt{\mathbf{A}}}  u_0\right)\right\|_{\dot{H}^1} \leq b \left\|\mathbf{U}(t) \mathbf{W} u_0 - {\mathrm e}^{{\mathrm i} t \sqrt{\mathbf{A}}}  u_0\right\|_{L^2(\Rm^+)} \rightarrow 0.
\]
As a result, we have 
\[
 \lim_{t\rightarrow +\infty} \left\|\mathbf{P}_{[0,b]} {\mathrm e}^{{\mathrm i} t \sqrt{\mathbf{A}}}  u_0\right\|_{\dot{H}^1} = \lim_{t\rightarrow +\infty} \left\|\mathbf{P}_{[0,b]} \mathbf{U}(t) \mathbf{W} u_0 \right\|_{\dot{H}^1} = \|\mathbf{W} u_0\|_{\dot{H}^1} = \|u_0\|_{\dot{\mathcal{H}}_{\mathbf A}^1}. 
\]
Since we always have 
\[
 \left\|\mathbf{P}_{[0,b]} {\mathrm e}^{{\mathrm i} t \sqrt{\mathbf{A}}}  u_0 \right\|_{\dot{H}^1}^2 + \left\| \mathbf{P}_{(b,+\infty)} {\mathrm e}^{{\mathrm i} t \sqrt{\mathbf{A}}}  u_0 \right\|_{\dot{H}^1}^2 = \left\|{\mathrm e}^{{\mathrm i} t \sqrt{\mathbf{A}}}  u_0 \right\|_{\dot{H}^1}^2 \leq  \left\|{\mathrm e}^{{\mathrm i} t \sqrt{\mathbf{A}}}  u_0 \right\|_{\dot{H}_{\mathbf{A}}^1}^2 = \|u_0\|_{\dot{\mathcal{H}}_{\mathbf A}^1}^2. 
\]
Therefore
\[
 \left\|\mathbf{P}_{(b,+\infty)} {\mathrm e}^{{\mathrm i} t \sqrt{\mathbf{A}}}  u_0\right\|_{\dot{H}^1}^2  \rightarrow 0.
\]
It is clear that $\mathbf{P}_{(b,+\infty)} \mathbf{U}(t) \mathbf{W} u_0 = 0$. Collecting both the low frequency part and the high frequency part, we conclude that 
\[
 \lim_{t\rightarrow +\infty} \left\|\mathbf{U}(t) \mathbf{W} u_0 - {\mathrm e}^{{\mathrm i} t \sqrt{\mathbf{A}}}  u_0\right\|_{\dot{H}^1} = 0,
\]
which immediately gives \eqref{to prove H1} and finishes the proof. 
\end{proof}

\paragraph{Negative time direction} Next we consider the wave operator in the negative time direction. We take the conjugate on \eqref{half wave to prove m1} and obtain 
\[
  \lim_{t\rightarrow+\infty} \left\|\overline{\mathbf{U}}(t) \overline{\mathbf{W}} \overline{u_0} - {\mathrm e}^{-{\mathrm i} t \sqrt{\mathbf{A}}}  \overline{u_0} \right\|_{L^2(\Rm^+)} = 0. 
\]
Here 
\begin{align*}
 &\overline{\mathbf{U}}(t) f = \frac{2}{\pi} \int_0^\infty \sin kx \left({\mathrm e}^{-{\mathrm i}k t - {\mathrm i} P(k,t)} (\mathcal{F}_0 f)(k)\right) {\rm d} k; & 
 &\overline{\mathbf{W}} v = \mathcal{F}_0^{-1} \frac{\mathcal{F} v}{A(k)}.
\end{align*}
Similarly we have the $\dot{H}^1$ convergence. We summarize 
\[
 \overline{\mathbf{W}} = \mathop{{\rm s-}\lim}\limits_{t\rightarrow +\infty} \overline{\mathbf{U}}(t)^{-1} {\mathrm e}^{{-\mathrm i} t \sqrt{\mathbf{A}}} 
\]
is a unitary operator from $L^2(\Rm^+)$ to itself, and from $\dot{\mathcal{H}}_{\mathbf A}^1$ to $\dot{H}^1(\Rm^+)$. 

\subsection{Wave equations} 

Now we consider the wave equation 
\[
 \left\{\begin{array}{l} \partial_t^2 u + \mathbf{A} u = 0; \\ (u,u_t)|_{t=0} = (u_0,u_1)\in \dot{\mathcal{H}}_{\mathbf{A}}^1 \times L^2(\Rm^+). \end{array}\right.
\]
By functional calculus the corresponding solution can be given by 
\[
 (u, u_t) = \left(\cos (t \sqrt{\mathbf{A}}) u_0 + \frac{\sin (t \sqrt{\mathbf{A}})}{\sqrt{\mathbf{A}}} u_1, -\sqrt{\mathbf{A}} \sin (t \sqrt{\mathbf{A}}) u_0 + \cos (t \sqrt{\mathbf{A}}) u_1\right). 
\]
Using the $L^2$ and $\dot{H}^1$ level approximations given in the previous subsection (as $t\rightarrow +\infty$)
\begin{align*}
 &\mathbf{U}(t) \mathbf{W} \approx {\mathrm e}^{{\mathrm i} t \sqrt{\mathbf{A}}};& &\overline{\mathbf{U}}(t) \overline{\mathbf{W}} \approx {\mathrm e}^{-{\mathrm i} t \sqrt{\mathbf{A}}};
\end{align*}
we deduce 
\begin{equation} \label{approximation of uut}
 \begin{pmatrix} u\\ u_t \end{pmatrix} = \mathcal{F}_0^{-1} \begin{pmatrix} \cos \xi(k,t) & k^{-1} \sin \xi(k,t)\\ -k \sin \xi(k,t) & \cos \xi(k,t) \end{pmatrix}  \frac{\mathcal{F}}{|A(k)|} \begin{pmatrix} u_0 \\ u_1 \end{pmatrix} + o(1);
\end{equation}
with 
\[
 \xi(k,t) = kt + P(k,t) + \arg A(k). 
\]
Here $o(1)$ represents a remainder term whose $\dot{H}^1 \times L^2(\Rm^+)$ norm converges to zero as $t\rightarrow +\infty$. Next we define 
\begin{align} \label{def of vec U}
 &\vec{\mathbf{U}}(t) = \mathcal{F}_0^{-1} \begin{pmatrix} \cos \eta (k,t) & k^{-1} \sin \eta(k,t) \\ -k \sin \eta(k,t) & \cos \eta(k,t) \end{pmatrix} \mathcal{F}_0; &
 & \eta(k,t) = kt + P(k,t). 
\end{align}
be a family of unitary operators in the energy space $\dot{H}^1\times L^2$, whose inverses can be given by 
\begin{align*}
 \vec{\mathbf{U}}(t)^{-1}  = \mathcal{F}_0^{-1} \begin{pmatrix} \cos \eta (k,t) & -k^{-1} \sin \eta(k,t) \\ k \sin \eta(k,t) & \cos \eta(k,t) \end{pmatrix} \mathcal{F}_0.  
\end{align*}
 In view of \eqref{approximation of uut}, it is not difficult to see 
\[
 \lim_{t\rightarrow +\infty} \vec{\mathbf{U}}(t)^{-1} \begin{pmatrix} u\\ u_t \end{pmatrix} = \mathcal{F}_0^{-1} \begin{pmatrix} \cos (\arg A(k)) & k^{-1} \sin (\arg A(k)) \\ -k \sin (\arg A(k)) & \cos (\arg A(k)) \end{pmatrix} \frac{\mathcal{F}}{|A(k)|} \begin{pmatrix} u_0 \\ u_1 \end{pmatrix}.
\] 
Thus the wave operator $\vec{\mathbf{W}}$ exists and can be given explicitly in term of $A(k)$. In summary we have 
\begin{proposition} \label{modified wave operator type I general} 
 Assume that the potential $q(x)$ is a type I repulsive potential with decay rate $\beta>1/3$. Let $\vec{\mathbf{S}}_q (t)$ be the corresponding wave propagation operator of the wave equation 
 \[
  u_{tt} - u_{xx} +q(x) u = 0, \qquad (x,t) \in \Rm^+ \times \Rm
 \]
 and $\vec{\mathbf{U}}(t)$ be the unitary operators defined in \eqref{def of vec U}. Then the wave operator given by a strong limit in $\dot{H}^1 \times L^2$
 \[
  \vec{\mathbf{W}} \doteq \mathop{{\rm s-}\lim}\limits_{t\rightarrow +\infty} \vec{\mathbf{U}}(t)^{-1} \vec{\mathbf{S}}_q (t)
 \]
 is a well-defined unitary operator from the space $\dot{\mathcal{H}}_{\mathbf{A}}^1 \times L^2$ to $\dot{H}^1 \times L^2$. In addition, the wave operator can be given explicitly by the Fourier transforms and the functions $A(k)$ associated to the wave functions of the operator $-{\rm d}^2/{\rm d}x^2 + q(x)$:
 \[
  \vec{\mathbf{W}} =\mathcal{F}_0^{-1} \begin{pmatrix} \cos (\arg A(k)) & k^{-1} \sin (\arg A(k)) \\ -k \sin (\arg A(k)) & \cos (\arg A(k)) \end{pmatrix} \frac{\mathcal{F}}{|A(k)|}. 
 \]
\end{proposition} 

\subsection{Fast decaying case}

In this subsection, we consider type I repulsive potential $q(x)$ with a higher decay rate and show that we may define a simpler version of phase shift function such that the modified wave operator still exists. We consider two cases, i.e. decay rate $\beta > 1/2$ and $q(x) \in L^1(\Rm^+)$ respectively.  Before we start the discussion, we first give a technical lemma. 

\begin{lemma} \label{simple inter}
 Let $R_1 (k,t)$ and $R_2(k,t)$ be two real-valued functions. For any given finite interval $[a,b]\subset (0,+\infty)$, their difference $R_1(k,t)-R_2(k,t)$ converges to a real-valued function $\Delta R(k)$ uniformly for all $k\in [a,b]$ as $t\rightarrow +\infty$. Let $\vec{\mathbf{U}}_j(t)$ be unitary operators ($j=1,2$)
 \begin{align*}
  &\vec{\mathbf{U}}_j(t) =  \mathcal{F}_0^{-1} \begin{pmatrix} \cos \xi_j (k,t) & k^{-1} \sin \xi_j (k,t) \\ -k \sin \xi_j (k,t) & \cos \xi_j (k,t) \end{pmatrix} \mathcal{F}_0;& &\xi_j(k,t) = kt + R_j(k,t).
 \end{align*}
 Then we have the following strong limit
 \[
  \mathop{{\rm s-}\lim}\limits_{t\rightarrow +\infty} \vec{\mathbf{U}}_2(t)^{-1} \vec{\mathbf{U}}_1 (t) = \mathcal{F}_0^{-1} \begin{pmatrix} \cos \Delta R(k) & k^{-1} \sin \Delta R(k) \\ -k \sin \Delta R(k) & \cos \Delta R(k) \end{pmatrix} \mathcal{F}_0.
 \]
\end{lemma}
\begin{proof}
 A direct calculation shows that 
\[
 \vec{\mathbf{U}}_2(t)^{-1} \vec{\mathbf{U}}_1 (t) = \mathcal{F}_0^{-1} \begin{pmatrix} \cos (R_1(k,t)-R_2(k,t)) & k^{-1} \sin (R_1(k,t)-R_2(k,t)) \\ -k \sin (R_1(k,t)-R_2(k,t)) & \cos (R_1(k,t)-R_2(k,t)) \end{pmatrix} \mathcal{F}_0.
\]
If the Fourier transform of the data are compactly supported in $[a,b]\subset (0,+\infty)$ and smooth, the limit clearly exists by the uniform convergence assumption. The general case then follows from smooth approximation techniques and the fact that $\vec{\mathbf{U}}_j(t)$ are all unitary operators.   
\end{proof}

\paragraph{Intermediate case} We recall the definition of phase shift function 
\[
 P(k,t) = \frac{1}{2k} Q_1(t) - \frac{1}{4k^3} q(t) Q_1(t) + \frac{1}{8k^3} Q_2(t) + \frac{1}{16k^5} \int_0^\infty q^3(x) {\rm d} x
\]
We observe that if $q(x)$ is a type I repulsive potential with decay rate $\beta > 1/2$, then all terms except for the first one converges as $t\rightarrow +\infty$. Thus we redefine the phase shift function 
\[
 P_1(k,t) = \frac{1}{2k} Q_1(t) = \frac{1}{2k} \int_0^t q(s) {\rm d} s. 
\]
and the approximated wave propagation operator 
\begin{align} \label{def of vec U1}
 &\vec{\mathbf{U}}_1 (t) = \mathcal{F}_0^{-1} \begin{pmatrix} \cos \eta_1 (k,t) & k^{-1} \sin \eta_1 (k,t) \\ -k \sin \eta_1 (k,t) & \cos \eta_1 (k,t) \end{pmatrix} \mathcal{F}_0; &
 & \eta_1(k,t) = kt + \frac{1}{2k}\int_0^t q(s){\rm d} s. 
\end{align}
By Lemma \ref{simple inter} we have the strong limit 
\begin{align*}
 \mathop{{\rm s-}\lim}\limits_{t\rightarrow +\infty} \vec{\mathbf{U}}_1(t)^{-1} \vec{\mathbf{U}} (t) & = \mathcal{F}_0^{-1} \begin{pmatrix} \cos \theta_1(k) & k^{-1} \sin \theta_1(k) \\ -k \sin \theta_1(k) & \cos \theta_1(k) \end{pmatrix} \mathcal{F}_0;\\
 \theta_1(k) & = \frac{1}{8k^3} \int_0^\infty q^2(x) {\rm d} x + \frac{1}{16k^5} \int_0^\infty q^3(x) {\rm d} x. 
\end{align*}
Since we may write $\vec{\mathbf{U}}_1(t)^{-1} \vec{\mathbf{S}}_q (t)$ as a composition
\[
 \vec{\mathbf{U}}_1(t)^{-1} \vec{\mathbf{S}}_q (t) = \left(\vec{\mathbf{U}}_1(t)^{-1} \vec{\mathbf{U}} (t)\right) \left(\vec{\mathbf{U}} (t)^{-1} \vec{\mathbf{S}}_q (t)\right),
\]
we may combine these two strong limits and immediately obtain 
\begin{proposition}
 Assume that the potential $q(x)$ is a type I repulsive potential with decay rate $\beta>1/2$. Let $\vec{\mathbf{S}}_q (t)$ be the corresponding wave propagation operator of the wave equation 
 \[
  u_{tt} - u_{xx} +q(x) u = 0, \qquad (x,t) \in \Rm^+ \times \Rm
 \]
 and $\vec{\mathbf{U}}_1 (t)$ be the unitary operators defined in \eqref{def of vec U1}. Then the wave operator defined by the strong limit
 \[
  \vec{\mathbf{W}}_1 \doteq \mathop{{\rm s-}\lim}\limits_{t\rightarrow +\infty} \vec{\mathbf{U}}_1 (t)^{-1} \vec{\mathbf{S}}_q (t)
 \]
 is a well-defined unitary operator from $\dot{\mathcal{H}}_{\mathbf{A}}^1 \times L^2$ to $\dot{H}^1 \times L^2$.
 \end{proposition} 

\paragraph{High decay rate case} Now let us consider a type I repulsive potential $q(x)$ satisfying 
\[
 \int_0^\infty q(x) {\rm d} x < +\infty. 
\]
Please note that in this case the monotonicity of $q$ implies that $q(x)$ must satisfy $q(x) \lesssim x^{-1}$. In this case the phase shift function $P(k,t)$ converges as $t\rightarrow +\infty$. Thus we may choose the usual wave propagation operator 
\begin{align} \label{def of vec U0}
  \vec{\mathbf{U}}_0 (t) = \mathcal{F}_0^{-1} \begin{pmatrix} \cos kt & k^{-1} \sin kt \\ -k \sin kt & \cos kt \end{pmatrix} \mathcal{F}_0. 
\end{align}
A similar argument as above then shows 
\begin{align*}
 \mathop{{\rm s-}\lim}\limits_{t\rightarrow +\infty} \vec{\mathbf{U}}_0(t)^{-1} \vec{\mathbf{U}} (t) & = \mathcal{F}_0^{-1} \begin{pmatrix} \cos \theta_0(k) & k^{-1} \sin \theta_0(k) \\ -k \sin \theta_0(k) & \cos \theta_0(k) \end{pmatrix} \mathcal{F}_0;\\
 \theta_0 & = \frac{1}{2k} \int_0^\infty q(x) {\rm d} x +\frac{1}{8k^3} \int_0^\infty q^2(x) {\rm d} x + \frac{1}{16k^5} \int_0^\infty q^3(x) {\rm d} x. 
\end{align*}

\begin{proposition}
 Assume that the potential $q(x)\in L^1(\Rm^+)$ is a type I repulsive potential. Let $\vec{\mathbf{S}}_q (t)$ be the corresponding wave propagation operator of the wave equation 
 \[
  u_{tt} - u_{xx} +q(x) u = 0, \qquad (x,t) \in \Rm^+ \times \Rm
 \]
 and $\vec{\mathbf{U}}_0 (t)$ be the usual wave propagation operator. Then the wave operator defined by the strong limit
 \[
  \vec{\mathbf{W}}_0 \doteq \mathop{{\rm s-}\lim}\limits_{t\rightarrow +\infty} \vec{\mathbf{U}}_0 (t)^{-1} \vec{\mathbf{S}}_q (t)
 \]
must be a unitary operator from $\dot{\mathcal{H}}_{\mathbf{A}}^1 \times L^2$ to $\dot{H}^1 \times L^2$. This implies that for any finite-energy solution $u$ of the equation above, there exists a finite-energy solution $v$ to the classic wave equation $v_{tt} - v_{xx} = 0$ in $\Rm^+$ (with boundary condition $v(0,t)=0$) such that 
\[
 \lim_{t\rightarrow +\infty} \left\|(u(\cdot,t),u_t(\cdot,t)) - (v(\cdot,t),v_t(\cdot,t))\right\|_{\dot{H}^1\times L^2(\Rm^+)} = 0. 
\]
 \end{proposition} 

\subsection{Type II and III repulsive potentials}

Now we consider type II or III repulsive potentials, which may have a strong singularity near the zero, and complete the proof of Theorem \ref{thm main one-dimensional}. Please note that type I repulsive potential must be a type II repulsive potential as well. Let $q(x)$ be such a potential with decay rate $\beta>1/3$. We define a new potential $q^\ast$ by 
\[
 q^\ast (x) = \left\{\begin{array}{ll} q(x), & x\in [1,+\infty);\\ q(1)+(x-1) q'(1), & x\in [0,1]. \end{array}\right.
\]
Clearly $q^\ast(x)$ is a type I repulsive potential. Let $\vec{\mathbf{S}}_q(t)$ and $\vec{\mathbf{S}}_{q^\ast} (t)$ be the wave propagation operators of the wave equations $w_{tt} - w_{xx} + q(x) w = 0$ and $w_{tt} - w_{xx} + q^\ast(x) w = 0$, respectively. An application of Proposition \ref{modified wave operator type I general} gives the existence and the unitary property of the strong limit
\[
 \mathop{{\rm s-}\lim}\limits_{t\rightarrow +\infty} \mathcal{F}_0^{-1} \begin{pmatrix} \cos \eta^\ast (k,t) & k^{-1} \sin \eta^\ast (k,t) \\ -k \sin \eta^\ast (k,t) & \cos \eta^\ast (k,t) \end{pmatrix}  \mathcal{F}_0 \vec{\mathbf{S}}_{q^\ast} (t).
\]
Here 
\[
 \eta^\ast (k,t) = k t +\frac{1}{2k} \int_0^t q^\ast(s) {\rm d} s - \frac{1}{4k^3} q^\ast (t) \int_0^t q^\ast(s) {\rm d} s+ \frac{1}{8k^3} \int_0^t q^\ast (s)^2 {\rm d} s + \frac{1}{16k^5} \int_0^\infty q^\ast (s)^3 {\rm d} s. 
\]
Combining this with Lemma \ref{simple inter}, the following strong limit exists and must be a unitary operator from the energy space of $w_{tt} - w_{xx} + q^\ast(x) w = 0$ to $\dot{H}^1 \times L^2(\Rm^+)$:
\begin{equation} \label{q ast limit}
 \mathop{{\rm s-}\lim}\limits_{t\rightarrow +\infty} \mathcal{F}_0^{-1} \begin{pmatrix} \cos \eta (k,t) & k^{-1} \sin \eta (k,t) \\ -k \sin \eta (k,t) & \cos \eta (k,t) \end{pmatrix}  \mathcal{F}_0 \vec{\mathbf{S}}_{q^\ast} (t), 
\end{equation} 
with 
\[
 \eta(k,t) =kt+ \frac{1}{2k} \int_1^t q(s) {\rm d} s - \frac{1}{4k^3} q (t) \int_1^t q(s) {\rm d} s+ \frac{1}{8k^3} \int_1^t q^2 (s) {\rm d} s.
\]
Now we recall Proposition \ref{local difference A} to deduce that the following strong limit exists and is a unitary operator between the corresponding energy spaces
\[
  \mathop{{\rm s-}\lim}\limits_{t\rightarrow +\infty} \vec{\mathbf{S}}_{q^\ast} (t)^{-1} \vec{\mathbf{S}}_q (t).  
\]
By a composition of operators, we may combine this with \eqref{q ast limit}, as well as the fact the operators in \eqref{q ast limit} are uniform bounded to obtain the first conclusion of Theorem \ref{thm main one-dimensional}. When $q$ satisfies stronger decay assumption $\beta > 1/2$, the same conclusion holds if we substitute $\eta(k,t)$ by 
\[
 kt + \frac{1}{2k} \int_1^t q(s) {\rm d} s,
\] 
thanks to Lemma \ref{simple inter}. The case $q \in L^1(1,+\infty)$ can be dealt with in the same manner. 

\begin{remark} \label{family uniform app}
 Let $u_{tt} - u_{xx} + \mu x^{-2} u + q_0(x) u = 0$ be a family of wave equations on $\Rm^+$. Here $\mu \in \{0\}\cup [3/4,+\infty)$ is a parameter and $q_0(x)$ is a fixed type II repulsive potential with a decay rate $\beta > 1/3$. We use the notation $\vec{\mathbf{S}}_\mu (t)$ and $\dot{\mathcal{H}}_\mu^1$ for the corresponding wave propagation operator and Sobolev space.  If we choose a family of unitary operators in $\dot{H}^1(\Rm^+) \times L^2(\Rm^+)$ 
 \[
  \vec{\mathbf{U}}(t) = \mathcal{F}_0^{-1} \begin{pmatrix} \cos \eta (k,t) & k^{-1} \sin \eta (k,t) \\ -k \sin \eta (k,t) & \cos \eta (k,t) \end{pmatrix}  \mathcal{F}_0,
 \]
 independent of $\mu$ with 
 \[
  \eta(k,t) =kt+ \frac{1}{2k} \int_1^t q_0(s) {\rm d} s - \frac{1}{4k^3} q_0 (t) \int_1^t q_0(s) {\rm d} s+ \frac{1}{8k^3} \int_1^t q_0^2 (s) {\rm d} s;
 \]
 then the following strong limit exists and is a unitary operator from the energy space $\dot{\mathcal{H}}_{\mu}^1\times L^2$ to $\dot{H}^1\times L^2$ for each $\mu$
 \[
   \mathop{{\rm s-}\lim}\limits_{t\rightarrow +\infty} \vec{\mathbf{U}} (t)^{-1} \vec{\mathbf{S}}_\mu (t). 
 \]
 This follows from a combination of Theorem \ref{thm main one-dimensional} and Lemma \ref{simple inter}. As a result, given any finite-energy solution $u$ to a wave equation above, we may always find a pair of initial data $(v_0,v_1)\in \dot{H}^1\times L^2(\Rm^+)$ such that 
 \[
  \lim_{t\rightarrow +\infty} \left\|\begin{pmatrix} u\\ u_t \end{pmatrix} - \vec{\mathbf{U}} (t) \begin{pmatrix} v_0\\ v_1\end{pmatrix} \right\|_{\dot{H}^1\times L^2} = 0. 
 \]
 If $q_0$ comes with a decay rate $\beta > 1/2$, then the same conclusion holds if we substitute $\eta (k,t)$ by $k t + \frac{1}{2k} \int_1^t q_0(s) {\rm d} s$. 
\end{remark}

\section{Wave equation with potential in dimensions $d\geq 3$}

In this section we consider the asymptotic behaviour of finite-energy solutions to the wave equation 
\begin{equation}\label{wave equation high}
 \partial_t^2 u - \Delta u + V(x) u = 0, \qquad (x,t) \in \Rm^d \times \Rm,
\end{equation}
and prove Theorem \ref{thm main high-dimensional}. Here $d\geq 3$ and $V(x) = q(|x|)$ is a radially symmetric type II repulsive potential. The major tool is the following spherical harmonic decomposition. 

\subsection{Spherical harmonic decomposition} 

In this subsection we utilize spherical harmonic decomposition to transform the higher dimensional problem to a family of half-line problems. 

\paragraph{Harmonic polynomials} We start by introducing the harmonic polynomials. We recall that the eigenfunctions of the Laplace-Beltrami operator on $\mathbb{S}^{d-1}$ are exactly the homogeneous harmonic polynomials of the variables $x_1, x_2, \cdots, x_d$. Such a polynomial $\Phi$ of degree $\nu$ satisfies
\[
 -\Delta_{\mathbb{S}^{d-1}} \Phi = \nu(\nu+d-2) \Phi. 
\]
We choose a Hilbert basis $\{\Phi_j(\theta)\}_{j\geq 0}$ of the operator $-\Delta_{\mathbb{S}^{d-1}}$ on the sphere $\mathbb{S}^{d-1}$ and denote the degree of the harmonic polynomial $\Phi_j$ by $\nu_j$. In particular we assume $\nu_0 = 0$ and $\nu_j > 0$ if $j \geq 1$. For more details, please refer to M\"{u}ller \cite{sharmonics}. 

\paragraph{Spherical harmonic decomposition} Now we may decompose each function $u \in L^2(\Rm^d)$ orthogonally in the following form 
\[ 
 u(x) = u_0(|x|) \Phi_0(x/|x|) + u_1(|x|) \Phi_1(x/|x|) + u_2(|x|) \Phi_2(x/|x|) + \cdots. 
\]
Here $u_j(r)$ can be calculated by 
\[
 u_j(r) = \int_{\mathbb S^{d-1}} u(r\theta) \Phi_j(\theta) {\rm d} \theta. 
\]
By the spherical coordinates we have
\[
 \|u_j(|x|) \Phi_j (x/|x|)\|_{L^2(\Rm^d)} = \left(\int_0^\infty r^{d-1} |u_j(r)|^2 {\rm d} r\right)^{1/2}. 
\]
Thus we may let $w_j(r) = r^{\frac{d-1}{2}} u_j(r)$ and rewrite 
\begin{align*}
 &u(x) = \sum_{j=0}^\infty |x|^{-\frac{d-1}{2}} w_j(|x|) \Phi_j(x/|x|); & &\|u\|_{L^2(\Rm^d)} = \sum_{j=0}^\infty \|w_j\|_{L^2(\Rm^+)}^2.
\end{align*}
 As a result, we may decompose $L^2(\Rm^d)$ orthogonally as below
\[
 L^2(\Rm^d) = \mathcal{H}_0 \oplus \mathcal{H}_1 \oplus \mathcal{H}_2 \oplus \cdots 
\]
Each component of the space can be given by 
\[
 \mathcal{H}_j = \left\{ |x|^{-\frac{d-1}{2}} w(|x|) \Phi_j(x/|x|): w \in L^2(\Rm^+)\right\},
\]
which can be viewed as a copy of $L^2(\Rm^+)$ by the following isometric bijection 
\[
 |x|^{-\frac{d-1}{2}} w(|x|) \Phi_k(x/|x|) \leftrightarrow w(r).
\]
The orthogonal projection $\mathbf{P}_j: L^2(\Rm^d) \rightarrow \mathcal{H}_j$ is given by 
\[
 \mathbf{P}_j u =  \left(\int_{\mathbb S^{d-1}} u(|x|\theta) \Phi_j(\theta) {\rm d} \theta\right) \Phi_j(x/|x|). 
\]
If we incorporate the natural isometric bijection above, the orthogonal projection $\tilde{\mathbf{P}}_j : L^2(\Rm^d) \rightarrow L^2(\Rm^+)$ can be given by 
\[
 w_j (r) = (\tilde{\mathbf{P}}_j u)(r) = r^{\frac{d-1}{2}} \int_{\mathbb S^{d-1}} u(r\theta) \Phi_j(\theta) {\rm d} \theta. 
\]
Please note that $\nu_0 = 0$ means that $\Phi_0(\theta)$ is a constant, thus the first component of the decomposition is a radial function. 

\paragraph{Decomposition of self-adjoint operators} Assume that $V(x) = q(|x|)$ is a radial type II repulsive potential. We let $\mathbf{H} = -\Delta + V(x)$.  A straight-forward calculation shows that for suitable functions $w$ we have
\begin{align} 
 \mathbf{H} \left(|x|^{-\frac{d-1}{2}} w(|x|) \Phi_j(x/|x|)\right) & =  \left(-\partial_r^2 - \frac{d-1}{r} \partial_r - \frac{1}{r^2} \Delta_\theta + q(r)\right) \left(r^{-\frac{d-1}{2}} w(r) \Phi_j(\theta)\right) \nonumber\\
 & = r^{-\frac{d-1}{2}} \left(-w''(r) + \frac{d-1}{r} w'(r) -\frac{d^2-1}{4r^2} w(r) \right) \Phi_j(\theta) \nonumber \\
 & \qquad + r^{-\frac{d-1}{2}} \left(- \frac{d-1}{r} w'(r) +\frac{(d-1)^2}{2r^2} w(r) \right) \Phi_j(\theta) \nonumber \\
 & \qquad \qquad + r^{-\frac{d-1}{2}} \left(\frac{\nu_j(d-2+\nu_j)}{r^2}w(r)+q(r) w(r)\right) \Phi_j(\theta) \nonumber\\
 & = r^{-\frac{d-1}{2}} \left(-w''(r) + \frac{\mu_j}{r^2} w(r) + q(r) w(r) \right) \Phi_j(\theta). \label{HA inter}
\end{align}
Here the constant $\mu_j$ is defined by 
\[
 \mu_j = \frac{(d-1+2\nu_j)(d-3+2\nu_j)}{4} \in \{0\}\cup [3/4,+\infty). 
\]
Thus if we view $\mathcal{H}_j$ as a copy of $L^2(\Rm^+)$, then the restriction of $\mathbf{H}$ on it is exactly the self-adjoint operator $\mathbf{A}_{\mu_j} = -{\rm d}^2/{\rm d} r^2 + \mu_j r^{-2} + q(r)$. The spherical harmonic decomposition above actually gives a decomposition of the self-adjoint operator $\mathbf{H}$. In addition, we have the energy identity 
\begin{align}
 \int_{\Rm^d} \left(|\nabla u|^2 + V(|x|) |u|^2\right) {\rm d} x & = \langle u, \mathbf{H} u\rangle 
  = \sum_{j=0}^\infty \langle w_j, \mathbf{A}_{\mu_j} w_j \rangle  \nonumber \\
  & = \sum_{j=0}^\infty \int_0^\infty \left(|w'_j(r)|^2 + \frac{\mu_j}{r^2} |w_j(r)|^2 +  q(r) |w_j(r)|^2 \right) {\rm d} r. \label{energy identity sd}
\end{align}
In particular, if we choose $V(x) = q(|x|) = 0$, we have the identity 
\begin{equation}  \label{energy identity V0}
 \int_{\Rm^d} |\nabla u|^2 {\rm d} x =  \sum_{j=0}^\infty \int_0^\infty \left(|w'_j(r)|^2 + \frac{\mu_j}{r^2} |w_j(r)|^2 \right) {\rm d} r.
\end{equation}
\begin{remark}
Strictly speaking, when we talk about a self-adjoint operator, it is necessary to specify its domain in some way. In fact we may choose the domain of $\mathbf{H}$ to be 
\[
 D(\mathbf{H}) = \left\{\sum_{j=0}^\infty |x|^{-\frac{d-1}{2}} w_j(|x|) \Phi_j(x/|x|):  w_j \in D(\mathbf{A}_{\mu_j}), \sum_{j=0}^\infty \left(\|w_j\|_{L^2(\Rm^+)}^2 + \|\mathbf{A}_{\mu_j} w_j\|_{L^2}^2\right) < \infty \right\}.
\]
If the potential $q$ is small, it is not difficult to see the domain of $\mathbf{H}$ defined above is exactly the same as $D(-\Delta)$, as one may expect. For simplicity we do not want to give the details about the domain of $\mathbf{H}$ for a general type II repulsive potential. But we would like to mention that by the corresponding quadratic forms of the self-adjoint operators, the energy identity as given in \eqref{energy identity sd} always holds.
\end{remark}

\paragraph{Wave propagation} The orthogonal decomposition of self-adjoint operators immediately yields the corresponding decomposition of free waves. Let $u$ be a finite-energy solution to the wave equation $u_{tt} - \Delta u +q(|x|) u = 0$. Then there exist a family of finite-energy solutions $w_j$ to the wave equation $\partial_t^2 w_j - \partial_x^2 w_j + \mu_j w_j + q(x) w_j = 0$ on the half line $\Rm^+$ with zero boundary condition, such that the following orthogonal decomposition holds for any $t\in \Rm$
\[
 \begin{pmatrix} u(\cdot,t) \\ u_t(\cdot,t) \end{pmatrix} = \sum_{j=0}^\infty \begin{pmatrix} |x|^{-\frac{d-1}{2}} w_j (|x|,t) \Phi_j(x/|x|)\\
 |x|^{-\frac{d-1}{2}} \partial_t w_j(|x|,t) \Phi_j(x/|x|)\end{pmatrix}. 
\]
Here the orthogonality not only holds in the energy space $\dot{\mathcal{H}}_V^1(\Rm^d) \times L^2(\Rm^d)$, but also in $\dot{H}^1 \times L^2(\Rm^d)$. Please note that each term in the decomposition above is also a finite-energy solution to the wave equation $u_{tt} - \Delta u +q(|x|) u = 0$. Next we claim that the potential energy converges to zero as $t\rightarrow +\infty$, i.e. 
\begin{equation} \label{decay of potential energy} 
 \lim_{t\rightarrow +\infty} \int_{\Rm^d} q(|x|) |u(x,t)|^2 {\rm d} x = 0.
\end{equation}
By the orthogonal decomposition given above, it suffices to show that 
\[
 \lim_{t\rightarrow +\infty} \int_{\Rm^d} q(|x|) \left||x|^{-\frac{d-1}{2}} w_j (|x|,t) \Phi_j(x/|x|)\right|^2 {\rm d} x = \lim_{t\rightarrow +\infty} \int_{\Rm^+} q(r) \left|w_j (r,t)\right|^2 {\rm d} r = 0,
\]
which immediately follows from the inward/outward energy theory. 
\subsection{A radial example}

Let us consider a special example in 3-dimensional case with radially symmetric data. This radially symmetric assumption implies that the corresponding spherical harmonic decomposition comes with only one non-trivial term 
\[
 u (x,t) = |x|^{-1} w(|x|,t) \Phi_0(x/|x|) = \frac{1}{2\sqrt{\pi}} |x|^{-1} w(|x|,t). 
\]
Here $w$ is a finite energy solution of the one-dimensional wave equation 
\[
 w_{tt} - w_{rr} + q(r) w = 0. 
\]
A direct calculation shows that
\begin{align*}
 \hat{u}(\xi,t) = \frac{1}{(2\pi)^{3/2}} \int_{\Rm^3} {\mathrm e}^{-{\mathrm i} x \cdot \xi} u(x,t) {\rm d} x = \frac{1}{2^{5/2}\pi^2} \int_{\Rm^3} {\mathrm e}^{-{\mathrm i} x\cdot \xi} |x|^{-1} w(|x|,t) {\rm d} x. 
\end{align*}
We then utilize spherical coordinates and rewrite 
\begin{align*}
 \hat{u}(\xi,t) & = \frac{1}{2^{5/2}\pi^2} \int_0^\infty \int_0^\pi \int_0^{2\pi}{\mathrm e}^{-{\mathrm i}r|\xi| \cos \theta} r^{-1} w(r,t) \cdot r^2 \sin \theta {\rm d} \varphi {\rm d} \theta {\rm d} r \\
 & = \frac{1}{2^{3/2}\pi} \int_0^\infty \int_0^\pi {\mathrm e}^{-{\mathrm i}r|\xi| \cos \theta} w(r,t) r \sin \theta {\rm d} \theta {\rm d} r\\
 & = \frac{1}{2^{1/2}\pi} \int_0^\infty \frac{\sin r|\xi|}{r|\xi|} w(r,t) r {\rm d} r \\
 & = \frac{1}{2^{1/2}\pi} |\xi|^{-1} (\mathcal{F}_0 w(\cdot,t)) (|\xi|). 
\end{align*}
Combining this identity of Fourier transforms and the corresponding wave operator on the half-line, we may deduce that the following strong limit exists in the space $\dot{H}^1\times L^2$
\begin{equation*} 
 \mathop{{\rm s-}\lim}\limits_{t\rightarrow +\infty} \mathcal{F}_0^{-1} \begin{pmatrix} \cos \eta (|\xi|,t) & |\xi|^{-1} \sin \eta (|\xi|,t) \\ -|\xi| \sin \eta (|\xi|,t) & \cos \eta (|\xi|,t) \end{pmatrix}  \mathcal{F}_0 \begin{pmatrix} u\\ u_t\end{pmatrix}, 
\end{equation*} 
with 
\[
 \eta(|\xi|,t) =|\xi| t+ \frac{1}{2|\xi|} \int_1^t q(s) {\rm d} s - \frac{1}{4|\xi|^3} q (t) \int_1^t q(s) {\rm d} s+ \frac{1}{8|\xi|^3} \int_1^t q^2 (s) {\rm d} s.
\]
From this example we guess that the approximated wave propagation operator in higher dimensional case can be defined in the same manner as in the half-line case, with the same phase shift function $P(|\xi|,t)$. The details are given in the subsequent subsection. 

\subsection{The general case}

In this subsection we prove Theorem \ref{thm main high-dimensional}. More precisely, we give the modified wave operator in dimensions $d\geq 3$ for the wave equation 
\[
 \partial_t^2 u - \Delta u + q(|x|) u = 0
\]
under the assumption that $q$ is a type II repulsive potential satisfying one of the following: 
\begin{itemize} 
 \item[(a)] $q$ comes with a decay rate $\beta > 1/2$ at the infinity;  
 \item[(b1)] $q$ comes with a decay rate $\beta \in (1/3,1/2]$ at the infinity; there exists a large number $R>0$ such that $q\in \mathcal{C}^2([R,+\infty))$ and that $|q''(x)| \lesssim x^{-2}$ for $x\geq R$; 
 \item[(b2)] $q$ comes with a decay rate $\beta \in (1/3,1/2]$ at the infinity; there exists a large number $R>0$ such that $q\in \mathcal{C}^2([R,+\infty))$ and that $q''(x)>0$ for $x\geq R$.
\end{itemize} 
Please note that if $\beta\leq 1/2$, the decay assumption at the infinity has already given the inequality $|q'(x)|\lesssim x^{-1-\beta}$. We then define the phase shift function $P(k,t)$ and the approximated wave propagation operators $\vec{\mathbf{U}}(t)$ accordingly 
\[
 P(k,t) = \left\{\begin{array}{ll} \displaystyle \frac{1}{2k} \int_1^t q(s) {\rm d} s, & \hbox{case (a)}; \\
\displaystyle  \frac{1}{2k} \int_1^t q(s) {\rm d} s - \frac{1}{4k^3} q(t) \int_1^t q(s) {\rm d} s + \frac{1}{8k^3} \int_1^t q^2(s) {\rm d} s, & \hbox{case (b1) or (b2)}. \end{array}\right.
\]
\begin{align*}
  & \vec{\mathbf{U}}(t) = \mathcal{F}_0^{-1} \begin{pmatrix} \cos \eta(|\xi|,t) & |\xi|^{-1} \sin \eta(|\xi|,t) \\ -|\xi| \sin \eta(|\xi|,t) & \cos \eta(|\xi|,t) \end{pmatrix} \mathcal{F}_0; & & \eta(|\xi|,t) = |\xi| t + P(|\xi|,t). 
\end{align*}
As usual we define 
\[
 Q_1(t) = \int_1^t q(s) {\rm d} s. 
\]
We start by giving a lemma concerning an ordinary differential equation:
\begin{lemma} \label{lemma ex}
 Assume that $q$ is a repulsive potential and satisfies either (b1) or (b2) above. Fix $\xi_0>0$. Then there exists a large time $T_0 = T_0(q, \xi_0)$ such that if $v$ is a solution to the ordinary differential equation 
 \[
  v_{tt} + \left(\xi^2 + q(t) - \frac{q'(t) Q_1(t)}{2\xi^2 + 2q^2(t)}\right) v = h, \qquad t\geq T.
 \]
 with $\xi>\xi_0$ and $h\in L^1([T,+\infty))$, then there exists two constants $A, B$ such that 
  \begin{align*}
  v(t) & = A \cos \eta(\xi,t) + B\sin \eta(\xi,t) + \mathcal{E}_0(\xi,t); \\
  v_t(t) & = -A \xi \sin \eta(\xi,t) + B \xi \cos \eta(\xi,t) + \mathcal{E}_1(\xi,t)
 \end{align*}
 holds for $t>\max\{T,T_0\}$ with $\eta(\xi,t) = \xi t + P(\xi,t)$ and 
 \begin{align*}
  \left| \mathcal{E}_0(\xi,t)\right| & \lesssim \xi^{-1} \int_t^\infty |h(s)| {\rm d}s + \left(\xi^{-2} t^{-\beta} + \xi^{-5} t^{-(3\beta-1)}\right) (|A|+|B|); \\
  \left| \mathcal{E}_1(\xi,t)\right| & \lesssim \int_t^\infty |h(s)| {\rm d}s + \left(\xi^{-1} t^{-\beta} + \xi^{-4} t^{-(3\beta-1)}\right) (|A|+|B|).
 \end{align*}
The implicit constant in the inequality only depends on the potential $q$ and $\xi_0$. 
\end{lemma}
\begin{proof}
 A direct calculation shows that $\cos(\xi t + P(\xi,t))$ and $\sin (\xi t + P(\xi,t))$ are both solutions to the ordinary differential equation 
 \[
  \mathbf{H} w = w_{tt} - b(\xi,t) w_t + \left(\xi^2 + q(t) - \frac{q'(t)Q_1(t)}{2\xi^2 + 2q(t)} + c(\xi,t)\right) w = 0,
 \]
 with the coefficients 
 \begin{align*}
  b(\xi,t) & = \frac{P_{tt} (\xi,t)}{\xi + P_t (\xi,t)} = \frac{\displaystyle \frac{q'(t)}{2\xi} - \frac{1}{4\xi^3}q''(t) Q_1(t) - \frac{1}{2\xi^3}q'(t) q(t)}{\displaystyle \xi + \frac{q(t)}{2\xi} - \frac{1}{4\xi^3}q'(t) Q_1(t) - \frac{1}{8\xi^3} q^2(t)}; \\
  c(\xi,t) & = - \frac{q'(t) q(t) Q_1(t)}{2\xi^2(\xi^2+q(t))} - \frac{q(t)}{\xi} \left(\frac{1}{4\xi^3} q'(t) Q_1(t) + \frac{1}{8\xi^3} q^2(t)\right) + \left(\frac{1}{4\xi^3} q'(t) Q_1(t) + \frac{1}{8\xi^3} q^2(t)\right)^2;
 \end{align*}
 and Wronskian 
\[
 W(\xi,t) = \begin{vmatrix} \cos(\xi t + P(\xi,t)) & \sin(\xi t + P(\xi,t)) \\ -(\xi + P_t(\xi,t)) \sin (\xi t + P(\xi,t)) & (\xi + P_t(\xi,t)) \cos(\xi t + P(\xi,t)) \end{vmatrix} =  \xi + P_t(\xi,t).
\] 
 The coefficients $b(\xi,t)$, $c(\xi,t)$ and Wronskian satisfy
 \begin{align}
  \int_R^\infty |b(\xi,t)| {\rm d} t & \lesssim \xi^{-2} R^{-\beta}, & & R>R_0; \label{ex decay b}\\
  \int_R^\infty |c(\xi,t)| {\rm d} t & \lesssim \xi^{-4} R^{-(3\beta-1)}, & &R>R_0; \label{ex decay c}\\
  \left|W(\xi,t) - \xi\right| = |P_t(\xi,t)|& \leq \xi/2, & & t>R_0. \label{ex constant W}
 \end{align}
 Here the implicit constant and $R_0$ can be chosen uniformly for all $\xi > \xi_0$. Since $v$ solves the differential equation 
 \begin{equation} \label{ex trans equation} 
  \mathbf{H} v = c(\xi,t) v - b(\xi,t) v_t + h,
 \end{equation}
 we expect that $(v,v_t)$ solves the integral equation ($\eta(\xi,t) = \xi t + P(\xi,t)$)
 \begin{align*}
    v(t) & = \cos \eta(\xi,t) \left(A + \int_t^\infty \sin \eta(\xi,s) \frac{c(\xi,s) v(s) - b(\xi,s) v_t(s)  + h(s)}{W(\xi,s)} {\rm d} s\right)\\
  &\quad  + \sin \eta(\xi,t) \left(B-\int_t^\infty \cos \eta(\xi,s) \frac{c(\xi,s) v(s) - b(\xi,s) v_t(s)  + h(s)}{W(\xi,s)}  {\rm d} s\right);\\
  v_t(t) & =  -(\xi + P_t(\xi,t))\sin \eta(\xi,t) \left(A + \int_t^\infty \sin \eta(\xi,s) \frac{c(\xi,s) v(s) - b(\xi,s) v_t(s)  + h(s)}{W(\xi,s)} {\rm d} s\right)\\
  & + (\xi + P_t(\xi,t))\cos \eta(\xi,s) \left(B-\int_t^\infty \cos \eta(\xi,s) \frac{c(\xi,s) v(s) - b(\xi,s) v_t(s)  + h(s)}{W(\xi,s)}  {\rm d} s\right);
  \end{align*}
 or equivalently, $(v,v_t)$ is a fixed point of the map $\mathbf{T}$ defined by the right hand side of integral equation above, where $A,B$ are both constants. Now we consider the space (the parameter $R_1> R_0$ to be determined later)
 \[
  X= \left\{(f,g)\in C([R_1,+\infty))^2: \sup_{t\geq R_1} \left(|\xi f(t)| + |g(t)|\right) < +\infty \right\}
 \]
 with norm 
 \[
  \|(f,g)\|_X = \sup_{t\geq R_1} \left(|\xi f(t)| + |g(t)|\right). 
 \]
 In view of \eqref{ex decay b}, \eqref{ex decay c} and \eqref{ex constant W}, we have
 \begin{align*}
  &\|\mathbf{T} (f,g)\|_{X} \leq 3 \xi (|A|+|B|) + 10 \int_{R_1}^\infty |h(s)| {\rm d} s+C_1\left(\xi^{-2} R_1^{-\beta} + \xi^{-5} R_1^{-(3\beta-1)}\right) \|(f,g)\|_{X};\\
  &\|\mathbf{T} (f_1,g_1) - \mathbf{T}(f_2,g_2)\|_X \leq C_1\left(\xi^{-2} R_1^{-\beta} + \xi^{-5} R_1^{-(3\beta-1)}\right) \|(f_1,g_1)-(f_2,g_2)\|_{X}. 
 \end{align*}
 Here $C_1$ depends on $q$ and $\xi_0$ only. We may choose a large number $T_0 = T_0(q,\xi_0)\geq R_0$  such that
 \[
  C_1\left(\xi_0^{-2} T_0^{-\beta} + \xi_0^{-5} T_0^{-(3\beta-1)}\right) < 1/2.
 \]
 This implies that $\mathbf{T}$ is a contraction map in $X$ if we choose $R_1 = \max\{T,T_0\}$. The unique fixed point $(v,v_t)$ is exactly a solution to \eqref{ex trans equation}, thus to the original differential equation. In addition, it is not difficult to see from the inequalities above that 
 \[
  \sup_{t\geq r} \left(|\xi v(t)| + |v_t(t)|\right) \leq 6\xi (|A|+|B|) + 20 \int_r^\infty |h(s)| {\rm d} s, \qquad r\geq R_1.
 \]
 The error estimate for $\mathcal{E}_0, \mathcal{E}_1$ then follows this upper bound and the integral equation. Finally it is not difficult to see this family of solutions with two parameters $A$, $B$ covers all possible solutions to the differential equation, which finishes the proof. 
 \end{proof}

\begin{remark} \label{remark ex} 
 A similar result holds for repulsive potentials with a decay rate $\beta > 1/2$. In this case we consider a solution to the differential equation 
 \[
  v_{tt} + (\xi^2 + q(t)) v = h, \qquad t\geq T
 \]
 with $h \in L^1([T,+\infty))$. Fix $\xi_0>0$. Then there exists  a number $T_0 = T_0 (q,\xi_0)$ such that any solution to the differential equation above satisfies 
  \begin{align*}
  &\left\{\begin{array}{l} v(t) = A \cos \eta(\xi,t) + B \sin \eta(\xi,t) + \mathcal{E}_0(t) ; \\ v_t(t) = -A \xi \sin \eta(\xi,t) + B \xi \cos(\xi,t) + \mathcal{E}_1(t);    \end{array}\right.& &\eta(\xi, t) = \xi t + P(\xi,t) = \xi t + \frac{1}{2\xi} \int_1^t q(s) {\rm d}s.
 \end{align*}
 with error term estimates 
 \begin{align*}
  & \left\{\begin{array}{l} \left|\mathcal{E}_0(t)\right| \lesssim \xi^{-1} \int_t^\infty |f(s)| {\rm d} s + (\xi^{-2}t^{-\beta} + \xi^{-3} t^{1-2\beta}) (|A|+|B|); \\ 
  \left|\mathcal{E}_1(t)\right| \lesssim  \int_t^\infty |f(s)| {\rm d} s + (\xi^{-1}t^{-\beta} + \xi^{-2} t^{1-2\beta}) (|A|+|B|). \end{array} \right. & & \forall  t > \max\{T_0,T\}. 
 \end{align*}
 Here the implicit constant only depends on the potential $q$. 
\end{remark}

\begin{lemma} \label{lemma bridge}
Let $q(x)$ be a type II repulsive potential satisfying the assumption (a), (b1) or (b2) given at the beginning of this subsection.  Assume that $\hat{w}_0$ and $\hat{w}_1$ are smooth functions of $k$ supported in an interval $[a,b]\subset (0,+\infty)$. In addition, $\Phi$ is a homogeneous harmonic polynomial with order $\nu$. Let
 \begin{align}
  w(r, t) & = \mathcal{F}_0^{-1} \left(\cos \eta(k,t) \hat{w}_0 + \sin \eta(k,t) \hat{w}_1\right), & & \eta(k,t) = kt + P(k,t); \\
  u(r\theta,t) & = r^{-\frac{d-1}{2}} w(r,t) \Phi(\theta), & & (r,\theta) \in \Rm^+ \times \mathbb{S}^{d-1}. 
 \end{align}
 Then the following limit exists in $\dot{H}^1(\Rm^d) \times L^2(\Rm^d)$ as $t\rightarrow +\infty$
 \[
  \lim_{t\rightarrow +\infty} \vec{\mathbf{U}}(t)^{-1} \begin{pmatrix} u(t)\\ u_t(t) \end{pmatrix}. 
 \]
 Here $\vec{\mathbf{U}}(t)$ and $P(k,t)$ are defined at the beginning of this subsection. 
\end{lemma} 
\begin{proof}
 We focus on the case when $q$ satisfies (b1) or (b2), which is more complicated. Before the conclusion of the proof we shall explain how to deal with case (a).  Without loss of generality we assume $\Phi =\Phi_j$ is one of the Hilbert basis $\{\Phi_j\}$. A straight-forward calculation shows that $w$ is $\mathcal{C}^2$ function in $\Rm^+\times \Rm$. We let $\hat{w}$ be the Fourier transform of $w$. We define an auxiliary function 
\[
 g(r,t) = \mathcal{F}_0^{-1} \left(\frac{\hat{w}}{k^2 +q(t)}\right) = \mathcal{F}_0^{-1} \left(\frac{\cos \eta(k,t) \hat{w}_0 + \sin \eta(k,t) \hat{w}_1}{k^2+q(t)}\right).
\]
It is clear that 
\[
 (-{\rm d}^2/{\rm d} r^2 + q(t)) g(r,t) = w(r,t). 
\]
In addition, a basic calculation of Fourier transforms shows that 
\begin{align*}
 w_{tt} - w_{rr} + q(t) w & = \mathcal{F}_0^{-1} \left[(\hat{w}_1 \cos \eta - \hat{w}_0 \sin \eta)\left(\frac{q'(t)}{2k} - \frac{1}{4k^3}q''(t) Q_1(t) -\frac{1}{2k^3}q'(t) q(t)\right)\right]\\
 & \qquad + \mathcal{F}_0^{-1} \left[\hat{w}\left(\frac{1}{2k^2}q'(t) Q_1(t) + \frac{1}{4k^4} q'(t) q(t) Q_1(t) + \frac{1}{8k^4} q^3(t)\right)\right]\\
 & \qquad \quad - \mathcal{F}_0^{-1} \left[\hat{w}\left(\frac{1}{16k^6} q'(t) q^2(t) Q_1(t) + \frac{1}{16k^6}q'(t)^2 Q_1^2(t) + \frac{1}{64k^6} q^4(x)\right)\right]. 
\end{align*}
We also have
\[
 \frac{1}{2k^2}q'(t) Q_1(t) \hat{w} = \frac{1}{2} q'(t) Q_1(t) \mathcal{F}_0 g + \frac{q'(t) q(t) Q_1(t)}{2k^2(k^2+q(t))} \hat{w}. 
\]
Combining these two identities we have
\begin{equation} \label{ex w equation}
 w_{tt} - w_{rr} + q(t) w = \frac{1}{2} q'(t) Q_1(t) g + f(r,t)
\end{equation}
with $f \in L^1 L^2([T,+\infty) \times \Rm^+)$. Here $T$ represents a large time. Furthermore, the following decay estimate holds for large time $t>T$ by an integration by parts
\begin{equation} \label{ex non cone decay}
 |w(r,t)| + |w_r(r,t)| + |w_t(r,t)| + |g(r,t)| + |g_r(r,t)| + |g_t(r,t)| \lesssim |r-t|^{-1}, \qquad |r-t|>t/6. 
\end{equation} 
In addition, we may also obtain the following uniform upper bound by Fourier transform
\begin{align}\label{ex upper bound by Fourier} 
 &\sup_{t>0} \|w\|_{H^1(\Rm^+)} < +\infty;& & \sup_{t>0} \|w_t\|_{L^2\cap \dot{H}^{-1} (\Rm^+)} < +\infty. 
\end{align}
In order to avoid the technical difficulty near $r=0$, we introduce a centre cut-off version of $w$ defined by
 \[
  \tilde{w} (r,t) = \phi(r/t) w(r,t). 
 \]
 Here $\phi: \Rm^+ \rightarrow [0,1]$ is a smooth cut-off function satisfying 
 \[
  \phi(x) = \left\{\begin{array}{ll} 0, & x\leq 1/3; \\ 1, & x\geq 2/3.  \end{array}\right.
 \]
 We have
\begin{align*}
  \tilde{w}_{tt} - \tilde{w}_{rr} + q(t) \tilde{w}  = & \phi(r/t) \left(w_{tt} - w_{rr} + q(t) w\right) -\frac{2r}{t^2} \phi'(r/t) w_t + \frac{r^2}{t^4} \phi''(r/t)w +2 \frac{r}{t^3} \phi'(r/t)w\\
  & \qquad - \frac{2}{t} \phi'(r/t) w_r - \frac{1}{t^2} \phi''(r/t) w\\
  = &  \frac{1}{2} q'(t) Q_1(t) \tilde{g} + \bar{f},
 \end{align*}
 with $\tilde{g} = \phi(r/t) g$ and $\bar{f} \in L^1 L^2([T,+\infty)\times \Rm^+)$ by \eqref{ex non cone decay} and \eqref{ex upper bound by Fourier}. It also follows from these estimates that
 \begin{equation} \label{ex approximation tilde w}
  \lim_{t\rightarrow +\infty} \|w-\tilde{w}\|_{\dot{H}^1(\Rm^+)} + \|w_t - \tilde{w}_t\|_{L^2(\Rm^+)} = 0. 
 \end{equation}
 Hardy's inequality immediately gives 
 \begin{equation} \label{ex approximation tilde w 2}
  \lim_{t\rightarrow +\infty} \int_0^\infty \frac{|w-\tilde{w}|^2}{r^2} {\rm d} r = 0. 
 \end{equation} Similarly we have
 \begin{align*}
  -\tilde{g}_{rr} + q(t) \tilde{g} & = \phi(r/t) \left(-g_{rr} + q(t) g\right) - \frac{2}{t} \phi'(r/t) g_r - \frac{1}{t^2} \phi''(r/t) g\\
  & = \tilde{w} + \bar{h}.
 \end{align*}
Here $\bar{h} \in L^1 L^2([T,+\infty)\times \Rm^+)$. By \eqref{ex upper bound by Fourier} and the support of $\tilde{w}$, it is easy to see 
 \[
  \left\|\frac{\tilde{w}(\cdot,t)}{r^2}\right\|_{L^2(\Rm^+)} \lesssim t^{-2} \quad \Rightarrow \quad \frac{\tilde{w}}{r^2} \in L^1([T,\infty); L^2(\Rm^+)). 
 \]
Similarly we have $\tilde{g} /r^2 \in L^1 L^2([T,+\infty)\times \Rm^+)$. Therefore we may write 
\begin{align}
 \tilde{w}_{tt} - \tilde{w}_{rr} + q(t) \tilde{w}  + \mu_{d+2\nu} \frac{\tilde{w}}{r^2} & = \frac{1}{2} q'(t) Q_1(t) \tilde{g} + \tilde{f}, & & \tilde{f} \in L^1 L^2([T,+\infty)\times \Rm^+); \label{ex tilde w equation}\\
 - \tilde{g}_{rr} + q(t) \tilde{g} + \mu_{d+2\nu} \frac{\tilde{g}}{r^2} & = \tilde{w} + \tilde{h}, & &  \tilde{h} \in L^1 L^2([T,+\infty)\times \Rm^+).
\end{align}
Here the constant $\lambda_{d+2\nu}$ is defined by 
 \[
  \mu_{d+2\nu} = \frac{(d+2\nu-1)(d+2\nu-3)}{4}. 
 \]
Now we define $\tilde{u}(r\theta,t)  = r^{-\frac{d-1}{2}} \tilde{w}(r,t) \Phi(\theta)$. It follows from \eqref{energy identity V0}, \eqref{ex approximation tilde w} and the Hardy inequality that 
\[
 \lim_{t\rightarrow +\infty} \left\|(\tilde{u}, \tilde{u}_t) - (u,u_t)\right\|_{\dot{H}^1\times L^2(\Rm^d)} = 0. 
\]
Therefore in order to finish the proof of the lemma, it suffices to show the following limit exists in $\dot{H}^1\times L^2(\Rm^d)$
\[
  \lim_{t\rightarrow +\infty} \vec{\mathbf{U}}(t)^{-1} \begin{pmatrix} \tilde{u}(t)\\ \tilde{u}_t(t) \end{pmatrix}. 
 \]
We let $g_1 (r\theta,t) = r^{-\frac{d-1}{2}} \tilde{g}(r,t) \Phi(\theta)$ and define $f_1, h_1 \in L^1 L^2([T,+\infty)\times \Rm^d)$ accordingly. A straight-forward computation similar to \eqref{HA inter} yields 
 \begin{align*}
  (\partial_t^2 - \Delta + q(t)) \tilde{u} & = \left(\partial_t^2 - \partial_r^2 - \frac{d-1}{r}\partial_r - \frac{1}{r^2} \Delta_\theta + q(t)\right)\left(r^{-\frac{d-1}{2}} \tilde{w}(r,t) \Phi(\theta)\right) \\
  & = r^{-\frac{d-1}{2}} \left(\tilde{w}_{tt} - \tilde{w}_{rr} + \mu_{d+2\nu} \frac{\tilde{w}}{r^2} + q(t) \tilde{w}\right) \Phi(\theta)\\
  & = \frac{1}{2} q'(t) Q_1(t) g_1 + f_1. 
 \end{align*}
Similarly we have 
\[ 
 (-\Delta + q(t)) g_1 = \tilde{u} + h_1. 
\]
Next we observe that 
\begin{equation} \label{ex uniform boundedness}
 \sup_{t>T} \left\|(\tilde{u},\tilde{u}_t)\right\|_{H^1(\Rm^d)\times (L^2\cap \dot{H}^{-1}(\Rm^d))} < +\infty. 
\end{equation}
This is a combination of the following facts:
\begin{itemize}
 \item The solutions $(w,w_t)$ are uniformly bounded in $H^1(\Rm^+) \times (L^2\cap \dot{H}^{-1})(\Rm^+)$; 
 \item The pair
 \[
   (\tilde{w}, \tilde{w}_t) = \left(\rho(r/t) w(r,t), \rho(r/t) w_t (r,t) - \frac{r}{t^2}\rho'(r/t) w(r, t)\right)
 \]
 is still uniformly bounded in $H^1(\Rm^+) \times L^2(\Rm^+)$ for large time $t\gg 1$ by a direct calculation. A duality argument also shows that $\tilde{w}_t$ is uniformly bounded in $\dot{H}^{-1}$. 
 \item Similarly the map $\rho \rightarrow r^{-\frac{d-1}{2}} \rho(r) \Phi(\theta)$ are bounded operators from $\dot{H}^s(\Rm^+)$ to $\dot{H}^s(\Rm^d)$ for $s=0,\pm 1$. The situation $s=0,1$ can be verified by a direct calculation and the case $s=-1$ follows a duality argument.
\end{itemize}
We may apply the Fourier transform with respect to $x$ and obtain ($\hat{u} = \mathcal{F}_0 \tilde{u}$) 
\begin{align*}
 &\hat{u}_{tt} + (|\xi|^2+q(t)) \hat{u} = \frac{1}{2} q'(t) Q_1(t) \hat{g}_1 + \hat{f}_1;& & (|\xi|^2 + q(t)) \hat{g}_1 = \hat{u} + \hat{h}_1. 
\end{align*}
Inserting the second equation into the first one, we have 
\begin{equation} \label{ex hat u}
 \hat{u}_{tt} + (|\xi|^2+q(t)) \hat{u} = \frac{q'(t) Q_1(t)}{2(|\xi|^2 + q(t))} \hat{u} + \frac{q'(t)Q_1(t)}{2(|\xi|^2 + q(t))}\hat{h}_1 + \hat{f}_1.
\end{equation}
Here the error term satisfies ($R>0$ is an arbitrary constant)
\[
  \hat{h}_2 \doteq \frac{q'(t)Q_1(t)\hat{h}_1}{2(|\xi|^2 + q(t))} + \hat{f}_1 \in L^1 L^2([T,+\infty)\times \{\xi: |\xi|>R\}) \hookrightarrow L^2(\{\xi: |\xi|>R\}; L^1([T,+\infty))). 
\]
By Lemma \ref{lemma ex}, there exists functions $A(\xi)$ and $B(\xi)$ such that for almost everywhere $\xi$ we have 
  \begin{align}
  \hat{u}(\xi, t) & = A(\xi) \cos \eta(|\xi|,t) + B(\xi)\sin \eta(|\xi|,t) + \mathcal{E}_0(\xi,t); \label{ex hat u}\\
  \hat{u}_t(\xi,t) & = -A(\xi) |\xi| \sin \eta(|\xi|,t) + B(\xi) |\xi| \cos \eta(|\xi|,t) + \mathcal{E}_1(\xi,t); \label{ex hat u t}
 \end{align}
 for $t>T(\xi)$ with 
 \begin{align*}
  \left| \mathcal{E}_0(\xi,t)\right| & \lesssim |\xi|^{-1} \int_t^\infty |\hat{h}_2 (s)| {\rm d}s + \left(|\xi|^{-2} t^{-\beta} + |\xi|^{-5} t^{-(3\beta-1)}\right) (|A(\xi)|+|B(\xi)|); \\
  \left| \mathcal{E}_1(\xi,t)\right| & \lesssim \int_t^\infty |\hat{h}_2 (s)| {\rm d}s + \left(|\xi|^{-1} t^{-\beta} + |\xi|^{-4} t^{-(3\beta-1)}\right) (|A(\xi)|+|B(\xi)|).
 \end{align*}
 Here $T(\xi)$ and the implicit constant in the inequalities above can be chosen independent of $\xi$ as long as $|\xi|$ is bounded from below by a positive constant. We may solve $A(\xi)$ and $B(\xi)$ from \eqref{ex hat u} and \eqref{ex hat u t}:
 \begin{align*}
  |\xi| A(\xi) & = \lim_{t\rightarrow +\infty} \left[|\xi|\hat{u}(\xi,t) \cos \eta(|\xi|,t) - \hat{u}_t(\xi,t) \sin \eta(|\xi|,t)\right]; \\
  |\xi| B(\xi) & = \lim_{t\rightarrow +\infty} \left[|\xi| \hat{u}(\xi,t) \sin \eta(|\xi|,t) + \hat{u}_t (\xi,t) \cos \eta(|\xi|,t)\right].
 \end{align*}
 This implies that $A(\xi)$ and $B(\xi)$ are both measurable functions of $\xi$. Fatou's lemma also guarantees that 
 \begin{align*}
  \int_{\Rm^d} \left(|\xi|^2 |A(\xi)|^2 + |\xi|^2 |B(\xi)|^2\right) {\rm d} \xi & \lesssim \liminf_{t\rightarrow \infty} \int_{\Rm^d} \left(|\xi|^2 |\hat{u}(\xi,t)|^2 + |\hat{u}_t(\xi,t)|^2 \right) {\rm d} \xi\\
  & \lesssim \liminf_{t\rightarrow \infty} \int_{\Rm^d} \left(|\nabla \tilde{u}(x,t)|^2 + |\tilde{u}_t(x,t)|^2 \right) {\rm d} x < +\infty. 
 \end{align*}
 A combination of this with the fact $\hat{h}_2 \in L^2(\{\xi\in \Rm^d: |\xi|>R\}; L^1([T,+\infty)))$ and the upper bound estimate of $(\mathcal{E}_0, \mathcal{E}_1)$ implies that 
 \begin{equation} \label{ex e12 decay}
  \lim_{t\rightarrow +\infty} \int_{|\xi|>R} \left(|\xi|^2 |\mathcal{E}_0(\xi,t)|^2 + |\mathcal{E}_1(\xi,t)|^2 \right) {\rm d} \xi  = 0, \qquad \forall R > 0.  
 \end{equation}
 A direct calculation gives 
 \begin{align}
  \vec{\mathbf{U}}(t)^{-1} \begin{pmatrix} \tilde{u}(t) \\ \tilde{u}_t (t) \end{pmatrix} & = \mathcal{F}_0^{-1} \begin{pmatrix} \cos \eta & -|\xi|^{-1} \sin \eta\\ |\xi| \sin \eta & \cos \eta \end{pmatrix} \begin{pmatrix} \hat{u} (\xi,t) \\ \hat{u}_t (\xi,t) \end{pmatrix}  \nonumber\\
  & = \mathcal{F}_0^{-1} \left[\begin{pmatrix} A(\xi)\\ |\xi| B(\xi)\end{pmatrix} + \begin{pmatrix} \mathcal{E}_0(\xi,t) \cos \eta - |\xi|^{-1} \mathcal{E}_1(\xi,t) \sin \eta \\ |\xi|\mathcal{E}_0(\xi,t) \sin \eta + \mathcal{E}_1(\xi,t) \cos \eta \end{pmatrix}\right]. \label{ex difference}
 \end{align}
 Next we let 
 \[
  \begin{pmatrix} \tilde{v}_0 \\ \tilde{v}_1\end{pmatrix} = \mathcal{F}_0^{-1} \begin{pmatrix} A(\xi)\\ |\xi| B(\xi)\end{pmatrix}  \in \dot{H}^1(\Rm^d) \times L^2(\Rm^d). 
 \]
 A combination of \eqref{ex e12 decay} and \eqref{ex difference} yields that 
 \[
  \lim_{t\rightarrow +\infty} \left\|\mathbf{P}_{|\xi|>R}\left(\vec{\mathbf{U}}(t)^{-1} \begin{pmatrix} \tilde{u}(t) \\ \tilde{u}_t (t) \end{pmatrix}  - \begin{pmatrix} \tilde{v}_0 \\ \tilde{v}_1\end{pmatrix}\right)\right\|_{\dot{H}^1\times L^2} = 0, \qquad \forall R>0.
 \]
 Here $\mathbf{P}$ is the regular frequency cut-off function. This gives the high frequency convergence. To deal with the low frequency part, we recall the uniform upper bound \eqref{ex uniform boundedness} and deduce  
 \[
  \lim_{R\rightarrow 0^+} \left(\limsup_{t\rightarrow +\infty} \left\|\mathbf{P}_{|\xi|\leq R} \begin{pmatrix} \tilde{u}(t) \\ \tilde{u}_t (t) \end{pmatrix} \right\|_{\dot{H}^1 \times L^2}\right) = 0.
 \]
 A combination of the high and low frequency parts finishes the proof for the potentials satisfying (b1) or (b2). Next we explain how to deal with case (a). In this case $w$ satisfies the following equation
\begin{align*}
 w_{tt} - w_{rr} + q(t) w & = \mathcal{F}_0^{-1} \left(-\frac{q(t)^2}{4k^2}\hat{w} - \frac{q'(t)}{2k} \sin \eta(k,t) \hat{w}_0 + \frac{q'(t)}{2k} \cos\eta(k,t) \hat{w}_1\right) \nonumber \\
  & \qquad \in L_t^1([T,+\infty);L_x^2(\Rm^+)). 
\end{align*}
The same centre cut-off gives $(\tilde{w}, \tilde{w}_t)$, which satisfies 
\[
 \lim_{t\rightarrow +\infty} \|(w,w_t) - (\tilde{w},\tilde{w}_t)\|_{\dot{H}^1(\Rm^+) \times L^2(\Rm^+)} = 0
\]
and 
\[
 \tilde{w}_{tt} - \tilde{w}_{rr} + \mu_{d+2\nu} \frac{\tilde{w}}{r^2} + q(t) \tilde{w} \in L_t^1([T,+\infty);L_x^2(\Rm^+)). 
\]
Next we define $(\tilde{u},\tilde{u}_t) = r^{-\frac{d-1}{2}} (\tilde{w}(r,t), \tilde{w}_t(r,t)) \Phi(\theta)$, which satisfies 
 \[
  (\partial_t^2 - \Delta + q(t)) \tilde{u} \in L^1 L^2([T,+\infty)\times \Rm^d). 
 \]
 Again it suffice to show the following limit exists in $\dot{H}^1 \times L^2(\Rm^d)$: 
\[
  \lim_{t\rightarrow +\infty} \vec{\mathbf{U}}(t)^{-1} \begin{pmatrix} \tilde{u}(t)\\ \tilde{u}_t(t) \end{pmatrix}. 
 \]
 Applying the Fourier transform and making use of Remark \ref{remark ex}, we may give the limit $(\tilde{v}_0,\tilde{v}_1)\in \dot{H}^1\times L^2$ in term of $A(\xi)$ and $B(\xi)$ and show the convergence in the high frequency part. Finally the low frequency part can be dealt with by the uniform boundedness of $(\tilde{u}, \tilde{u}_t)$ in $L^2 \times \dot{H}^{-1}$.  
\end{proof}

\begin{proof}[Proof of Theorem \ref{thm main high-dimensional}]
 Now we are ready to prove the main theorem in higher dimensional case $d\geq 3$. We start by considering the strong limit $\vec{\mathbf{U}}(t)^{-1} \vec{\mathbf{S}}_V(t)$. Given a free wave $u$, We start by recalling the orthogonal decomposition in the energy space $\mathcal{H}_V^1(\Rm^d) \times L^2(\Rm^d)$ given by the spherical harmonic decomposition: 
 \[
  \begin{pmatrix} u(\cdot,t) \\ u_t(\cdot,t) \end{pmatrix} = \sum_{j=0}^\infty \begin{pmatrix} r^{-\frac{d-1}{2}} w^j (r,t) \Phi_j (\theta) \\ r^{-\frac{d-1}{2}} w_t^j (r,t) \Phi_j(\theta) \end{pmatrix}.  
 \]
 Here $w^j$ is a finite-energy free wave of the wave equation $w_{tt} - w_{rr} + \mu_{d+2\nu_j} r^{-2} w + q w =0$. Since the series above converges uniformly in the energy space for all $t\in \Rm$, it suffices to show that the following limit in $\dot{H}^1\times L^2$ exists for any $j \geq 0$:
 \begin{equation} \label{to prove Rmd}
  \lim_{t\rightarrow +\infty} \vec{\mathbf{U}}(t)^{-1}  \begin{pmatrix} r^{-\frac{d-1}{2}} w^j (r,t) \Phi_j (\theta) \\ r^{-\frac{d-1}{2}} w_t^j (r,t) \Phi_j(\theta) \end{pmatrix}. 
 \end{equation}
 According to Remark \ref{family uniform app}, there exists a pair of initial data $(v_0^j, v_1^j) \in \dot{H}^1(\Rm^+) \times L^2(\Rm^+)$ such that 
 \[
  \lim_{t\rightarrow +\infty} \left\|\mathcal{F}_0^{-1} \begin{pmatrix} \cos \eta(k,t) & k^{-1} \sin \eta(k,t) \\ -k\sin \eta(k,t) & \cos \eta(k,t) \end{pmatrix} \mathcal{F}_0 \begin{pmatrix} v_0^j\\ v_1^j \end{pmatrix} - \begin{pmatrix} w^j\\ w_t^j \end{pmatrix} \right\|_{\dot{H}^1\times L^2(\Rm^+)} = 0.
 \]
 Here $\eta(k,t) = kt + P(k,t)$. By the standard cut-off and smooth approximation techniques, we may find a sequence of smooth functions with compact support $\{(\hat{v}_{0,\ell}, \hat{v}_{1,\ell})\}_{\ell \geq 1}$, such that 
 \[
  \lim_{\ell \rightarrow +\infty} \int_0^\infty \left(k^2 |\hat{v}_{0,\ell} (k) - (\mathcal{F}_0 v_0^j)(k)|^2 + |\hat{v}_{1,\ell}(k) - (\mathcal{F}_0 v_1^j)(k)|^2 \right) {\rm d} k  = 0. 
 \]
 Therefore we have 
  \begin{equation} \label{long limit v}
  \lim_{\ell \rightarrow +\infty} \left(\limsup_{t\rightarrow +\infty} \left\|\mathcal{F}_0^{-1} \begin{pmatrix} \cos \eta(k,t) & k^{-1} \sin \eta(k,t) \\ -k\sin \eta(k,t) & \cos \eta(k,t) \end{pmatrix} \begin{pmatrix} \hat{v}_{0,\ell} \\ \hat{v}_{1,\ell} \end{pmatrix} - \begin{pmatrix} w^j\\ w_t^j \end{pmatrix} \right\|_{\dot{H}^1\times L^2(\Rm^+)}\right) = 0.
 \end{equation} 
 Let $v^\ell = \mathcal{F}_0^{-1} \left[\hat{v}_{0,\ell}\cos \eta(k,t) + (k^{-1} \hat{v}_{1,\ell}) \sin \eta(k,t)\right]$. A basic calculation shows that 
 \[
  \lim_{t\rightarrow +\infty} \left\|v_t^\ell - \mathcal{F}_0^{-1} \left[\hat{v}_{0,\ell} (-k \sin \eta(k,t)) +  \hat{v}_{1,\ell} \cos \eta(k,t)\right]\right\|_{L^2(\Rm^+)} = 0.
 \]
 Inserting this into \eqref{long limit v}, we obtain 
 \[
  \lim_{\ell \rightarrow +\infty} \left(\limsup_{t\rightarrow +\infty} \left\| \begin{pmatrix} v^\ell  \\ v_t^\ell \end{pmatrix}  - \begin{pmatrix} w^j\\ w_t^j \end{pmatrix} \right\|_{\dot{H}^1\times L^2(\Rm^+)}\right) = 0.
 \]
 Since the map $g(r) \rightarrow r^{-\frac{d-1}{2}} g(r) \Phi_j (\theta)$ is bounded from $L^2(\Rm^+)$ to $L^2(\Rm^d)$ and from $\dot{H}^1(\Rm^+)$ to $\dot{H}^1(\Rm^d)$, we deduce 
 \begin{equation} \label{limit convergence v ell}
  \lim_{\ell \rightarrow +\infty} \left(\limsup_{t\rightarrow +\infty} \left\| \begin{pmatrix} r^{-\frac{d-1}{2}} v^\ell(r,t) \Phi_j (\theta)  \\ r^{-\frac{d-1}{2}} v_t^\ell(r,t) \Phi_j (\theta) \end{pmatrix}  - \begin{pmatrix} r^{-\frac{d-1}{2}}w^j (r,t) \Phi_j(\theta)\\ r^{-\frac{d-1}{2}} w_t^j (r,t) \Phi_j(\theta)\end{pmatrix} \right\|_{\dot{H}^1\times L^2(\Rm^d)}\right) = 0.
 \end{equation}
Lemma \ref{lemma bridge} guarantees that the following limits exist in $\dot{H}^1 \times L^2$:
\[
 \lim_{t\rightarrow +\infty} \vec{\mathbf{U}}(t)^{-1} \begin{pmatrix} r^{-\frac{d-1}{2}} v^\ell(r,t) \Phi_j (\theta)  \\ r^{-\frac{d-1}{2}} v_t^\ell(r,t) \Phi_j (\theta) \end{pmatrix}, \qquad \forall \ell \in \mathbb N. 
\]
A combination of this with \eqref{limit convergence v ell} then verifies the existence of the limit \eqref{to prove Rmd}, thus the existence of the strong limit $\vec{\mathbf{U}}(t)^{-1} \vec{\mathbf{S}}_V(t)$ as $t\rightarrow +\infty$. The unitary property follows the fact that $\vec{\mathbf{U}}(t)$ is unitary in $\dot{H}^1\times L^2$ and the limit 
\[
 \lim_{t\rightarrow +\infty} \|\vec{\mathbf{S}}_V(t) (u_0,u_1)\|_{\dot{H}^1\times L^2(\Rm^d)} = \|(u_0,u_1)\|_{\dot{\mathcal{H}}_V^1 \times L^2}.
\]
Here we utilize the decay of potential energy given in \eqref{decay of potential energy}. The remaining task is to show that we may substitute $\vec{\mathbf{U}}(t)$ by the classic wave propagation operator $\vec{\mathbf{S}}_0(t)$ if $q \in L^1(1,+\infty)$. In this case, our repulsive assumption implies that $q$ comes with a decay rate $\beta \geq 1$. Thus the strong limit $\vec{\mathbf{U}}(t)^{-1} \vec{\mathbf{S}}_V(t)$ exists as $t\rightarrow +\infty$ and is a unitary operator from $\dot{\mathcal{H}}_V^1 \times L^2$ to $\dot{H}^1 \times L^2$, where 
\[
 \vec{\mathbf{U}}(t) = \mathcal{F}_0^{-1} \begin{pmatrix} \cos \left(|\xi| t + \frac{1}{2|\xi|} \int_1^t q(s) {\rm d} s\right) & |\xi|^{-1} \sin  \left(|\xi| t + \frac{1}{2|\xi|} \int_1^t q(s) {\rm d} s\right) \\ -|\xi| \sin  \left(|\xi| t + \frac{1}{2|\xi|} \int_1^t q(s) {\rm d} s\right) & \cos  \left(|\xi| t + \frac{1}{2|\xi|} \int_1^t q(s) {\rm d} s\right)\end{pmatrix} \mathcal{F}_0.
\]
A similar argument to Lemma \ref{simple inter} shows that the strong limit $\vec{\mathbf{S}}_0 (t)^{-1} \vec{\mathbf{U}}(t)$ exists and is unitary operator on $\dot{H}^1 \times L^2$. Finally we observe that 
\[
 \vec{\mathbf{S}}_0 (t)^{-1} \vec{\mathbf{S}}_V (t) = \left(\vec{\mathbf{S}}_0 (t)^{-1} \vec{\mathbf{U}}(t)\right)\left( \vec{\mathbf{U}} (t)^{-1} \vec{\mathbf{S}}_V (t)\right)
\]
and finish the proof.
\end{proof}

\section{Dispersion Rate} 

In this section we consider the energy distribution property of solutions to wave equation with repulsive potential. As usual, for a potential $q(x)$ we define 
\[
 Q_1(t) = \int_1^t q(x) {\rm d} x. 
\]

\begin{lemma}[dispersion rate in half line] \label{main dispersion}
 Assume that $w$ is a finite-energy solution to the wave equation $u_{tt} - u_{xx} + q(x) u = 0$ on the half line, where $q(x)$ is a type I, II or III repulsive potential with a decay rate $\beta>1/3$ and satisfies $Q_1 (t) \rightarrow +\infty$ as $t\rightarrow +\infty$. Then given any $\varepsilon > 0$, there exist two constants $0<c_1<c_2<+\infty$ such that the following inequalities 
 \begin{align*}
  &\int_0^{t-c_2 Q_1(t)} e(x,t) {\rm d} x < \varepsilon;& &\int_{t-c_1 Q_1(t)}^\infty e(x,t) {\rm d} x < \varepsilon
 \end{align*}
 hold for sufficiently large time $t\gg 1$. Here $e(x,t)$ is the energy density function 
 \[
  e(x,t) = |w_x(x,t)|^2 + |w_t(x,t)|^2 + q(x) |w(x,t)|^2. 
 \]
 In addition, if $\ell(t)$ is a function satisfying $\ell(t)/Q_1 (t) \rightarrow 0$ as $t\rightarrow +\infty$, then 
 \[
  \lim_{t\rightarrow +\infty} \left(\sup_{r>0} \int_r^{r+\ell(t)} e(x,t) {\rm d} x\right) = 0. 
 \]
 \end{lemma}
 \begin{proof}
  By our inward/outward energy theory (see Corollary \ref{cor limit asymptotic}), we always have 
  \[
   \lim_{t\rightarrow +\infty} \int_0^\infty q(x) |w(x,t)|^2 {\rm d} x = 0.
  \]
  Thus it suffices to consider the corresponding integral of $e_0(x,t) = |w_x(x,t)|^2 + |w_t(x,t)|^2$. By Theorem \ref{thm main one-dimensional} and smooth approximation techniques we only need to the prove the corresponding result for 
  \[
   \tilde{e}_0 (x,t) = |\tilde{w}_0 (x,t)|^2 + |\tilde{w}_1(x,t)|^2
  \]
  with
  \begin{align*}
   \tilde{w}_0 (x,t) & = \frac{2}{\pi}\int_0^\infty \cos (kx) \left[kf(k)\cos (kt + P(k,t)) + g(k) \sin (kt+P(k,t))\right] {\rm d} k;\\
   \tilde{w}_1 (x,t) & = \frac{2}{\pi}\int_0^\infty \sin (kx) \left[-kf(k)\sin (kt + P(k,t)) + g(k) \cos (kt+P(k,t))\right] {\rm d} k.
  \end{align*}
  Here $f$ and $g$ are smooth functions supported in an interval $[a,b] \subset (0,+\infty)$ and the phase shift function $P(k,t)$ is defined by 
  \[
  P(k,t) = \frac{1}{2k} \int_1^t q(s) {\rm d} s - \frac{1}{4k^3} q(t) \int_1^t q(s) {\rm d} s + \frac{1}{8k^3} \int_1^t q^2(s) {\rm d} s.
  \]
  We claim that given any function $f \in \mathcal{C}_0^\infty([a,b])$, the following four functions 
  \[
   \tilde{w}_{\pm,\pm} (x,t) = \int_0^\infty f(k) \exp\left[\pm {\mathrm i}\left(\pm kx + kt + P(k,t)\right)\right] {\rm d} k
  \] 
  satisfy the following inequalities when $t$ is sufficiently large:
  \begin{align*}
   |\tilde{w}_{\pm,\pm} (x,t)| & \lesssim \left(t - \frac{1}{2a^2} Q_1(t) -x \right)^{-1}, & &\forall x < t - \frac{1}{a^2} Q_1(t);\\
   |\tilde{w}_{\pm,\pm} (x,t)| & \lesssim \left(x - t + \frac{1}{2b^2} Q_1(t) \right)^{-1}, & &\forall x > t - \frac{1}{4b^2} Q_1(t);\\
   |\tilde{w}_{\pm,\pm} (x,t)| & \lesssim [Q_1(t)]^{-1/2}, & & \forall x > 0. 
  \end{align*}
  It is not difficult to see that these inequalities imply that the desired result holds. The remaining task is to prove these inequalities. By taking the complex conjugate we see that it suffices to prove the inequality for 
 \[
   \tilde{w}_{+,\pm} (x,t) = \int_0^\infty f(k) \exp\left[ {\mathrm i}\left(\pm kx + kt + P(k,t)\right)\right] {\rm d} k. 
 \]
 We observe that 
  \[
   \frac{\rm d}{{\rm d} k} \exp\left[{\mathrm i}\left(\pm kx + kt + P(k,t)\right)\right] = {\mathrm i}\left(\pm x + t + \partial_k P(k,t)\right)\exp\left[{\mathrm i}\left(\pm kx + kt + P(k,t)\right)\right].
  \]
  Here if $k\in [a,b]$ and $t$ is sufficiently large, the derivatives of $P(k,t)$ always satisfy 
  \begin{align*}
   \partial_k P(k,t) & = -\frac{1}{2k^2} \int_1^t q(s) {\rm d} s + \frac{3}{4k^4} q(t) \int_1^t q(s) {\rm d} s - \frac{3}{8k^4} \int_1^t q^2(s) {\rm d} s\\
   & \in \left[-\frac{3}{4a^2} Q_1(t), -\frac{3}{8b^2} Q_1(t)\right]; \\
   \left|\partial_k^j P(k,t)\right| & \simeq Q_1(t), \qquad j\geq 1.
  \end{align*}
  An integration by parts immediately shows that $|\tilde{w}_{+,+} (x,t)| \lesssim (x+t)^{-1}$ for all $x>0$ and verifies the first two inequalities for $|\tilde{w}_{+,-} (x,t)|$. Finally we prove the third inequality for $|\tilde{w}_{+,-}(x,t)|$. This follows an argument of stationary phase. By the first two inequalities, without loss of generality, we may assume that 
  \[
   t - \frac{1}{a^2} Q_1(t) \leq x \leq t - \frac{1}{4b^2} Q_1(t). 
  \]
 Next we observe that if $t\gg 1$, then $\partial_k^2 P(k,t) \gtrsim Q_1(t)$ for $k \in [a/2,2b]$ and 
 \begin{align*}
  & - x + t + \partial_k P(a/2,t) < 0; & &- x + t + \partial_k P(2b,t) > 0. 
 \end{align*}
 Therefore there exists exactly one number $k_0= k_0(x,t)\in (a/2,2b)$, such that $-x + t + \partial_k P(k_0,t) = 0$. We then let 
 \[
  J(x,t) = \left(k_0 - Q_1(t)^{-1/2}, k_0 + Q_1(t)^{-1/2}\right),
 \]
 and write 
 \begin{align*}
  \tilde{w}_{+,-} (x,t) & = \int_{J(x,t)} f(k) \exp\left[{\mathrm i}\left(- kx + kt + P(k,t)\right)\right] {\rm d} k \\
  & \qquad + \int_{[a,b]\setminus J(x,t)} f(k) \exp\left[{\mathrm i}\left(- kx + kt + P(k,t)\right)\right] {\rm d} k = J_1 + J_2.
 \end{align*}
 Clearly the first integral $J_1$ is dominated by $Q_1(t)^{-1/2}$; while for $k \in [a,b] \setminus J(x,t)$ we have 
 \[
  \left|- x + t + \partial_k P(k,t)\right| \geq \left(\inf_{k'\in [a/2,2b]} \partial_k^2 P(k',t)\right) Q_1(t)^{-1/2} \gtrsim Q_1(t)^{1/2}. 
 \]
 An integration by parts then shows that 
 \begin{align*}
  |J_2| & = \left|\int_{[a,b]\setminus J(x,t)} \frac{f(k)}{-x + t + \partial_k P(k,t)} {\rm d} \exp\left[{\mathrm i}\left(- kx + kt + P(k,t)\right)\right] \right|\\
   & \leq O\left(Q_1(t)^{-1/2}\right) + \left| \int_{[a,b]\setminus J(x,t)} f(k) \frac{\partial_k^2 P(k,t)}{(-x + t + \partial_k P(k,t))^2} \exp\left[{\mathrm i}\left(- kx + kt + P(k,t)\right)\right] {\rm d} k \right|\\
  & \leq O\left(Q_1(t)^{-1/2}\right) + (\sup |f|) \int_{[a,b]\setminus J(x,t)} \frac{\partial_k^2 P(k,t)}{(-x + t + \partial_k P(k,t))^2} {\rm d} k \\
  & \leq O\left(Q_1(t)^{-1/2}\right) . 
 \end{align*}
 Combining the upper bounds of $J_1$ and $J_2$, we finish the proof.
 \end{proof}

\begin{proof}[Proof of Corollary \ref{main dispersion high}]
 By the spherical harmonic decomposition, it suffices to prove the corresponding result for 
 \[
 \begin{pmatrix} u(\cdot,t) \\ u_t(\cdot,t) \end{pmatrix} = \sum_{j=0}^N \begin{pmatrix} r^{-\frac{d-1}{2}} w^j (r,t) \Phi_j (\theta) \\ r^{-\frac{d-1}{2}} w_t^j (r,t) \Phi_j(\theta) \end{pmatrix}.
 \]
 Here $w^j(r,t)$ is a finite-energy solution to the wave equation 
 \[
  w_{tt} - w_{rr} + \left(q(r) + \frac{\lambda_{d+2\nu_j}}{r^2}\right) w = 0. 
 \]
 Indeed those solutions are dense in the energy space $\dot{\mathcal{H}}_V^1 \times L^2(\Rm^d)$.  Since the sum above is finite, it follows that we only need to prove the corresponding result for a single term 
 \[
   \begin{pmatrix} u(\cdot,t) \\ u_t(\cdot,t) \end{pmatrix} = \begin{pmatrix} r^{-\frac{d-1}{2}} w^j (r,t) \Phi_j (\theta) \\ r^{-\frac{d-1}{2}} w_t^j (r,t) \Phi_j(\theta) \end{pmatrix}.
 \]
 By smooth approximation technique, we may also assume $w^j(r,t) \in \mathcal{C}^2(\Rm^+\times \Rm)$ by Corollary \ref{smoothness away from t}. Next we recall \eqref{HA inter} and make a direct calculation 
 \begin{align*}
  & \int_{a<|x|<b} \left(|\nabla u(x,t)|^2 + |u_t(x,t)|^2 + q(|x|) |u(x,t)|^2 \right) {\rm d} x \\
  = & \int_{a<|x|<b} \left(-\Delta u + q(x) u\right) \bar{u} {\rm d} x + \int_{|x|=b} \bar{u} u_r {\rm d} S - \int_{|x|=a} \bar{u} u_r {\rm d} S + \int_a^b |w^j_t (r,t)|^2 {\rm d} r \\
  = & \int_a^b \left(-w_{rr}^j + q(r) w^j + \frac{\lambda_{d+2\nu_j}}{r^2} w^j \right)\overline{w^j} {\rm d} r+ \left[r^\frac{d-1}{2} \frac{\partial}{\partial r} \left(r^{-\frac{d-1}{2}} w^j\right) \overline{w^j}\right]_{r=a}^b + \int_a^b |w^j_t (r,t)|^2 {\rm d} r \\
  = & \int_a^b \left(|w_r^j|^2 + q(r) |w^j|^2 + \frac{\lambda_{d+2\nu_j}}{r^2} |w^j|^2 + |w_t^j|^2 \right) {\rm d} r - \frac{d-1}{2b} |w^j(b,t)|^2 + \frac{d-1}{2a} |w^j(a,t)|^2.
 \end{align*} 
 In particular, if $a=0$ or $b=+\infty$, then we may simply ignore the boundary value, since we always have 
 \[
  \lim_{r\rightarrow 0^+} \frac{|w(r)|^2}{r} = \lim_{r\rightarrow \infty} \frac{|w(r)|^2}{r} = 0, \qquad \forall w \in \dot{H}^1(\Rm^+). 
 \]
 Combining the identity above with the asymptotic behaviour
 \[
  \lim_{t\rightarrow +\infty} \sup_{r>0} \frac{|w^j(r,t)|^2}{r} =0, 
 \]
 as given in Proposition \ref{asymptotic b}, we obtain for any functions $0\leq a(t) < b(t) \leq \infty$ that 
 \begin{align*}
  \lim_{t\rightarrow +\infty} \bigg[ & \int_{a(t)<|x|<b(t)} \left(|\nabla u(x,t)|^2 + |u_t(x,t)|^2 + q(|x|) |u(x,t)|^2 \right) {\rm d} x \\
  & -  \int_{a(t)}^{b(t)} \left(|w_r^j(r,t)|^2 + q(r) |w^j(r,t)|^2 + \frac{\lambda_{d+2\nu_j}}{r^2} |w^j(r,t)|^2 + |w_t^j(r,t)|^2 \right) {\rm d} r \bigg] = 0. 
 \end{align*}
 The conclusion of Theorem \ref{main dispersion high} then follows from the corresponding result for $w^j(r,t)$ given in Lemma \ref{main dispersion}.
\end{proof}

\section*{Appendix} 

In this section we verify that the operator $\mathbf{A} = -{\rm d}^2/{\rm d}x^2 + q_0(x) + \lambda x^{-2}$ is indeed a self-adjoint operator in $L^2(\Rm^+)$ and prove the inequalities given in Remark \ref{C1 continuity}. Here the actual domain $D(\mathbf{A})$ and the assumption of $q_0(x)$ are given at the beginning of Section \ref{sec: inward}. 

An integration by parts immediately gives the symmetry property of $\mathbf{A}$. It suffices to shows that $\mathbf{A}^*\subseteq \mathbf{A}$. Let $v \in D(\mathbf{A}^*)$, then the following identity holds for any $u \in C_0^\infty(0,+\infty)$
\begin{equation*}
 \int_0^\infty u \overline{\mathbf{A}^\ast v} {\rm d} x = \int_0^\infty (-u''(x)+q(x) u) \bar{v} {\rm d} x \; \Rightarrow \; 
 \int_0^\infty u \left(\overline{\mathbf{A}^\ast v - q v}\right) {\rm d} x = \int_0^\infty (-u''(x)) \bar{v} {\rm d} x. 
\end{equation*}
This implies that $v$ comes with a second-order weak derivation $v''(x) = qv - \mathbf{A}^\ast v\in L_{\rm loc}^2(0,+\infty)$. Since we are working in the one-dimensional case, it immediately follows that $v \in AC^2(\Rm^+)$, i.e. $v$ and $v'$ are both absolutely continuous in any compact interval $[a,b]\subset \Rm^+$ with 
\[
 \mathbf{A}^\ast v = -v''(x) + q(x) v \in L^2(\Rm^+). 
\]
Since $v\in L^2(\Rm^+)$ and $q\in L^\infty(1,+\infty)$, we immediately obtain $v''(x) \in L^2(1,+\infty)$. The remaining task is to investigate the behaviour of $v$ near zero and verify $v \in D(\mathbf{A})$. Without loss of generality we assume that $\kappa$ is slightly smaller than $2$. We consider two cases separately:

\paragraph{Case I} If $\lambda = 0$, then we may write 
\begin{align*}
 v(x) &= v(1) - \int_x^1 v'(y) {\rm d} y = v(1) - \int_x^1 \left(v'(1) - \int_y^1 v''(s) {\rm d} s\right) {\rm d} y\\
 & = v(1) +(x-1)v'(1) + \int_x^1 (s-x) v''(s) {\rm d} s\\
 & = v(1) +(x-1)v'(1) + \int_x^1 (s-x) \left(q(s) v(s) - (\mathbf{A}^\ast v)(s)\right) {\rm d} s. 
\end{align*}
Thus we have 
\begin{align}
 |v(x)| \leq |v(1)|+|v'(1)|+\|\mathbf{A}^\ast v\|_{L^2(0,1)} + C \int_x^1 s^{1-\kappa} |v(s)| {\rm d} s. \label{iteration one}
\end{align}
By Cauchy-Schwarz we have $|v(x)| \leq C_1 x^{3/2-\kappa}$ for $x\in (0,1)$. Inserting this into \eqref{iteration one} yields $|v(x)| \leq C_2 x^{7/2-2\kappa}$. Please note that $7/2-2\kappa > 3/2-\kappa$. We may iterate this argument and finally obtain 
\[
 \sup_{x\in (0,1)} |v(x)| \leq C_3. 
\]
By the identity 
\begin{equation} \label{iteration two}
 v'(x) = v'(1) - \int_x^1 v''(s) {\rm d} s = v'(1) + \int_x^1 \left((\mathbf{A}^\ast v)(s) - q(s) v(s)\right) {\rm d} s, 
\end{equation}
we have $|v'(x)| \leq C_4 x^{1-\kappa}$ for $x\in (0,1)$, which implies that the limit 
\[
 v(0) \doteq \lim_{x\rightarrow 0^+} v(x) 
\]
is well-defined. Next we show $v(0)=0$. Let 
\begin{align*}
 &g(x) = \sum_{j=0}^N g_j(x); & &g_0(x) = x;& & g_{j+1}(x) = \int_0^x \int_0^y g_j(s) q(s) {\rm d} s{\rm d} y. 
 \end{align*}
An induction shows that $|g_{j}(x)| \lesssim |x|^{1+j(2-\kappa)}$ and 
\[
 -g''(x) + q(x) g(x) = q(x) g_N(x), \qquad x>0.
\]
We may choose a sufficiently large $N$ such that $g\in \mathcal{C}^2(\Rm^+)$ satisfies
\begin{align*}
 &-g'' + q g\in L^2(0,10);& &g, g'\in L^2(0,10);& &|g(x)| \lesssim |x|, \; x\in (0,1);& &\lim_{x\rightarrow 0^+} g'(x) = 1.&
\end{align*}
We may cut off the part near infinity smooth and construct a function $\tilde{g} \in D(\mathbf{A})$ with $\tilde{g}(x) = g(x)$ for $x\in (0,1)$. By the definition of $\mathbf{A}^\ast$ we have the inner product identity $\langle \tilde{g}, \mathbf{A}^\ast v\rangle = \langle \mathbf{A} \tilde{g}, v \rangle$. An integration by parts shows that 
\[
 \lim_{x \rightarrow 0^+} \left(\tilde{g}'(x) \overline{v(x)} - \tilde{g}(x) \overline{v'(x)}\right) = 0,
\]
which implies that $v(0)=0$. Therefore we may utilize the upper bound of $v'(x)$ and obtain that $|v(x)| \leq C_5 x^{2-\kappa}$ for $x\in (0,1)$. Inserting this into \eqref{iteration two} yields that $|v'(x)| \lesssim x^{3-2\kappa}$ thus $|v(x)|\leq C_6 x^{4-2\kappa}$. Repeating this process we conclude that $|v'(x)| \leq C_7$ and $|v(x)|\leq C_7 x$ for $x\in (0,1)$. This implies that 
\[
 q v \in L^1(0,1)\quad \Longrightarrow \quad v''(x) \in L^1(0,1). 
\]
A careful review of the argument above reveals that the constants $C_j$'s and all relevant norms in the definition of $D(\mathbf{A})$ are dominated by $\|v\|_{\mathcal{H}_{\mathbf{A}}^2}$ up to a constant $C=C(q)$. An integration by parts also shows that 
\[
 \int_0^\infty \left(|v'(x)|^2 + q(x) |v(x)|^2\right) = \langle \mathbf{A} v, v \rangle \lesssim \|v\|_{\mathcal{H}_{\mathbf A}^2}^2. 
\]
Thus $\|v'(x)\|_{L^2(\Rm^+)} \lesssim \|v\|_{\mathcal{H}_{\mathbf A}^2}$. The upper bound of $L^\infty$ norms of $v(x)$ and $v'(x)$ then follows from the Sobolev embedding. 

\paragraph{Case 2} If $\lambda\geq 3/4$, then we let $\gamma$ and $1-\gamma$ be two roots of the equation $z(z-1)= \lambda$, where $\gamma \geq 3/2$. Clearly all solutions to the differential equation $-u''(x) + \lambda x^{-2} u = 0$ are given by $C_1 x^\gamma + C_2 x^{1-\gamma}$. We then define the function  
\begin{align*}
 w(x) = x^\gamma \int_x^1 y^{1-\gamma} \cdot \frac{-q_0(y) v(y) + (\mathbf{A}^\ast v)(y)}{2\gamma-1} {\rm d} y + x^{1-\gamma} \int_0^x y^\gamma \cdot \frac{-q_0(y) v(y) + (\mathbf{A}^\ast v)(y)}{2\gamma-1} {\rm d} y. 
\end{align*}
A straight-forward calculation shows that $w$ solves the following differential equation in $(0,1)$ 
\[
 -w''(x) + \frac{\lambda}{x^2} w =  -q_0(x) v(x) + (\mathbf{A}^\ast v)(x) = - v''(x) + \frac{\lambda}{x^2} v.
\]
It immediately follows that 
\[
 v(x) = w(x) + C_1 x^\gamma + C_2 x^{1-\gamma}. 
\]
We recall the upper bound of $|q_0(x)| \lesssim |x|^{-\kappa}$ and apply the Cauchy-Schwarz to deduce  
\begin{align}
 |w(x)| &\lesssim_q x^\gamma \|y^{1-\gamma}\|_{L^2(x,1)} \|\mathbf{A}^\ast v\|_{L^2} + x^{1-\gamma} \|y^\gamma\|_{L^2(0,x)} \|\mathbf{A}^\ast v\|_{L^2} \nonumber \\
 & \qquad + x^\gamma \int_x^1 y^{1-\gamma-\kappa} |v(y)| {\rm d} y + x^{1-\gamma} \int_0^x y^{\gamma-\kappa} |v(y)| {\rm d} y \label{upper bound of w ap}\\
 & \lesssim_q J(x) \|\mathbf{A}^\ast v\|_{L^2}  + x^{3/2-\kappa} \|v\|_{L^2}. \label{upper bound of w ap1}
\end{align}
Here $J(x)$ is defined by 
\[
 J(x) = \left\{\begin{array}{ll} x^{3/2}, & \lambda>3/4; \\ x^{3/2}\left(|\ln x|^{1/2}+1\right), & \lambda = 3/4. \end{array}\right.
\]
 This implies that $w \in L^2(0,1)$. Since $v\in L^2(\Rm^+)$, it follows that $C_2 = 0$ thus 
\[
 v(x) = w(x) + C_1 x^\gamma. 
\]
By considering the $L^2$ norm of both sides, we also have $|C_1|\lesssim_q \|v\|_{L^2} + \|\mathbf{A}^\ast v\|_{L^2}$. As a result, we have 
\[
 |v(x)| \lesssim_q x^{3/2-\kappa} \left(\|v\|_{L^2} + \|\mathbf{A}^\ast v\|_{L^2}\right). 
\]
Inserting this into \eqref{upper bound of w ap}, we obtain 
\[
 |w(x)| \lesssim_q J(x) \|\mathbf{A}^\ast v\|_{L^2} + x^{7/2-2\kappa} \left(\|v\|_{L^2} + \|\mathbf{A}^\ast v\|_{L^2}\right). 
\]
Therefore we gain a little more regularity of $v$ near zero
\[
 |v(x)| \lesssim_q x^{7/2-2\kappa} \left(\|v\|_{L^2} + \|\mathbf{A}^\ast v\|_{L^2}\right). 
\] 
Repeating this process, we finally conclude that 
\[
 |v(x)| \lesssim_q J(x) \left(\|v\|_{L^2} + \|\mathbf{A}^\ast v\|_{L^2}\right). 
\]
This implies for any $p \in [1,2)$ that 
\[
 q v \in L^p(0,1) \quad \Longrightarrow \quad v'' \in L^p(0,1). 
\]
Thus $v'(0)$ and $v(0)$, namely the limits of $v(x)$ and $v'(x)$ at zero, exist. The inequality $|v(x)|\lesssim J(x)$ then yields $v(0)=v'(0)=0$. The remaining part of argument is similar to the case $\lambda =0$.

\section*{Acknowledgement}
Ruipeng Shen is financially supported by National Natural Science Foundation of China Project 12471230.


\begin{thebibliography}{99}
\bibitem{handbook} M. Abramowitz and I. A. Stegun(Editors). {``Handbook of Mathematical Functions, with Formulas, Graphs and Mathematical Tables.''} Dover Publications, Inc., New York, 1992.
\bibitem{abstractSchr} J. A. Barcel\'{o}, A. Ruiz and L. Vega. {``Some dispersive estimates for Schr\"{o}dinger equations with repulsive potentials.''} \textit{Journal of Functional Analysis} 236(2006): 1-24.
\bibitem{WPfastsmall} M. Beals and W. Strauss. {``$L^p$ estimates for the wave equation with a potential.''} \textit{Communications in Partial Differential Equations} 18(1993), no.7-8: 1365-1397. 
\bibitem{dispersive2} P. Brenner. {``On scattering and everywhere defined scattering operators for nonlinear Klein-Gordon equations.''} \textit{Journal of Differential Equations} 56(1985), no. 3: 310-344. 
\bibitem{BSchrWOAb} M. Beceanu. {``Structure of wave operators for a scaling-critical class of potentials.''} \textit{American Journal of Mathematics} 136, no. 2(2014): 255-308. 
\bibitem{BSSchrWO} M. Beceanu. and W. Schlag. {``Structure formulas for wave operators.''} \textit{American Journal of Mathematics} 142, no. 3(2020): 751-807.
\bibitem{WPKato3} T. A. Bui, X. T. Duong and Y. Hong. {``Dispersive and Strichartz estimates for the three-dimensional wave equation with a scaling-critical class of potentials.''} \textit{Journal of Functional Analysis} 271(2016), no. 8: 2215-2246. 
\bibitem{WISStrichartz1} N. Burq, F. Planchon, J. G. Stalker and A. S. Tahvildar-Zadeh. {``Strichartz estimates for the wave and Schr\"{o}dinger equations with the inverse-square potential.''} \textit{Journal of Functional Analysis} 203(2003), no. 2: 519-549.
\bibitem{WISStrichartz2} N. Burq, F. Planchon, J. G. Stalker and A. S. Tahvildar-Zadeh. {``Strichartz estimates for the wave and Schr\"{o}dinger equations with potentials of critical decay.''} \textit{Indiana University Mathematics Journal} 53(2004), no.6: 1665-1680.
 \bibitem{WPms} G. Chen. {``Multisolitons for the defocusing energy critical wave equation with potentials.''} \textit{Communications in Mathematical Physics} 364(2018), no.1: 45-82.
 \bibitem{StrichartzReversed} G. Chen. {``Wave equations with moving potentials.''} \textit{Communications in Mathematical Physics} 375(2020), no.2: 1503-1560.
 \bibitem{CKSchrWO} M. Christ. and A. Kiselev. {``Scattering and wave operators for one-dimensional Schr\"{o}dinger operators with slowly decaying nonsmooth potentials.''} \textit{Geometrical and Functional Analysis} 12, no. 6(2002): 1174-1234. 
 \bibitem{WPfast1d} O. Costin and M. Huang. {``Decay estimates for one-dimensional wave equations with inverse power potential.''} \textit{Transactions of the American Mathematical Society} 367(2015), no. 5: 3705-3732.
\bibitem{WPfast3d} S. Cuccagna. {``On the wave equation with a potential.''} \textit{Communications in Partial Differential Equations} 25(2000), no.7-8: 1549-1565. 
 \bibitem{WISsmalldata} W. Dai, D. Fang and C. Wang. {``Long-time existence for semi-linear wave equations with the inverse-square potential.''} \textit{Journal of Differential Equations} 309(2022): 98-141. 
 \bibitem{WMP1} P. D'ancona. {``Kato smoothing and Strichartz estimates for wave equations with magnetic potentials.''} \textit{Communications in Mathematical Physics} 335(2015): 1-16.
 \bibitem{WMLpSchr} P. D'ancona and L. Fanelli. {``$L^p$-boundedness of the wave operator for one-dimensional Schr\"{o}dinger operators.''} \textit{Communications in Mathematical Physics} 268, no. 2(2006): 415-438.
 \bibitem{WMP2} P. D'ancona and L. Fanelli. {``Decay estimates for the wave and Dirac equations with a magnetic potential.''} \textit{Communications in Pure and Applied Mathematics} 60(2007), no.3: 357-392.
 \bibitem{WPKato1} P. D'ancona and V. Pierfelice. {``On the wave equation with a large rough potential.''} \textit{Journal of Functional Analysis} 227(2005): 30-77.
 \bibitem{WPDen} S. A. Denisov. {``Wave equation with slowly decaying potential: asymptotics of solution and wave operator.''} \textit{Mathematical Modelling of Natural Phenomena} 5, no. 4(2010): 122-149. 
 \bibitem{SchrInversePower} V. D. Dinh. {``On nonlinear Schr\"{o}dinger equations with repulsive inverse-power potentials.''} \textit{Acta Applicandae Mathematicae} 171(2021): paper no. 14. 
 \bibitem{WIS1d} R. Donninger and W. Schlag. {``Decay estimates for the one-dimensional wave equation with an inverse power potential.''} \textit{International Mathematical Research Notices} 2010, no. 22: 4276-4300.
 \bibitem{WIS1d2} R. Donninger and J. Krieger. {``A vector field method on the distorted Fourier side and decay for wave equation with potentials.''} \textit{Memoirs of the American Mathematical Society} 241(2016), no. 1142. 
\bibitem{WMP3} L. Fanelli, J. Zhang and J. Zheng. {``Dispersive estimates for 2D-wave equations with critical potentials.''} \textit{Advances in Mathematics} 400(2022), paper no. 108333. 
\bibitem{WISGV} V. Georgiev and N. Visciglia. {``Decay estimates for the wave equation with potential.''} \textit{Communications in Partial Differential Equations} 28(2003), no.7-8: 1325-1369.
\bibitem{dispersive1} J. Ginibre and G. Velo. {``Time decay of finite energy solutions of the nonlinear Klein-Gordon and Schr\"{o}dinger equations.''}  Annales de l'I.H.P. Physique th\'{e}orique 43 (1985), no. 4: 399-442. 
\bibitem{strichartz} J. Ginibre, and G. Velo. {``Generalized Strichartz inequality for the wave equation.''} \textit{Journal of Functional Analysis} 133(1995): 50-68.
 \bibitem{WPfast2d} W. Green. {``Time decay estimates for the wave equation with potential in dimension two.''} \textit{Journal of Differential Equations} 257(2014): 868-919. 
 \bibitem{Goldberg1} M. Goldberg. {``Dispersive estimates for the three-dimensional Schr\"{o}dinger equation with rough potentials.''} \textit{American Journal of Mathematics} 128(2006), no. 3: 731-750.
 \bibitem{Goldberg2} M. Goldberg. {``Dispersive bounds for the three-dimensional Schr\"{o}dinger equation with almost critical potentials.''} \textit{Geometrical and Functional Analysis} 16(2006), no. 3: 517-536. 
 \bibitem{Goldberg3} M. Goldberg. {``Strichartz estimates for the Schr\"{o}dinger equation with time periodic $L^{n/2}$ potentials.''} \textit{Journal of Functional Analysis} 256(2009), no. 3: 718-746. 
\bibitem{LHSchrWO} L. H\"{o}rmander. {``The existence of wave operators in scattering theory.''} \textit{Mathematische Zeitschrift} 146, no. 1(1976): 69-91.
\bibitem{zeros} Y. Ikebe. {``The zeros of regular Coulomb wave functions and their derivatives.''} \textit{Mathematics of Computation} 29, no.131(1975): 878-887. 
 \bibitem{JKSchrOper} A. Jenson, and T. Kato. {``Spectral properties of Schr\"{o}dinger operators and time-decay of the wave functions.''} \textit{Duke Mathematical Journal} 46, no. 3(1979): 583-611. 
 \bibitem{JSchrOper} A. Jenson. {``Spectral properties of Schr\"{o}dinger operators and time decay of the wave functions results in $L^2(\Rm^m)$, $m\geq 5$.''} \textit{Duke Mathematical Journal} 47, no. 1(1980): 57-80.
 \bibitem{JSchrOperd4} A. Jenson. {``Spectral properties of Schr\"{o}dinger operators and time-decay of the wave functions. Results in $L^2(\Rm^4)$.''} \textit{Journal of Mathematical Analysis and Applications} 101, no. 2(1984): 397-422. 
 \bibitem{WPcm1} H. Jia, B. Liu and G. Xu. {``Long time dynamics of defocusing energy critical $3+1$ dimensional wave equation with potential in the radial case.''} \textit{Communications in Mathematical Physics} 339(2015): 353-384.
 \bibitem{WPcm2} H. Jia, B. Liu, W. Schlag and G. Xu. {``Generic and non-generic behavior of solutions to defocusing energy critical wave equation with potential in the radial case.''} \textit{International Mathematical Research Notices} 2017, no.19: 5977-6035.
 \bibitem{WPcm3} H. Jia, B. Liu, W. Schlag and G. Xu. {``Global center stable manifold for the defocusing energy critical wave equation with potential.''} \textit{American Journal of Mathematics} 142(2020), no.5: 1497-1557. 
 \bibitem{JSSdecaySchr} J.-L. Journ\'{e}, A. Soffer, and C. D. Sogge. {``Decay estimates for Schr\"{o}dinger operators.''} \textit{Communications in Pure and Applied Mathematics} 44, no. 5(1991): 573-604.
\bibitem{km} C. E. Kenig, and F. Merle. {``Nondispersive radial solutions to energy supercritical non-linear wave equations, with applications.''} \textit{American Journal of Mathematics} 133, No 4(2011): 1029-1065.
\bibitem{KS} R. Killip and B. Simon. {``Sum rules and spectral measures of Schr\"{o}dinger operators with $L^2$ potential.''} \textit{Annals of Mathematics} 170(2009): 739-782.
\bibitem{selfadjointbook} B. M. Levitan and I. S. Sargsjan. {``Introduction to spectral theory: Selfadjoint ordinary differential operators.''} American Mathematical Society, Providence, Rhode Island, 1975.
 \bibitem{WIScritical} C. Miao, J. Murphy and J. Zheng. {``The energy-critical nonlinear wave equation with an inverse-square potential.''} \textit{Annales de l'Institut Henri Poincar\'{e} C. Analyse Non Lin\'{e}aire} 37(2020), no.2: 417-456. 
 \bibitem{wavebook} C. Miao and R. Shen. {``Regularity and Scattering of Dispersive Wave Equations.''} \textit{De Gruyter Studies in Mathematics}, Volume 100, De Gruyter, Berlin-Boston 2025. 
 \bibitem{WISsubRadial} C. Miao, R. Shen and T. Zhao. {``Scattering theory for the subcritical wave equation with inverse square potential.''} \textit{Selecta Mathematica. New Series} 29(2023), no.3: paper no. 44. 
 \bibitem{SchrCoul} C. Miao, J. Zhang and J. Zheng. {``A nonlinear Schr\"{o}dinger equation with Coulomb potential.''} \textit{Acta Mathematica Scientia, Series B} 42(2022), no. 6: 2230-2256.   
 \bibitem{WISStrichartzRadial} C. Miao, J. Zhang and J. Zheng. {``Strichartz estimates for wave equation with inverse square potential.''} \textit{Communications in Contemporary Mathematics} 15(2013), no.6, 1350026. 
 \bibitem{SchrCritical} H. Mizutani. {``Remarks on endpoint Strichartz estimates for Schr\"{o}dinger equations with the critical inverse-square potential.''} \textit{Journal of Differential Equations} 263(2017), no. 7: 3832-3853.
 \bibitem{StrichartzSchrCoul} H. Mizutani. {``Strichartz estimates for Schr\"{o}dinger equations with slowly decaying potentials.''} \textit{Journal of Functional Analysis} 279(2020): 108789.  
 \bibitem{WISnonradial} H. Mizutani. {``Scattering theory in homogeneous Sobolev spaces for Schr\"{o}dinger and wave equation with rough potentials.''} \textit{Journal of Mathematical Physics} 61(2020), no.9: 091505. 
 \bibitem{Hsmooth1} H. Mizutani and X. Yao. {``Kato smoothing, Strichartz and uniform Sobolev estimates for fractional operators with sharp Hardy potentials.''} \textit{Communications in Mathematical Physics} 388(2021): 581-623.
\bibitem{sharmonics} C. M\"{u}ller. {``Spherical harmonics.''} Springer-Verlag, Berlin-New York, 1966. 
 \bibitem{enscatter1} K. Nakanishi. {``Unique global existence and asymptotic behaviour of solutions for wave equations with non-coercive critical nonlinearity.''} \textit{Communications in Partial Differential Equations} 24(1999): 185-221.
\bibitem{WPcompact} V. Petkov. {``Global Strichartz estimates for the wave equation with time-periodic potentials.''} \textit{Journal of Functional Analysis} 235(2006): 357-376.
\bibitem{WPKato2} V. Pierfelice. {``Decay estimate for the wave equation with a small potential.''} \textit{Nonlinear Differential Equations and Applications} 13(2007): 511-530. 
\bibitem{WISPST1} F. Planchon, J. G. Stalker and A. S. Tahvildar-Zadeh. {``$L^p$ estimates for the wave equation with the inverse-square potential.''} \textit{Discrete and Continuous Dynamical Systems. Series A} 9(2003), no.2: 427-442.
\bibitem{WISPST2} F. Planchon, J. G. Stalker and A. S. Tahvildar-Zadeh. {``Dispersive estimate for the wave equation with the inverse-square potential.''} \textit{Discrete and Continuous Dynamical Systems. Series A} 9(2003), no.6: 1387-1400.
\bibitem{Simon} M. Reed and B. Simon. {``Methods of Modern Mathematical Physics. I-IV''} Academic Press.   
\bibitem{TimedecaySchlag} I. Rodnianski and W. Schlag. {``Time decay for solutions of Schr\"{o}dinger equations with rough and time-dependent potentials.''} \textit{Inventiones Mathematicae} 155(2004), no. 3: 451-513. 
 \bibitem{ss2} J. Shatah, and M. Struwe. {``Well-posedness in the energy space for semilinear wave equations with critical growth''} \textit{International Mathematics Research Notices} 7(1994): 303-309.
 \bibitem{shenenergy} R. Shen. {``Energy distribution of radial solutions to energy subcritical wave equation with an application on scattering theory.''} \textit{Transactions of the American Mathematical Society} 374(2021): 3827-3857.
 \bibitem{shen3dnonradial} R. Shen {``Inward/outward Energy Theory of Non-radial Solutions to 3D Semi-linear Wave Equation''} \textit{Advances in Mathematics} 374(2020): 107384.
 \bibitem{shenhighradial} R. Shen. {``Long time behaviour of finite-energy radial solutions to energy subcritical wave equation in higher dimensions.''} \textit{Journal of Differential Equations} 368(2023): 26-57.
 \bibitem{Taylorbook} M. E. Taylor. {``Partial differential equations II. Qualitative studies of linear equations.. Second Edition.''} Springer, New York, 2011. 
\bibitem{WPfasthigher} G. Vodev. {``Dispersive estimates of solutions to the wave equation with a potential in dimensions $n\geq 4$.''} \textit{Communications in Partial Differential Equations} 31(2006), no. 10-12: 1709-1733.
\bibitem{WWkpSchr} R. Weder. {``The $W^{k,p}$-continuity of the Schr\"{o}dinger wave operators on the line.''} \textit{Communications in Mathematical Physics} 208, no. 2(1999): 507-520. 
\bibitem{YWkpSchr} K. Yajima. {``The $W^{k,p}$ continuity of wave operators for Schr\"{o}dinger operators.''} \textit{Journal of the Mathematical Society of Japan} 47, no. 3(1995): 551-581.
\bibitem{YLpSchr} K. Yajima. {``$L^p$-boundedness of wave operators for two-dimensional Schr\"{o}dinger operators.''} \textit{Communications in Mathematical Physics} 208, no. 1(1999): 125-152. 
\bibitem{SchrInverseSquare} J. Zhang and J. Zheng. {``Scattering theory for nonlinear Schr\"{o}dinger equations with inverse-square potential.''} \textit{Journal of Functional Analysis} 267(2014), no. 8: 2907-2932.
\end{thebibliography}
\end{document}